\documentclass{elsart}
\pdfoutput=1

\usepackage{verbatim,amsmath,amsfonts,subfigure}
\usepackage[pdftex]{graphicx}
\newcommand{\mlabel}[1]{\label{#1}}

\begin{document}

\sloppy

\def\R      {\mathbb R}
\def\C      {\mathbb C}
\def\ones   {\mathbf{1}}
\def\mathi      {\mathrm{i}}
\def\e      {\mathrm{e}}

\newtheorem{thrm}{Theorem}[section]
\newtheorem{crllr}{Corollary}[section]
\newtheorem{lmm}{Lemma}[section]
\newtheorem{definition}{Definition}[section]

\begin{frontmatter}



\title{Total--Variation--Diminishing Implicit--Explicit Runge--Kutta Methods
for the Simulation of Double-Diffusive Convection in Astrophysics}

\author[univie]{Friedrich Kupka\corauthref{cor1}\thanksref{support1}}
\ead{Friedrich.Kupka@univie.ac.at}
\ead[url]{http://www.mpa-garching.mpg.de/\~{}fk/}
\thanks[support1]{Supported by the Austrian Science Fund (FWF), project P21742-N16}

\author[univie]{Natalie Happenhofer\thanksref{support2}}
\ead{natalie.happenhofer@univie.ac.at}
\thanks[support2]{Supported by the Austrian Science Fund (FWF), project P20973}

\author[nav]{Inmaculada Higueras\thanksref{support3}}
\ead{higueras@unavarra.es}
\thanks[support3]{Supported by the Ministerio de Ciencia e Innovaci\'on, project  MTM2008-00785}

\author[tu]{Othmar Koch\thanksref{support1}}
\ead{othmar@othmar-koch.org}
\ead[url]{http://www.othmar-koch.org}

\corauth[cor1]{Corresponding Author}

\address[univie]{University of Vienna, Faculty of Mathematics,
Nordbergstra{\ss}e 15, A-1090 Wien, Austria}

\address[nav]{Universidad P\'ublica de Navarra, Departamento de Ingenier{\'i}a Matem\'atica e Inform\'atica, Campus de Arrosadia,
31006 Pamplona, Spain}

\address[tu]{Vienna University of Technology, Institute for Analysis and Scientific Computing,
A-1040 Wien, Austria. Formerly at University of Vienna, Faculty of Mathematics,
Nordbergstra{\ss}e 15, A-1090 Wien, Austria}

\begin{abstract}
We put forward the use of \emph{total--variation--diminishing} (or more generally,
\emph{strong stability preserving}) implicit--explicit Runge--Kutta methods
for the time integration of the equations of motion associated with the semiconvection problem
in the simulation of stellar convection. The fully compressible Navier--Stokes equation, augmented 
by continuity and total energy equations, and an equation of state describing the relation between 
the thermodynamic quantities, is semi-discretized in space by \emph{essentially non-oscillatory schemes}
and dissipative finite difference methods. It is subsequently integrated in time by Runge--Kutta methods
which are constructed such as to preserve the total variation diminishing (or strong stability) property 
satisfied by the spatial discretization coupled with the forward Euler method. We analyse the stability, accuracy 
and dissipativity of the time integrators and demonstrate that the most successful methods yield 
a substantial gain in computational efficiency as compared to classical explicit Runge--Kutta methods.
\end{abstract}

\begin{keyword}
hydrodynamics \sep stellar convection and pulsation
\sep double-diffusive convection \sep
numerical methods \sep total--variation--diminishing
\sep strong stability preserving \sep TVD \sep SSP
\MSC 65M06 \sep 65M08 \sep 65M20 \sep 65L05 \sep 76F65
\PACS 97.10.Cv \sep 97.10.Sj \sep 47.11.-j \sep 02.70.-c
\end{keyword}
\end{frontmatter}

\section*{Introduction}
\label{sec:intro}

Numerical hydrodynamical simulations are a common tool in
astrophysical research. Just as some of their counterparts in the
atmospheric sciences and in oceanography, astrophysical fluid flows
are characterized by a vast range of timescales which are present
in the solutions of the dynamical equations governing the temporal
evolution of such flows \cite{hujran01}. Large relative changes of
the solutions typically occur on the hydrodynamical timescale
$\tau_{\mathrm fluid}=\Delta x / |\mathbf{u}|$. Here, $\mathbf{u}$ is the
local flow velocity and $\Delta x$ is the local spatial resolution of the
simulation, which coincides with the grid size obtained from spatial
discretization of the governing partial differential equations.
However, some of the physically important processes can also operate
on much shorter timescales than $\tau_{\mathrm fluid}$. Examples
include radiative transfer, sound waves, magnetohydrodynamic
processes, and chemical or nuclear reactions (see \cite{hujran01}
for example, and references therein).

In stellar astrophysics the two most important among those timescales are that
of radiative energy exchange at the scale of a grid cell, $\tau_{\mathrm rad}$,
and the time $\tau_{\mathrm sound}$ a sound wave needs to cross such
a cell.

As long as sound waves are energetically or dynamically unimportant, a numerical
simulation can be advanced with much larger time steps by using semi-implicit time
integration methods, for example, by a fractional step approach (see \cite{ferziger02}
for a general introduction).

Similarly, if radiative transfer has the numerical characteristics of a stiff problem,
as is the case, for instance, for numerical simulations of the surface layers of
A-type stars \cite{kupkaetal09}, implicit time integration appears desirable as well.
Another important example where radiative diffusion can limit the time-step
is the numerical simulation of double-diffusive processes in stellar interiors.

Semiconvection is the most important special case of double-diffusive convection
in astrophysics. Models of stellar structure and evolution predict settings where
the heavier product of nuclear fusion provides stability to a zone which otherwise
would be unstable to convective overturning, because temperature sufficiently
rapidly decreases against the direction of gravity.
Such a zone would become convective if its composition were mixed. The question
whether such a zone should be treated as if it were mixed or not is referred to as the
\emph{semiconvection problem} (see \cite{kipweig94}, \cite{schwaharm58}, and
Chap.~13.3 and 13-A in \cite{weissetal04}, for example). A thorough physical analysis of the
semiconvection problem based on numerical simulations in two spatial dimensions for
a parameter set relevant to stellar astrophysics is given in \cite{zaussispruit11}. Further
discussions and reviews on this topic can be found, for instance, in
\cite{canuto09,canuto11,huplin79,spruit92,turner85}. For this problem, long total
integration times are required even when the microphysics is idealised, whereas
the time integration step is governed by $\tau_{\mathrm rad}$ and $\tau_{\mathrm sound}$.

For reasons outlined above, in this paper we discuss the advantages
of implicit--explicit (IMEX) Runge--Kutta methods for simulations of stellar
convection and diffusion in the parameter regime commonly associated with
semiconvection as discussed in \cite{zaussispruit11}. These methods
treat only part of the right-hand sides implicitly, where the
resulting (generally nonlinear) equations can be solved by means of a
generalized Poisson problem. It turns out that the \emph{total--variation--diminishing
(TVD)} property is essential for a numerical time integrator to be successful 
in simulations of the problems in our focus: to suppress spurious oscillations in the
spatial discretization (which in this paper we realize for the hyperbolic terms by
essentially non-oscillatory schemes and by dissipative centered
finite difference schemes for the parabolic terms), this property
has been demonstrated to be necessary for a stable integration in
\cite{gottliebetal09}. 

The TVD property is more generally referred to as \emph{strong
stability preserving (SSP)} or \emph{monotonicity} when norms other
than the total variation norm or even sublinear functionals are
considered. 
When the space discretization has the property that the functional
of the discrete spatial profile is decreased in the course of numerical
time propagation by the forward Euler method for a time-step
$\Delta t_{\mathrm FE}$, then an SSP method preserves this property
under a step-size restriction of the form 
$\Delta t\leq C\Delta t_{\mathrm FE}$ with $C>0$. Since the term
\emph{total--variation--diminishing} is more commonly used in the
context of astrophysical simulations, where the total--variation--seminorm
is the functional of interest, we mostly use these terms here synonymously.
We expect the SSP (or TVD) IMEX methods to be useful also for
other astrophysical problems where  a high radiative (conductive)
diffusivity of internal energy (temperature) restricts the time-step
of hydrodynamical simulations, such as simulations of stellar
surface convection with steep temperature gradients or at high
resolution.

The outline of the paper is as follows. First, we introduce 
SSP IMEX Runge--Kutta schemes and survey the related literature in Section~\ref{imex}.

Next, in Section~\ref{sec:antares} we specify the general set of equations to be solved
in numerical simulations of semiconvection and related flows and
describe the general solution techniques implemented for this class of problems
in the ANTARES code \cite{Muthsam2010} which we use for the numerical
examples discussed further below. Subsequently, we discuss how this framework
has to be modified when solving the dynamical equations with the IMEX
approach.

In Section~\ref{imex-candidates} we analyse several SSP IMEX methods
from the literature with respect to their radius of absolute
monotonicity, stability and dissipativity. We show that the methods
yield a significant advantage over classical explicit Runge--Kutta
schemes with respect to both efficiency and accuracy and we also
suggest a modification for one of the methods, which turns out to
improve its efficiency.

We then present numerical simulations of semiconvection in Section~\ref{num} to 
demonstrate the efficiency of the SSP IMEX methods as compared to the classical,
explicit SSP Runge--Kutta time integrators and some non-SSP IMEX methods from
the literature by giving numerical examples for a single layer in a physical scenario
similar to that one studied in \cite{zaussispruit11}. We conclude this paper by a summary
of the main properties of the IMEX methods, suggesting reasons for preferring particular
methods and providing an outlook on interesting applications, which appear
especially suited for this numerical approach.

\section{Implicit--Explicit Runge--Kutta Methods
for Semiconvection}\mlabel{imex}\mlabel{sec:semi}


To introduce our numerical methods in an abstract setting, we consider the ODE initial value problem
\begin{equation}\mlabel{de1}
\dot y(t) = F(y(t)) + G(y(t)),\qquad y(0)=y_0,
\end{equation}
where we assume that the vector fields $F$ and $G$ have different stiffness properties.
For this type of problems, \emph{partitioned Runge--Kutta schemes} \cite{haireretal87},
also called \emph{additive Runge--Kutta schemes}, are popular. Methods of this kind
use different Runge--Kutta formulae for the treatment of the two vector
fields. We will see that the spatial semi-discretization of (\ref{mnseimex}) below and the
associated boundary conditions give rise to this kind of system.

An $s$-stage partitioned Runge--Kutta method characterized by coefficient matrices
$A=(a_{i,j})$ and $\tilde A = (\tilde a_{i,j})$ defines one step
$y_\mathrm{old} \to y_\mathrm{new}$ by
\begin{eqnarray}
&& y_i = y_\mathrm{old} + \Delta t \sum_{j=1}^s a_{i,j} F(y_j) +
                          \Delta t \sum_{j=1}^s \tilde{a}_{i,j} G(y_j),\quad i=1,\dots,s,\mlabel{rk1}\\
&& y_\mathrm{new} = y_\mathrm{old} + \Delta t \sum_{j=1}^s b_{j} F(y_j) +
                          \Delta t \sum_{j=1}^s \tilde{b}_{j} G(y_j).\mlabel{rk2}
\end{eqnarray}
If ${a}_{i,j}=0$ for $j\geq i$, the method is referred to as an \emph{implicit--explicit (IMEX)\/}
method.

These methods have first been investigated with respect to the SSP
or TVD property in the context of hyperbolic systems with relaxation,
where $G = \frac1\varepsilon \hat G,\ \varepsilon \ll 1,$ in
\cite{parrus05}. Ibidem, the common specification for strong
stability preserving IMEX methods is introduced. An IMEX method is referred to as
`SSP$k(s,\sigma,p)$' when it has the following properties: $k$ is
the order of the method in the stiff limit $\varepsilon\to0$, which
is characterized by the coefficients for the explicit part. The
latter must necessarily be SSP and is referred to as the
\emph{asymptotically SSP scheme}. $s$ is the number of stages in the
implicit scheme and $\sigma$ the number of stages in the explicit
scheme. $p$ is the global order of the resulting combined method. It
is essential to observe that if the implicit scheme characterized by
$\tilde A=(\tilde a_{i,j})$ is a \emph{diagonally implicit
Runge--Kutta (DIRK)} method, then the explicit part is evaluated
only once in each stage, providing the desired computational
advantage \cite{parrus05}.

The analysis in \cite{parrus05} is valid only for $\varepsilon \ll 1$ \cite{higueras06}.
However, several useful examples of strong stability preserving IMEX Runge--Kutta methods are given,
see Section~\ref{imex-candidates}.

\cite{higueras06} develops a comprehensive theory of strong stability preserving
additive Runge--Kutta schemes which extends the concepts for standard Runge--Kutta
methods in a natural way:

Let $\tau, \ \tilde \tau$ be the step-size restrictions for
monotonicity of the explicit Euler method for the vector fields $F$ and
$G$, respectively. We define the \emph{region of absolute monotonicity}
\begin{equation}\mlabel{stabdef}
\mathcal{R}(A,\tilde A) = \{ (r,\tilde r) \in \R^2: (A,\tilde A) \mbox{ is absolutely
monotonous on } [-r,0] \times [-\tilde r,0]\},
\end{equation}
where the absolute monotonicity at a point $(r_0,\tilde r_0)$ is characterized
by algebraic relations for the matrices $A,\ \tilde A$. The boundary in the first
quadrant, $\partial \mathcal{R}(A,\tilde A) \cap \{(r,\tilde r): r,\tilde r \geq0\},$
is denoted as the \emph{curve of absolute monotonicity}.
The significance of the region $\mathcal{R}(A,\tilde A)$ is expressed in the following theorem \cite{higueras06}:
\begin{thrm}\mlabel{am1}
Let $(A,\tilde A)$ be absolutely monotonous at $(-r,-\tilde r)$ with step-size
coefficients $\tau,\ \tilde \tau$. Then for $h\leq
\min\left\{r \tau, \tilde r \tilde \tau \right\}$, diminishing
of the norm holds,
\begin{eqnarray*}
&& \|y_i\|\leq \|y_\mathrm{old}\|,\quad i=1,\dots,s, \qquad \|y_\mathrm{new}\| \leq \|y_\mathrm{old}\|.
\end{eqnarray*}
\end{thrm}
\cite{higueras09} gives order barriers for strong stability preserving additive Runge--Kutta methods
similarly to \cite{kraaijevanger91}. The order of an additive Runge--Kutta method $(A,\tilde A)$ is
bounded by the orders of $A$ and $\tilde A$, respectively. This implies for IMEX methods the
order barrier $p\leq 4$ \cite{kraaijevanger91}. Moreover, \cite{higueras09} gives a
simple algebraic criterion for a nontrivial region of absolute monotonicity in terms
of incidence matrices of $A,\ \tilde A$. Some examples of strong stability preserving
IMEX Runge--Kutta methods analysed in Section~\ref{imex-candidates}
can be seen in \cite{higueras09,kraaijevanger91}.

Some other relevant issues can be studied for IMEX methods. In this paper we focus
on stability regions, error constants and dissipativity analysis. We close this
section with a brief description of these concepts.

The stability region of IMEX Runge--Kutta methods is defined in
\cite{ascheretal97,ascheretal95} via the test equation of the form
(\ref{de1}), where
\begin{equation}\mlabel{testeq1}
F(u) = \mathi \beta u,\quad G(u) = \alpha u,\qquad \alpha \leq 0 < \beta.
\end{equation}
For this problem,
\begin{equation}\mlabel{stabfcn}
y_\mathrm{new} = R(z)y_\mathrm{old},\qquad z= \alpha \Delta t + \mathrm{i} \beta \Delta t,
\end{equation}
and the stability region is the part of the complex plane where $|R(z)|<1$.
We will perform the corresponding analysis of the stability function for the
methods considered in Section~\ref{imex-candidates}.

To determine error constants of the methods we have computed
the empirical convergence orders by solving the non-linear test problem
\begin{equation}\mlabel{eq:test1}
y'(t)=(1+\sin(y(t))) + (y^2(t)-\sin(y(t))),\qquad y(0)=0,
\end{equation}
with the known exact solution $y(t)=\tan(t)$. In this paper, the error constants are
determined from the errors at $t=1.3$.
Their size is vital for the comparison of the accuracy of the methods of
the same order and therefore the assessment of the work/precision relation.
Of course, this single example only gives a rough indication of the size
of the error constant and no rigorous estimate, but it seems sufficient
for our purpose of comparing the methods in our focus with respect to
accuracy, which we do in Sections~\ref{num} and \ref{conclusions} further 
below. The example was chosen such as to represent a nonlinear initial
value problem with known solution whose profile is arbitrarily unsmooth
as $t\to\frac\pi2$.

Finally, we study the dissipativity of the time integrators in
conjunction with suitable space discretizations. \cite{strikwerda04} gives a justification for considering only the
diffusion term in this context since the advection term becomes negligible asymptotically.
We will thus investigate the dissipativity of the implicit scheme specified by $\tilde A$.
To this end, we apply the spatial discretizations $L_{\Delta x}$ in our focus to the heat equation
$$ u_t = b u_{xx},$$
and associate for the spatial discretization $u_{j\pm k} \leftrightarrow \e^{\pm \mathi \theta k}$.
Thus, we compute the \textit{amplification factor}
$g(\theta,\mu) = R((L_{\Delta x} u)_j)$, with $\mu = b \frac{\Delta t}{(\Delta x)^2}$.
This represents the factor by which oscillations of frequency $\theta$ are amplified in each time step.
We pay particular attention to the case $\theta=\pi$ corresponding
to the mesh width. If $g=0$ or $|g|=1$ also for a smaller
$0<|\theta_0|<\pi$, then this value would represent the limit for a robust integration.
However, such a pathological behaviour is only conceivable for methods of higher order.

The spatial discretizations which were found to show a dissipative behaviour in
\cite{kochetal10a} are the second-order three-point scheme (the upper index refers to the time step)
\begin{equation}   \label{3pt}
u_{xx}(x_j,t_n) \approx \frac{u_{j+1}^n - 2 u_{j}^n + u_{j-1}^n}{(\Delta x)^2} =: (L_{\Delta x} u^n)_j,
\end{equation}
and the fourth-order stencil
\begin{equation}   \label{4th}
(L_{\Delta x} u^n)_j := \frac{-u_{j+2}^n + 16\,u_{j+1}^n - 30\, u_{j}^n + 16\,u_{j-1}^n - u_{j-2}^n}{12\, (\Delta x)^2}.
\end{equation}

These are the methods actually implemented in ANTARES, where (\ref{4th}) is the default (see also \cite{kochetal10a}).

\section{Solving the Hydrodynamical Equations with ANTARES}
\label{sec:antares}

In our model, the fundamental equation of motion is the \emph{fully compressible Navier--Stokes equation}
which describes momentum conservation:
\begin{equation}\mlabel{momeq1}
(\rho \mathbf{u})' + \nabla \cdot (\rho \mathbf{u} \otimes \mathbf{u} + p I) = \rho \mathbf{g} + \nabla \cdot \sigma.
\end{equation}
The state variables in the model equations generally depend on the spatial variables
$(x,y,z)$ and time $t$. In the simulations presented in Section~\ref{num} we solve problems in
two spatial variables, whence the variable $z$ will be dropped in the rest of the paper.
The (explicit) dependencies are stated in Table~\ref{tab:const}.
For simplicity, we omit the dependencies in the problem specification
(\ref{momeq1})--(\ref{rt1}).
The model is completed by the \emph{continuity equation}
\begin{equation}\mlabel{conteq1}
\rho' + \nabla \cdot(\rho \mathbf{u}) =0,
\end{equation}
which ensures conservation of mass, and the \emph{total energy equation}
\begin{equation}\mlabel{toteeq1}
e'+\nabla \cdot ( \mathbf{u}(e+p) ) = \rho (\mathbf{g} \cdot \mathbf{u}) +
\nabla \cdot (\mathbf{u} \cdot \sigma) + Q_\mathrm{rad},
\end{equation}
which describes conservation of the latter. In the case of a two-component fluid,
the system is augmented by the concentration equation of the second species,
\begin{equation}\mlabel{conceq}
(c \rho)' + \nabla \cdot (c \rho \mathbf{u}) = \nabla \cdot (\rho \kappa_c \nabla c).
\end{equation}

The variables and parameters
which appear in the model formulation are collected in Table~\ref{tab:const}.

\begin{table}[h!]
\begin{center}
\begin{tabular}{|l|l|}
\hline
$\rho=\rho(x,y,z,t)$                        & gas density \\
$c = c(x,y,z,t)$                            & concentration of second species \\
$\mathbf{u} =\mathbf{u}(x,y,z,t)=(u,v,w)^T$ & flow velocity \\
$\rho \mathbf{u}$                           & momentum density \\
$\mathbf{u} \otimes \mathbf{u}$             & Kronecker product \\
$p=p(T,\rho)$                               & gas pressure \\
$\mathbf{g}=(g,0,0)^T$                      & gravitational acceleration \\
$\sigma=\sigma(x,y,z,t)$                        & viscous stress tensor for zero bulk viscosity \\
$\eta$                                       & dynamic viscosity (appears in the definition of $\sigma$)\\
$e = e(x,y,z,t) = e_\mathrm{int} + e_\mathrm{kin}$   & total energy density \\
$T=T(x,y,z,t)$                              & temperature \\
$Q_\mathrm{rad}=Q_\mathrm{rad}(x,y,z,t)$    & radiative source term \\
$c_\mathrm{p}=c_\mathrm{p}(T,\rho,c)$   & specific heat at constant pressure \\
$\chi_\nu=\chi_\nu(T,\rho,c)$                 & (specific) opacity at frequency $\nu$ \\
$K=K(T,\rho,c)$   & radiative (or thermal) conductivity \\
$\kappa_T=  K / (c_\mathrm{p} \rho)$                     & radiative (or thermal) diffusivity \\
$\kappa_c=\kappa_c(T,\rho,c) $                  & diffusion coefficient for species c \\

$I_\nu = I_\nu(\mathbf{r}),\
\mathbf{r}=\mathbf{r}(x,y,z)$               & specific intensity along the ray of direction $\mathbf{r}$\\
$S_\nu = S_\nu(x,y,z)$                      & source function \\
\hline
\end{tabular}
\caption{Variables and parameters in the equations (\ref{momeq1})--(\ref{rt1}).\mlabel{tab:const}}
\end{center}
\end{table}

In general, the radiative source term $Q_\mathrm{rad}$ is determined as the
stationary limit of the \emph{radiative transfer equation}
\begin{equation}\mlabel{rt1}
\mathbf{r} \cdot \nabla I_\nu = \rho \chi_\nu (S_\nu-I_\nu),
\end{equation}
which is solved for all ray directions $\mathbf{r}$ and for all frequencies $\nu$,
resulting in the specific intensity $I_\nu(\mathbf{r})$, for details see \cite{weissetal04}.
$S_\nu$ here denotes the \emph{source function}.

The equations of hydrodynamics (\ref{momeq1}), (\ref{conteq1}) and (\ref{toteeq1})
are closed by the equation of state which describes the relation between the
thermodynamic quantities. For the particular choice, see \cite{Muthsam2010}.

For the initial condition, a slightly perturbed static model stellar atmosphere or
stellar envelope
is used which is equipped with a small seed velocity field or density
perturbation to start dynamics away from equilibrium.

In the framework discussed below, boundary conditions are based on the assumption that all quantities
are periodic in both horizontal directions. Moreover, for the hydrodynamical
equations, ``closed'' (Dirichlet) boundary conditions at the upper
and lower boundary of the computational domain are used.
A recent development is to replace these by ``open'' (Robin) boundary
conditions. These allow inflow and outflow of fluid along the vertical
direction which is defined by the direction of $\mathbf{g}$. For the radiative transfer equation (\ref{rt1}),
incoming radiation at the boundary of the computational domain must
be specified.
Since double-diffusive convection in stars takes place in regions which
are optically thick, the quantity $Q_\mathrm{rad}$ can accurately be
obtained by means of the diffusion approximation for radiative transfer,
$Q_\mathrm{rad} = \nabla \cdot F_\mathrm{rad} = \nabla \cdot (K \nabla T)$.
In this case, further knowledge about the intensity $I_\nu$ is not necessary.

The ANTARES code \cite{Muthsam2010} solves this system of equations
numerically in either one, two, or three spatial dimensions on a rectangular
grid (spherical coordinates with a logarithmically rectangular grid are also possible,
i.e., the grid may be locally rectangular with logarithmic grading in the radial component).

For the spatial discretization, ANTARES allows the definition of several
grids which can be nested inside each other to improve
resolution in regions of interest. At the moment, ANTARES provides up to
three levels of nested grids.
For the hyperbolic terms, discretizations of
ENO (\emph{essentially non-oscillatory} \cite{shu97}) type are implemented.
These comprise classical ENO methods, WENO (\emph{weighted essentially non-oscillatory})
methods \cite{shu97} (optionally in conjunction with
Marquina flux splitting \cite{marquina96}) and CNO (\emph{convex non-oscillatory}) schemes
\cite{liuosher98}. Each of the methods uses adaptive stencils which are chosen
such as to avoid spurious oscillations in the computed solution.
The spatial derivatives are calculated for each direction separately.

The parabolic terms are discretized by dissipative finite difference
schemes \cite{kochetal10a} of fourth order. The \emph{radiative heating rate}
is determined by the \emph{short characteristics method}, or by means of the
diffusion approximation for radiative transfer, where appropriate. For the time integration,
\emph{total variation diminishing} Runge--Kutta methods \cite{shu88,shu88a} are employed.

ANTARES implements two different
parallelization concepts. For architectures with distributed
memory, domain decomposition is used and realized by an MPI implementation.
In this approach each grid is split along the horizontal direction(s) and optionally,
also along vertical ones, into subdomains. The memory
required to store the computational variables for each
subdomain is provided by the resources available to the
CPU core performing the computations
necessary for that subdomain. In this way, each CPU core is
mapped to a specific geometrical volume.
However, since some supercomputers offer only a limited amount
of memory per CPU core and because the domain
decomposition approach is not very efficient on small grids, ANTARES
offers a second type of parallelization which can be used along with
or independently of the former. It is based on a shared
memory concept for each subdomain and is implemented
through OpenMP directives. Therefore, the most time consuming
operations which can also be performed independently of each
other are identified and computed in parallel. This approach scales 
only to a moderate number of CPU cores, but allows improvement 
of the scaling and the computational speed of the domain decomposition
based parallelization for a larger number of problems and
for a greater variety of computer architectures.

In the following, the dynamical evolution of the fluid is described by
the multispecies Navier--Stokes equations presented above. Additionally,
dimensionless quantities such as the Prandtl number $\mathrm{Pr} = c_\mathrm{p} \eta / K$,
the Lewis number $\mathrm{Le} = c_\mathrm{p} \rho \kappa_c / K$,
the Rayleigh number $\mathrm{Ra}$ and the stability parameter $R_{\rho}$ are
defined to determine the diffusivities $\kappa_T$, $\kappa_c$ and
the viscosity $\eta$. The former quantities arise in the definition
of the starting model but do not appear in the evolution equations
(\ref{mnseimex}) below.
Since we solve the dynamical equations for a compressible flow, we specify
the vertical extent of the simulation domain in multiples of the pressure scale height
$H_p = P / (\rho g)$. For the simulations presented in Section~\ref{num} the domain
always covers $1 H_p$. In the physical model, intermolecular forces are neglected, so the fluid is
assumed to be an ideal gas. The radiative source term $Q_\mathrm{rad}$ is modelled
using the diffusion approximation with a heat conductivity $K$ which is constant in time
and otherwise only a function of the vertical coordinate \cite{muthsamconvzones,muthetal95}.
We use a variant of this setting here, where not only $\mathrm{Pr}$, $\mathrm{Le}$, and 
$c_\mathrm{p}$, but also $K$, $\kappa_c$, and $\eta/\rho$ take constant values
\cite{zaussispruit11}. This setup simplifies studies of the basic physics while it is still useful
for extrapolations to astrophysically relevant cases (cf.\ \cite{zaussispruit11}).

For the model problem, the multispecies Navier--Stokes equations can be recast as

\begin{equation}
\underbrace{
\frac{d}{dt}
\left(\begin{array}[4]{c}
\rho \\
\rho c \\
\rho \mathbf{u} \\
e \\
\end{array}\right)
}_{\dot y(t)}
=
\underbrace{-
\nabla \cdot
\left(\begin{array}[4]{c}
\rho \mathbf{u} \\
\rho c \mathbf{u} \\
\rho \mathbf{u} \otimes \mathbf{u} + P - \sigma\\
e \mathbf{u} + P \mathbf{u} - \mathbf{u} \cdot \sigma\\
\end{array}\right)
-
\left(\begin{array}[4]{c}
0\\
0\\
\rho g \\
\rho g \mathbf{u} \\
\end{array}\right)
}_{F(y(t))}
+
\underbrace{
\nabla \cdot
\left(\begin{array}[4]{c}
0\\
\rho \kappa_c \nabla c \\
0 \\
K \nabla T \\
\end{array}\right)
}_{G(y(t))}.  \label{mnseimex}
\end{equation}

In the context of the problem (\ref{mnseimex}),
the $i$\textsuperscript{th} implicit stage of an IMEX Runge--Kutta method is typically of the form

\begin{equation}
 y_i = y^* + \Delta t\, \tilde{a}_{i,i}\, G(y_i),
\end{equation}
where $y^*$ is known from previous stages. This is a consequence of $\tilde a_{i,j}=0$
for $j>i$.

This translates to

\begin{eqnarray}
&& \rho_i = \rho^*, \\
&&\left(\rho c \right)_i = (\rho c)^* + \Delta t\, \tilde{a}_{i,i} \nabla \cdot (\rho_i \kappa_c \nabla c_i), \label{impc}\\
&&(\rho \mathbf{u})_i = (\rho \mathbf{u})^*, \label{eq:impmom} \\
&&e_i = e^* + \Delta t\, \tilde{a}_{i,i} \nabla \cdot (K \nabla T_i). \label{eq:impe}
\end{eqnarray}

Rearranging (\ref{impc}) leads to

\begin{equation}
\frac{\rho^*}{\Delta t\, \tilde{a}_{i,i}} c_i - \nabla \cdot (\rho^* \kappa_c \nabla c_i) =
\frac{\rho^* c^*}{\Delta t\, \tilde{a}_{i,i}}. \label{eq:genellc}
\end{equation}

Obviously, this is a general elliptic equation for $c_i$ of the form

\begin{equation}
g(x,y) \varphi(x,y) - \nabla \cdot (h(x,y) \nabla \varphi(x,y)) = f(x,y). \label{genell}
\end{equation}

\vspace{1em}

Due to model assumptions, (\ref{eq:impe}) can also be transformed to resemble a general elliptic equation. We start by recalling

\begin{eqnarray}
e &=& e_\mathrm{int}+ e_\mathrm{kin} \\
  &=& e_\mathrm{int}+ \frac{1}{2} \rho |\mathbf{u}|^2.
\end{eqnarray}

Bearing in mind equation (\ref{eq:impmom}), equation (\ref{eq:impe}) reads
\begin{equation}
e_{\mathrm{int}\text{ }i} = e_{\mathrm{int}}^* + \Delta t\, \tilde{a}_{i,i} \nabla \cdot (K \nabla T_i).
\end{equation}
The equation of state for an ideal gas\footnote{Note that, as is common in astrophysics, the
gas constant $R_{\mathrm{gas}}$ is taken to be relative to the atomic mass unit such that 
the dimensionless mean molecular weight can be used in the equation of state instead of the molar mass.}
relates the temperature $T$ and the internal energy $e_\mathrm{int}$ via
\begin{equation}
e_\mathrm{int}= \frac{3}{2} \frac{T \rho R_{\mathrm{gas}}}{m}, \label{eq:eos}
\end{equation}
if we assume the ratio of the specific heats at constant pressure and volume to equal $5/3$.
Here, $m$ denotes the mean molecular weight of the compound.

So we have at stage $i$
\begin{equation}
e_{\mathrm{int}\text{ }i} = \frac{3}{2} \frac{T_i \rho^* R_{\mathrm{gas}}}{m_i}
\end{equation}

and therefore, we arrive at

\begin{equation}
\frac{3}{2} \frac{R_\mathrm{gas} \rho^*}{m_i\, \Delta t\, \tilde{a}_{i,i}} T_i - \nabla \cdot (K \nabla T_i)
= \frac{e_{\mathrm{int}}^*}{\Delta t\, \tilde{a}_{i,i}}.
\end{equation}

Since $m_i$ is evaluated using the mass fraction $c_i$, it is necessary to solve equation (\ref{eq:genellc}) first.

Thus, the solution of an implicit stage translates to the solution of the generalized Poisson equations for
the mass fraction $c$ and the temperature $T$. In the ANTARES framework, finite elements are used for
the discretization of (\ref{genell}). The resulting linear system is solved by the conjugate gradient method.
For parallel computations, the Schur complement algorithm is applied. A detailed description is given in \cite{thesishgs}.

The above procedure applies without modification to the case of the fully compressible
Navier--Stokes equation. However, for low Mach number flows a splitting approach is preferable
where the terms containing pressure are treated separately in a post-processing step, i.\,e.\ after
evaluation of all other terms for the computation of the velocity fields. The latter are obtained for the
current step from an additional generalized Poisson equation for the pressure. This procedure
is described in detail in \cite{ferziger02}. Consequently, the explicit stage is evaluated here using a fractional
step method \cite{Kwatra09} as implemented in \cite{thesisnh}.

\section{Strong Stability Preserving IMEX Schemes from the Literature}\mlabel{imex-candidates}

In this section we study different SSP IMEX methods from the literature focusing
on the topics described in Section~\ref{imex}. The results are summarized
in Figure~\ref{stabcompareimex} and Table~\ref{summary}.

All the strong stability preserving IMEX schemes listed in the following subsections have DIRK
(diagonally implicit Runge--Kutta) methods
as the implicit scheme. This structure ensures that the stages can be solved successively,
and the explicit part only has to be evaluated once in each stage.

\subsection*{An SSP2(2,2,2) Method}
\cite{parrus05} gives an IMEX SSP2(2,2,2) method with nontrivial region
of absolute monotonicity ($\gamma = 1-\frac1{\sqrt2}$):
\begin{equation}\mlabel{hig1}
\begin{array}{c|cc}
0       & 0       & 0       \\[1mm]
1       & 1       & 0       \\[1mm] \hline\\[-3mm]
A       & \frac12 & \frac12
\end{array}
\hspace*{1cm}
\begin{array}{c|ccc}
\gamma   & \gamma    & 0      \\[1mm]
1-\gamma & 1-2\gamma & \gamma \\[1mm] \hline\\[-3mm]
\tilde A & \frac12   & \frac12
\end{array}.
\end{equation}
The coefficients imply $\mathcal{R}(A)=1$, \ $\mathcal{R}(\tilde A)= 1
+ \sqrt2$, and
$$ \mathcal{R}(A,\tilde A) = \{(r,\tilde r): 0 \leq r \leq 1,\ 0 \leq \tilde r \leq \sqrt2(1-r)\},$$
see \cite{higueras06}.

The stability region is entirely located in the left half plane,
tangent to the imaginary axis and unbounded as $\Re(z)\to-\infty$.
Hence, the schemes are $A(\alpha)$--stable with $\alpha=\frac\pi2$,
but not $A$-stable \cite{haiwan91}. Moreover, $\lim_{\Re(z)\to -\infty}R(z) =0$.
A plot of the stability region is given in Figure~\ref{fig1} (left).
\begin{figure}[h]
\begin{center}
\includegraphics[width=4.7cm,angle=270]{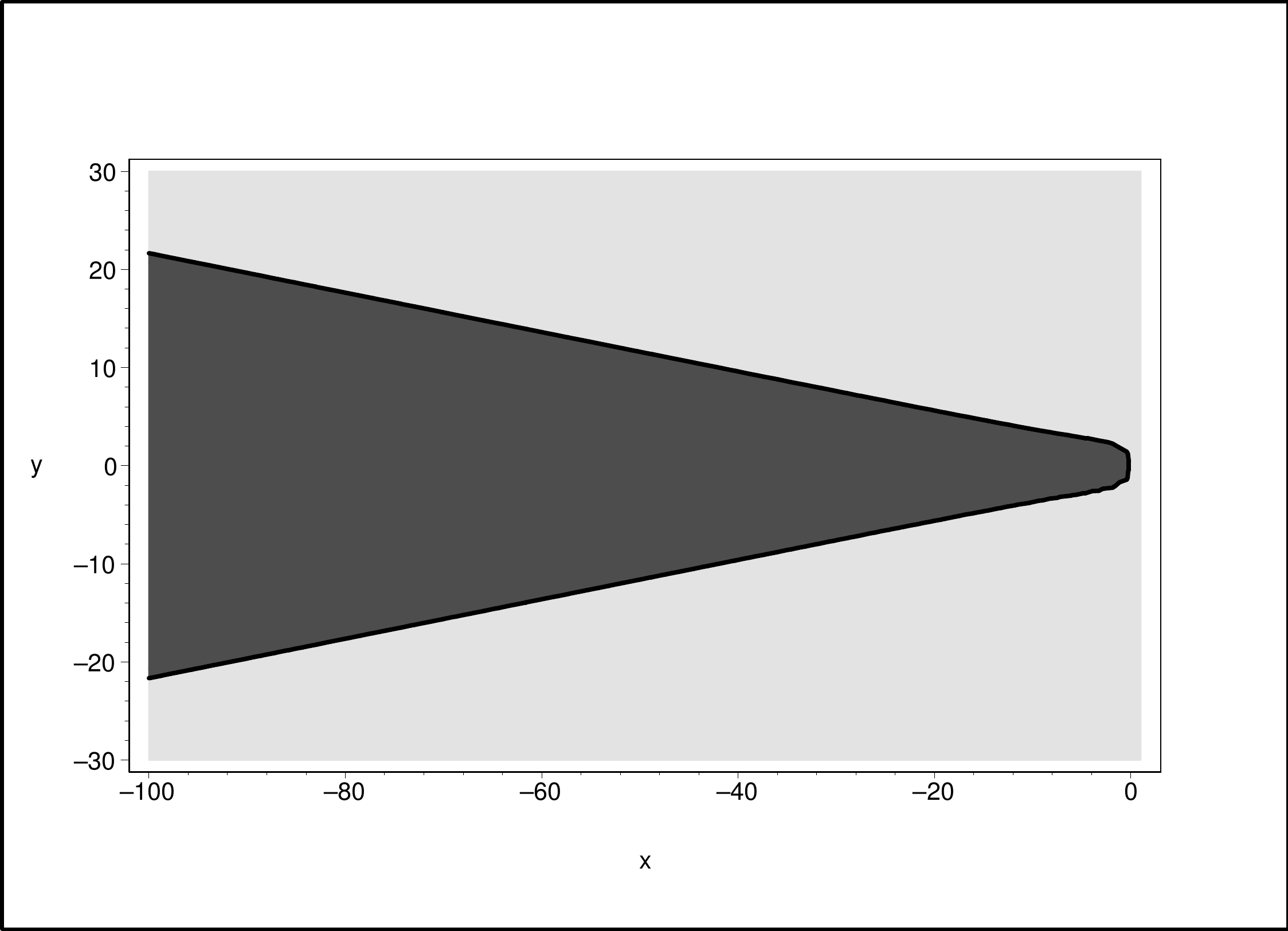}\hspace{0.5cm}
\includegraphics[width=4.7cm,angle=270]{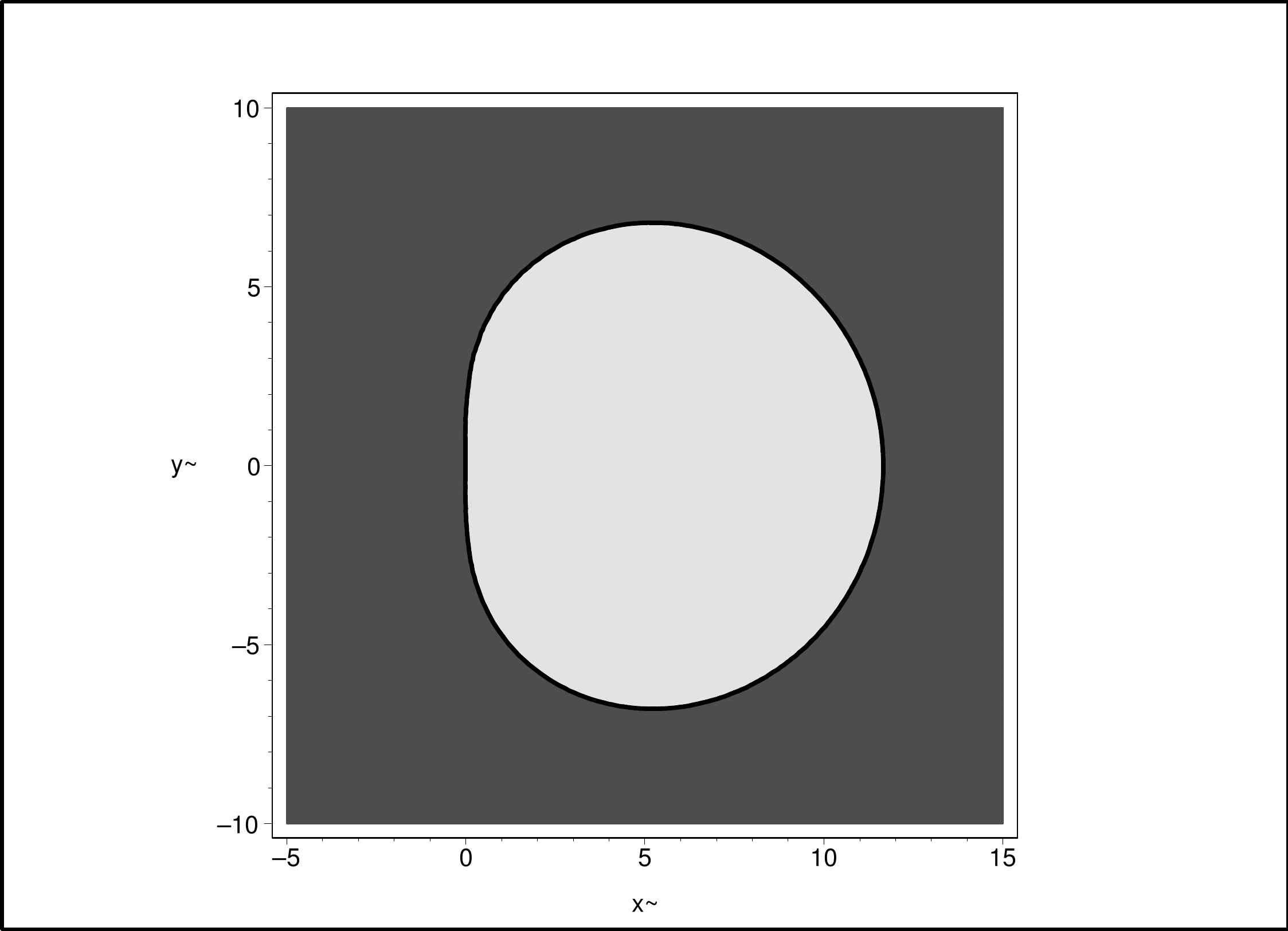}
\caption{Stability regions of IMEX method~(\ref{hig1}) (left) and $\tilde A$ (right).\mlabel{fig1}}
\end{center}
\end{figure}

The stability function of the implicit scheme $\tilde{A}$ is
\begin{equation}\mlabel{stabimplhig1}
R_{\tilde A}(z) = 2\,{\frac { \left( 1+\sqrt {2} \right)  \left( 1+z+\sqrt {2} \right) }
{ \left( 2-z+\sqrt {2} \right) ^{2}}}.
\end{equation}
A plot of the related stability region is shown in Figure~\ref{fig1} (right). The scheme $\tilde{A}$ appears to be $A$-stable
and satisfies $\lim_{\Re(z)\to-\infty}R_{\tilde A}(z)=0$, implying $L$-stability.

The dissipativity analysis for the implicit scheme defined by $\tilde{A}$ yields the
amplification factors for the standard
three-point space discretization (\ref{3pt}) and the fourth order stencil (\ref{4th}),
respectively. The amplification factors
are evaluated at the points $\theta\in\{0,\frac\pi4,\frac\pi2,\pi\}$
in Tables~\ref{evalsfimex13pt} and \ref{evalsfimex1trial}, respectively.
\begin{table}
\begin{center}
\begin{tabular}{||r|c||}
\hline
\multicolumn{1}{||c|}{$\theta$} & \multicolumn{1}{c||}{$g(\theta,\mu)$} \\
\hline\hline
$0$ &  $1$ \\ \hline
$\frac{\pi}{4}$ & $2\,{\frac { \left( 1+\sqrt {2} \right)  \left( 1-2\,\mu+\mu\,\sqrt {2}
+\sqrt {2} \right) }{ \left( -2-2\,\mu+\mu\,\sqrt {2}-\sqrt {2}\right) ^{2}}}$  \\ \hline
$\frac{\pi}{2}$ & $2\,{\frac { \left( 1+\sqrt {2} \right)  \left( 1-2\,\mu+\sqrt {2}
 \right) }{ \left( 2+2\,\mu+\sqrt {2} \right) ^{2}}}$  \\ \hline
$\pi$ & $2\,{\frac { \left( 1+\sqrt {2} \right)  \left( 1-4\,\mu+\sqrt {2}
 \right) }{ \left( 2+4\,\mu+\sqrt {2} \right) ^{2}}}$ \\ \hline
\end{tabular}
\caption{Values of $g(\theta,\mu)$ for some $\theta$, implicit scheme in (\ref{hig1}), three point space discretization (\ref{3pt}).\mlabel{evalsfimex13pt}}
\end{center}
\end{table}

\begin{table}
\begin{center}
\begin{tabular}{||r|c||}
\hline
\multicolumn{1}{||c|}{$\theta$} & \multicolumn{1}{c||}{$g(\theta,\mu)$} \\
\hline\hline
$0$ &  $1$ \\ \hline
$\frac{\pi}{4}$ & $12\,{\frac { \left( 1+\sqrt {2} \right)  \left( 6-15\,\mu+8\,\sqrt {2}
\mu+6\,\sqrt {2} \right) }{ \left( -12-15\,\mu+8\,\sqrt {2}\mu-6\,
\sqrt {2} \right) ^{2}}}$  \\ \hline
$\frac{\pi}{2}$ & $6\,{\frac { \left( 1+\sqrt {2} \right)  \left( 3-7\,\mu+3\,\sqrt {2}
 \right) }{ \left( 6+7\,\mu+3\,\sqrt {2} \right) ^{2}}}$  \\ \hline
$\pi$ & $6\,{\frac { \left( 1+\sqrt {2} \right)  \left( 3-16\,\mu+3\,\sqrt {2}
 \right) }{ \left( 6+16\,\mu+3\,\sqrt {2} \right) ^{2}}}$ \\ \hline
\end{tabular}
\caption{Values of $g(\theta,\mu)$ for some $\theta$, implicit scheme in (\ref{hig1}), fourth order space discretization (\ref{4th}).
\mlabel{evalsfimex1trial}}
\end{center}
\end{table}

The first positive zero of $g(\pi,\mu)$ is $\approx 0.6035$ for the three-point scheme (\ref{3pt}), where the function changes its sign,
and $|g(\pi,\mu)|$ never exceeds 1.
The first positive zero of $g(\pi,\mu)$ for the fourth order space discretization
(\ref{4th}) is $\approx 0.4526$, where the function changes its sign. The modulus never exceeds 1.

\subsubsection*{Modification of $\gamma$}
We may conceive of optimizing the method (\ref{hig1}) by adapting the value of the
parameter $\gamma$ in the definition of $\tilde A$ according to the
resulting stability, accuracy, and dissipativity properties.
The region of absolute monotonicity depends on $\gamma$ as follows \cite{higueras11private}:
\begin{eqnarray}
{\cal R}(\tilde A) &=&
\left\{\begin{array}{lll}
\frac{1}{1 -3 \gamma}, & \qquad &0\leq \gamma\leq \frac14,\\[2ex]
\frac{1-2 \gamma }{2 \gamma ^2-4 \gamma +1}-\sqrt{\frac{4 \gamma -1}{\left(2 \gamma ^2-4 \gamma
   +1\right)^2}}, & &1/4< \gamma < 1- \frac{1}{\sqrt{2}},\\[2ex]
 1+\sqrt{2}, & &\gamma=  1- \frac{1}{\sqrt{2}},  \\[2ex]
 \frac{1-2 \gamma }{2 \gamma ^2-4 \gamma +1}+\sqrt{\frac{4 \gamma -1}{\left(2 \gamma ^2-4 \gamma
   +1\right)^2}}, & &1- \frac{1}{\sqrt{2}}< \gamma \leq \frac12,
\end{array}\right.\mlabel{hig1modkra1}\\[2mm]
{\cal R}(A, \tilde A)&=&\left\{\begin{array}{ll} \left\{ (r,\tilde r): 0\leq r \leq 1 \, ,
0\leq \tilde r \leq \frac{1-r}{1-\gamma}\right\}\, ,   &  \quad 0\leq \gamma\leq \frac13, \\[2ex]
\left\{ (r,\tilde r): 0\leq r \leq \frac{1-2 \gamma}{\gamma} \, ,  0\leq \tilde r \leq \frac{1-r}{1-\gamma}\right\}\, ,
&  \quad \frac13\leq \gamma\leq \frac12.\end{array}\right.\mlabel{hig1modkra2}
\end{eqnarray}
A plot of the function ${\cal R}(\tilde A)$ in dependence of $\gamma$ is given in Figure~\ref{fig:ratildegamma}.

\begin{figure}[h]
\begin{center}
\includegraphics[width=6cm,angle=0]{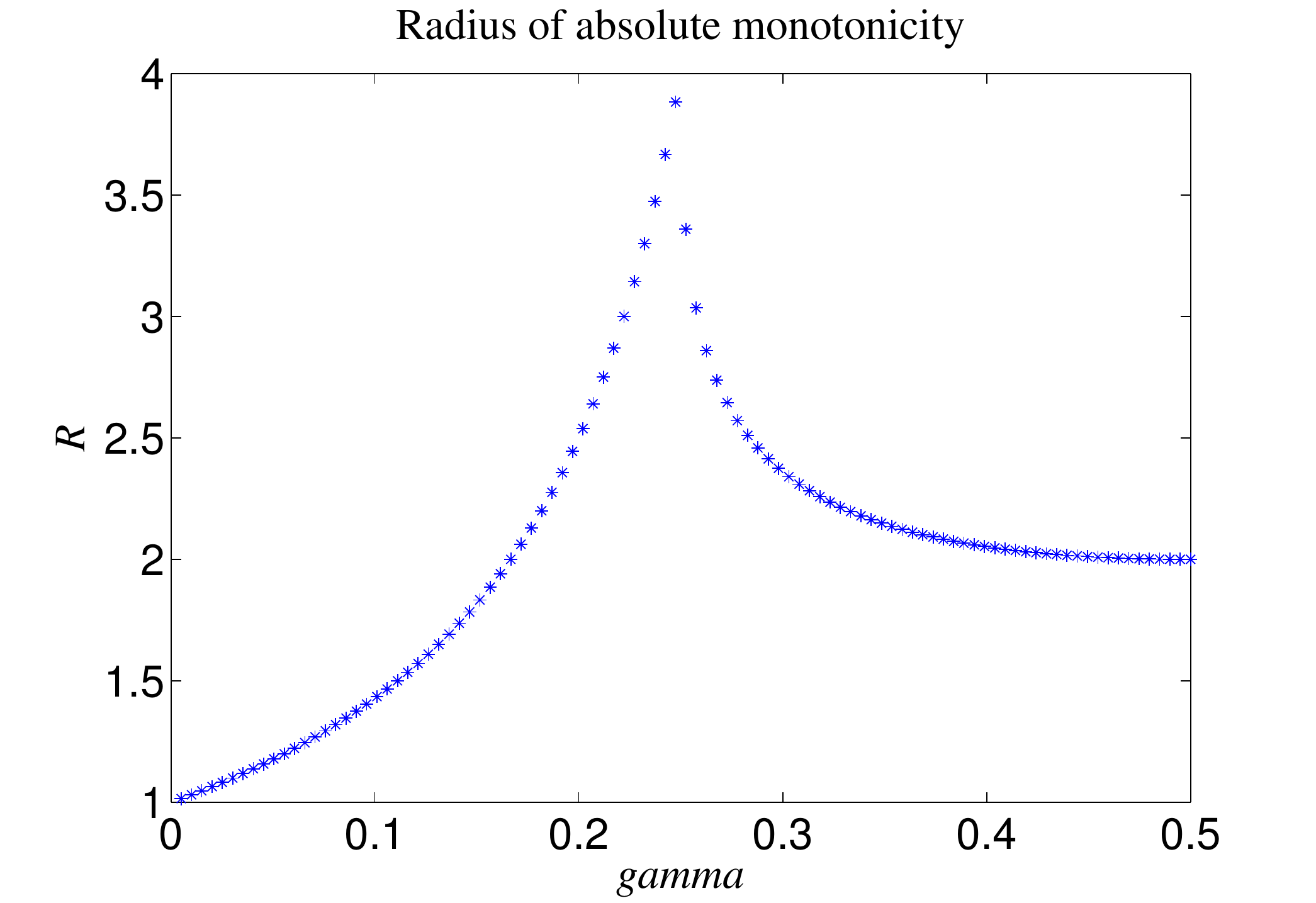}
\caption{Radius of absolute monotonicity ${\cal R}(\tilde A)$ as a function of $\gamma$ for (\ref{hig1}).\mlabel{fig:ratildegamma}}
\end{center}
\end{figure}

The regions of absolute monotonicity ${\cal R}(A, \tilde A)$ for the values $\gamma\in\{0.1,0.2,0.3\}$ are plotted
in Figure~\ref{fig:kraaimex1gamma}.

\begin{figure}[h]
\begin{center}
\includegraphics[width=4.3cm,angle=0]{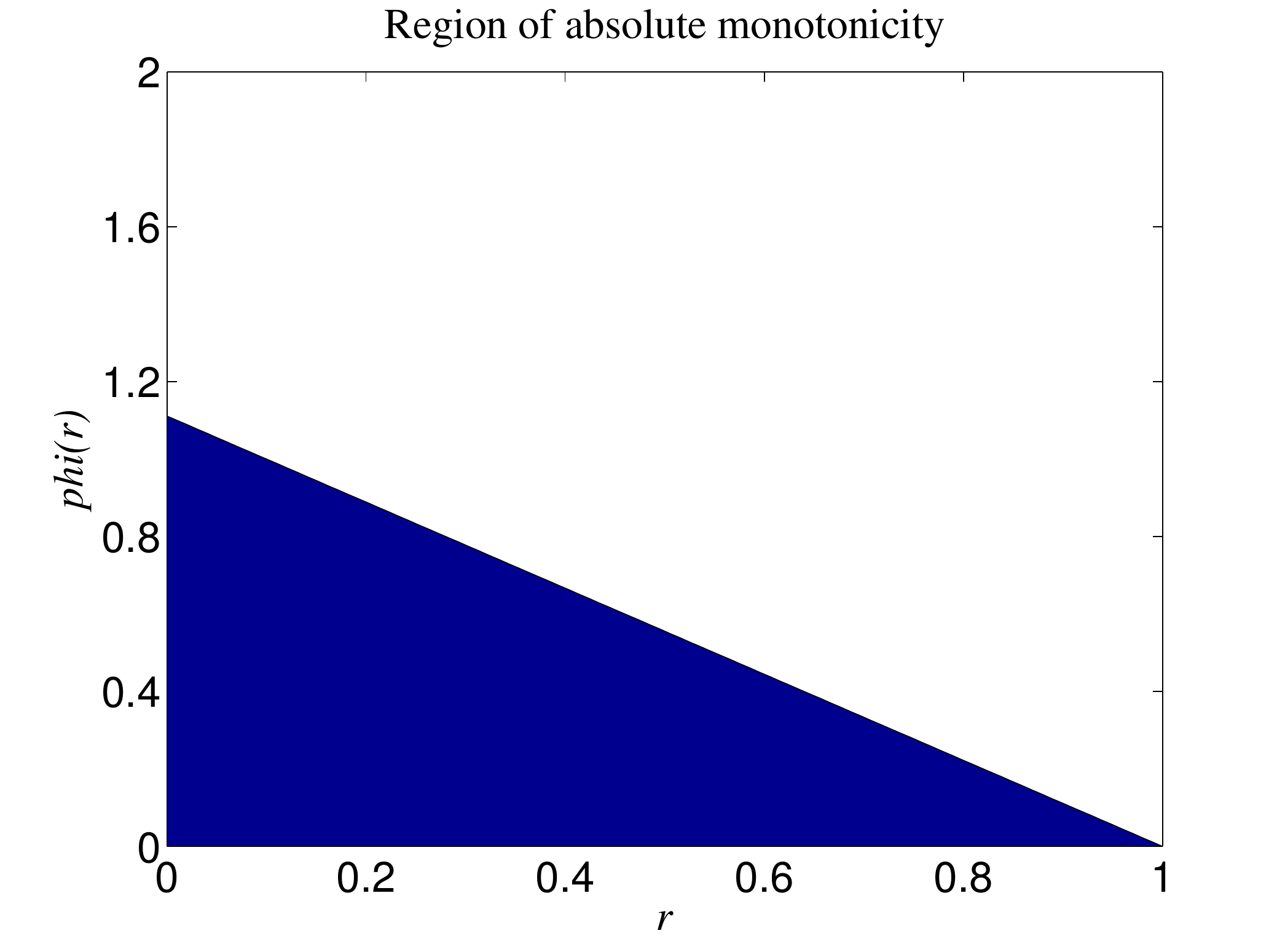}\hspace{3mm}
\includegraphics[width=4.3cm,angle=0]{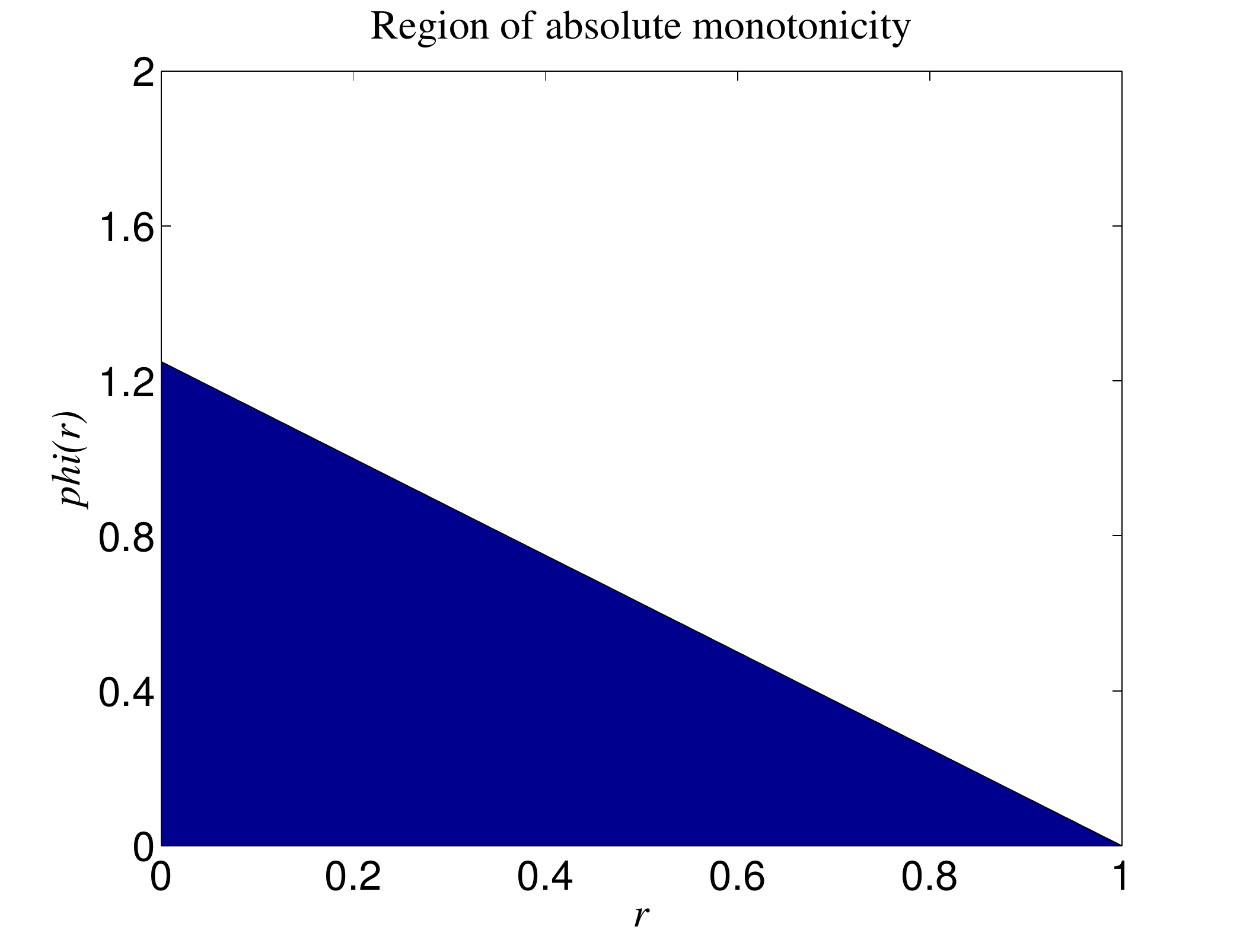}\hspace{3mm}
\includegraphics[width=4.3cm,angle=0]{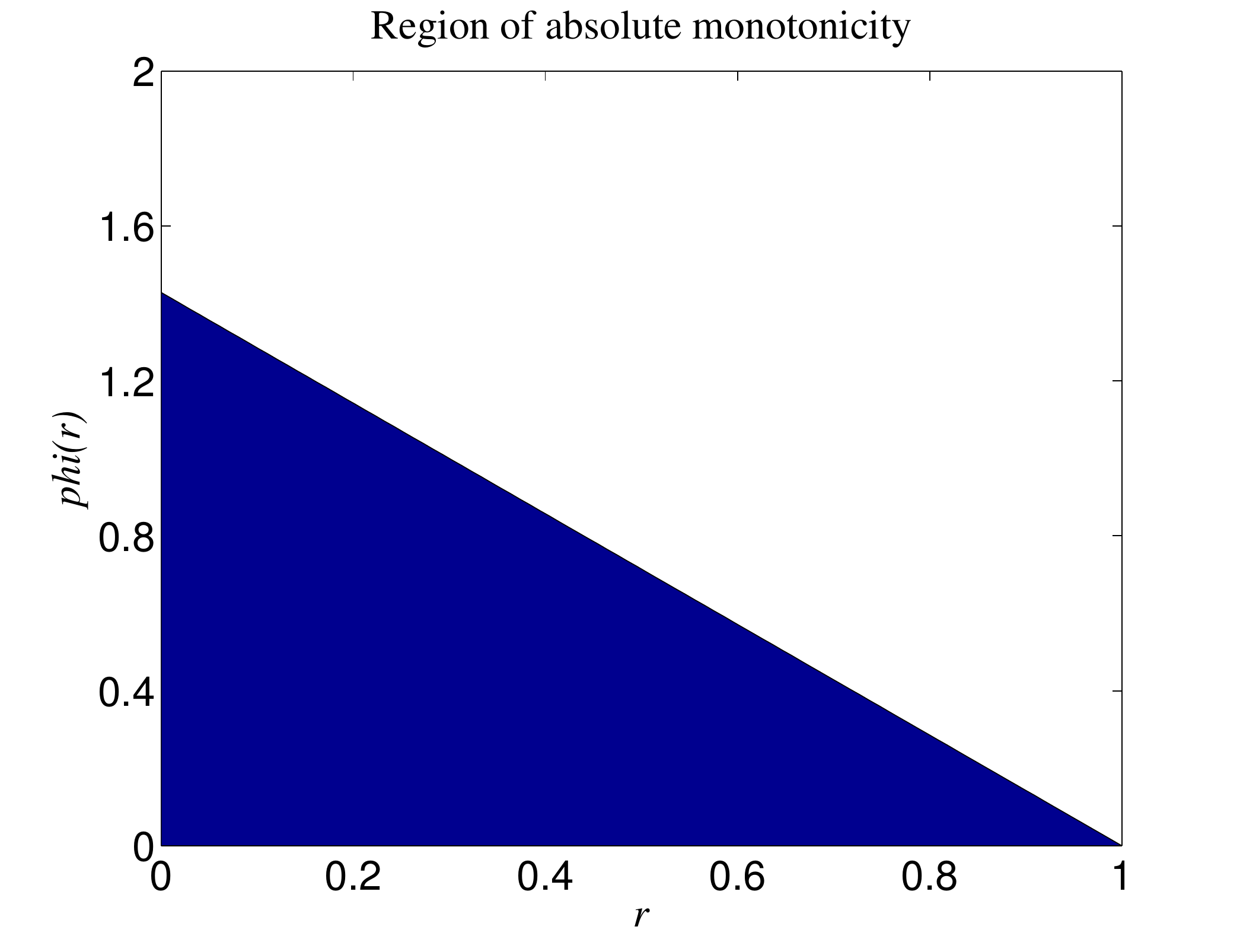}\hspace{3mm}
\caption{Regions ${\cal R}(A,\tilde A)$ for $\gamma\in\{0.1,0.2,0.3\}$
for (\ref{hig1}).\mlabel{fig:kraaimex1gamma}}
\end{center}
\end{figure}

The stability regions for the IMEX schemes for the different values of $\gamma$ cover bounded
subdomains of the left half plane for $\gamma<0.25$, while for
$\gamma\geq 0.25,$ the stability regions cover unbounded domains in the left
half plane. In fact, the left boundaries $z_\mathrm{left}$ satisfy
$$z_\mathrm{left} = \left\{ \begin{array}{ll} \frac2{4\gamma-1}, & \quad \gamma < 0.25,\\
-\infty, & \quad \gamma \geq 0.25. \end{array} \right.$$
However, even in the cases with unbounded stability regions, in general
$\lim_{\Re(z)\to-\infty}R(z) \neq 0$. The boundaries of the stability regions are plotted
in Figure~\ref{fig:stabimex1modboundary} (left), where values equal to 0 represent unbounded
stability regions. The implicit schemes $\tilde A$ show the same stability behaviour
concerning both the boundaries of the stability regions and the limits for $\Re(z)\to-\infty$.

\begin{figure}[h]
\begin{center}
\includegraphics[width=6cm,angle=0]{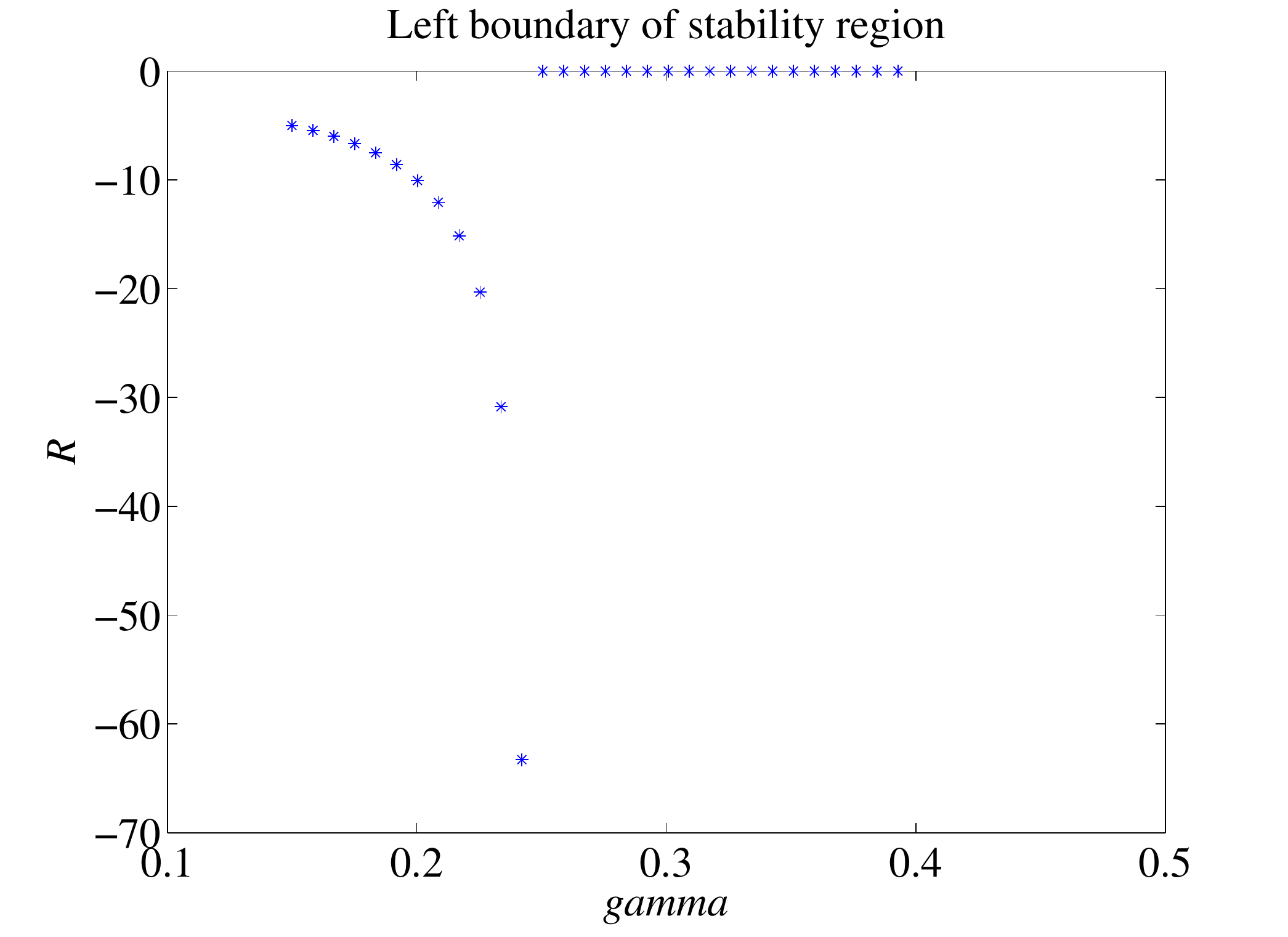} \hspace*{0.5cm}
\includegraphics[width=6cm,angle=0]{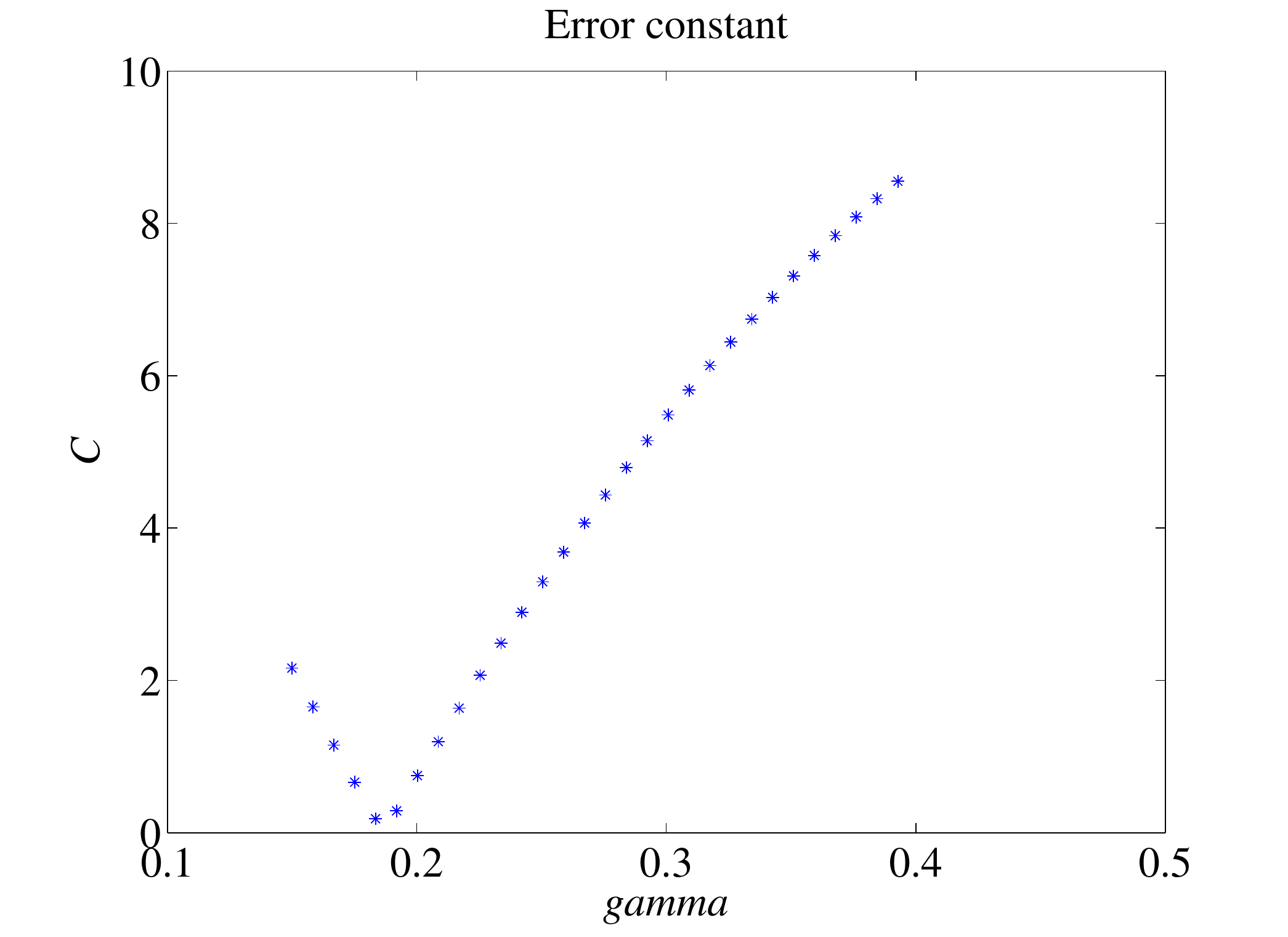}
\caption{Left boundary of the stability region (left) and error constant
computed for (\ref{eq:test1}) (right)
as a function of $\gamma$ for (\ref{hig1}).\mlabel{fig:stabimex1modboundary}}
\end{center}
\end{figure}

It was demonstrated with a \textsc{Matlab} implementation that the
convergence order two is retained also for the modified values of
$\gamma$. The error constant depends on $\gamma$, however.
In Figure~\ref{fig:stabimex1modboundary} (right) we plot the error constant as a
function of $\gamma$, where the error is determined at $t=1.3$ for the test
problem (\ref{eq:test1}).
We note that for small $\gamma$, the error constant decreases as
$\gamma$ grows, while for $\gamma>0.1833$ the constant grows
monotonically. This behaviour does not appear to be related
to the results we obtain for the dissipativity analysis,
where for $\gamma=0.25$ the behaviour changes.

\begin{figure}[h]
\begin{center}
\includegraphics[width=6cm,height=4cm,angle=0]{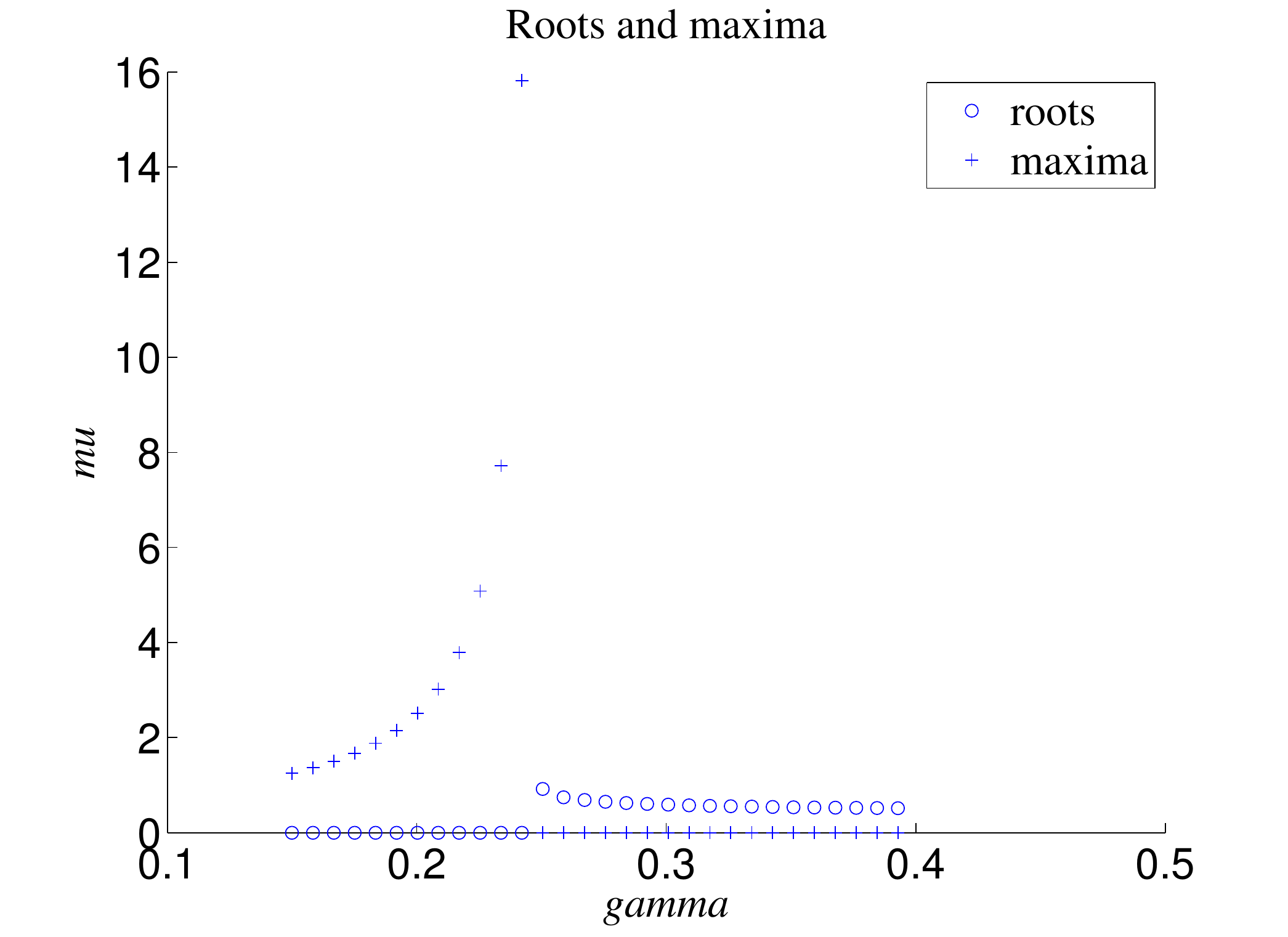} \hspace{0.5cm}
\includegraphics[width=6cm,height=4cm,angle=0]{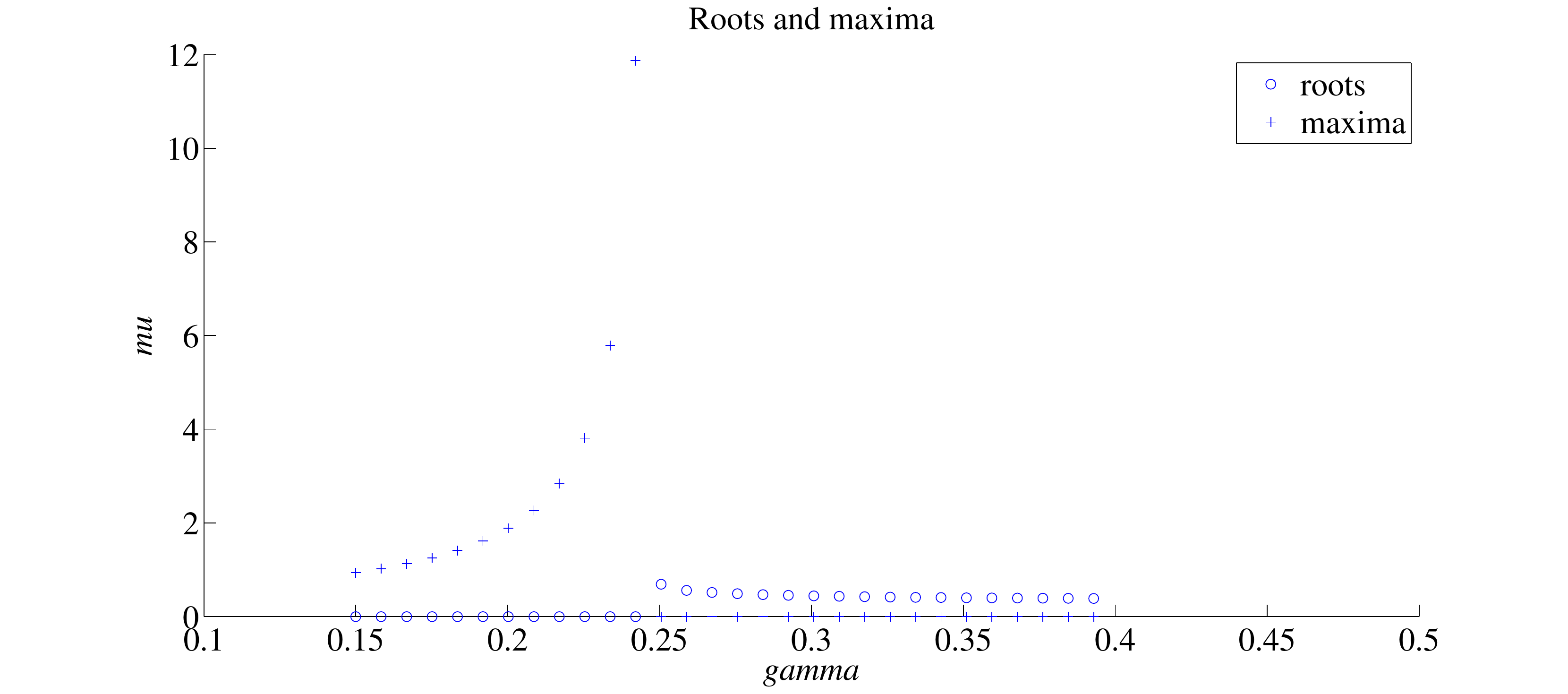}
\caption{Dissipativity analysis for (\ref{hig1}) for different values of $\gamma$,
three point space discretization (\ref{3pt}) (left) and fourth order discretization (\ref{4th}) (right).\mlabel{fig:varygamma1}}
\end{center}
\end{figure}

To assess the dissipativity properties of the modified scheme
we vary $\gamma$ and compute the first positive zero of $g(\pi,\mu)$
and the point $\mu$ where the modulus of this function exceeds 1.

For $\gamma\geq 0.25$, the amplification factor for the three-point space
discretization (\ref{3pt}) has a zero at
$\mu = \frac{1-2\gamma-\sqrt{4\gamma - 1}}{8\gamma^2-16\gamma+4}$, but there is no
real root for $\gamma< 0.25$.
Also for $\gamma< 0.25$, $|g|$ exceeds 1 at $\mu = \frac1{2-8\gamma}$,
and this behaviour does not occur for $\gamma\geq 0.25$.
For the discretization (\ref{4th}), the
behaviour is the same, but the scale with respect to $\mu$ is multiplied by 0.75
(see also Section~3.2 in \cite{kochetal10a}).

To illustrate this analysis, in Figure~\ref{fig:varygamma1} we vary $\gamma$ and compute for 30 values of $\gamma$
spaced equidistantly in the interval $[0.15,1.1-1/\sqrt{2}]$ the first positive zero of $g(\pi,\mu)$ and the point $\mu$
where the modulus of this function exceeds 1. This analysis is given for the three point
space discretization (\ref{3pt}) in Figure~\ref{fig:varygamma1}.
We conclude that it may be of interest to choose $\gamma<0.25$ to ensure $0<g(\pi,\mu)<1$.
This can be achieved by a value of $\gamma$ just slightly smaller than $0.25$.

Taking into account the results displayed in Fig.~\ref{fig:stabimex1modboundary} we suggest
to optimize $\gamma$ by taking it as large as necessary to avoid any linear stability restrictions
due to the term $G(y(t))$ in (\ref{de1}) and as small as possible to minimize the stability
constant $C$. Note that some special values for $\gamma$ are 0, $(1-1/\sqrt{3})/2 \approx 0.2113$,
and $1/4$, in addition to the originally proposed value of $1-1/\sqrt{2} \approx 0.2929$. They are readily
identified to yield the classical explicit SSPRK(2,2) method\footnote{We will henceforth use the common
specification `SSPRK($s$,$p$)' introduced in \cite{kraaijevanger91} for an $s$-stage order $p$ explicit
strong stability preserving Runge--Kutta method.} of order two by Shu \& Osher \cite{shu88a}, the optimal
implicit third order SSP method with two stages \cite{ferspi08a}, and the optimal implicit second order
SSP method with two stages \cite{ferspi08a}.
In our numerical tests reported in Section~\ref{num} below, we found that the choice
$\gamma=0.24$ yielded the most efficient time integrator.

\subsection*{An SSP2(3,3,2) Method}
\cite{higueras06} gives an IMEX SSP2(3,3,2) method with nontrivial region
of absolute monotonicity:
\begin{equation}\mlabel{hig2}
\begin{array}{c|ccc}
0       & 0       & 0       & 0 \\[1mm]
\frac12 & \frac12 & 0       & 0 \\[1mm]
1       & \frac12 & \frac12 & 0 \\[1mm] \hline\\[-3mm]
A       & \frac13 & \frac13 & \frac13
\end{array}
\hspace*{1cm}
\begin{array}{c|ccc}
\frac15    & \frac15    & 0       & 0 \\[1mm]
\frac3{10} & \frac1{10} & \frac15 & 0 \\[1mm]
1          & \frac13    & \frac13 & \frac13 \\[1mm] \hline\\[-3mm]
\tilde A   & \frac13    & \frac13 & \frac13
\end{array}
\end{equation}
This is a modification of a scheme from \cite{parrus05}, where the
latter turned out to have a trivial region of absolute monotonicity.
It holds that $\mathcal{R}(A)=2$ and $R(\tilde A)=\frac59(\sqrt{70}-4)$,
and
$$ \mathcal{R}(A,\tilde A) = \{(r,\tilde r): 0 \leq r \leq 1,\ 0 \leq \tilde r \leq \phi(r)\}, $$
where
%
%
$$ \phi(r) = \left\{ 
\frac14(-28+9r)+\frac14\sqrt{1264 - 984 r +201 r^2} \right\}.$$
We note that the latter is a correction with respect to \cite{higueras06}, since we have found
$r$ to be necessarily bound by $1$ in $ \mathcal{R}(A,\tilde A)$.
A plot of $\mathcal{R}(A,\tilde A)$ is given in Figure~\ref{kra2}.

\begin{figure}[h]
\begin{center}
\includegraphics[width=6cm,angle=0]{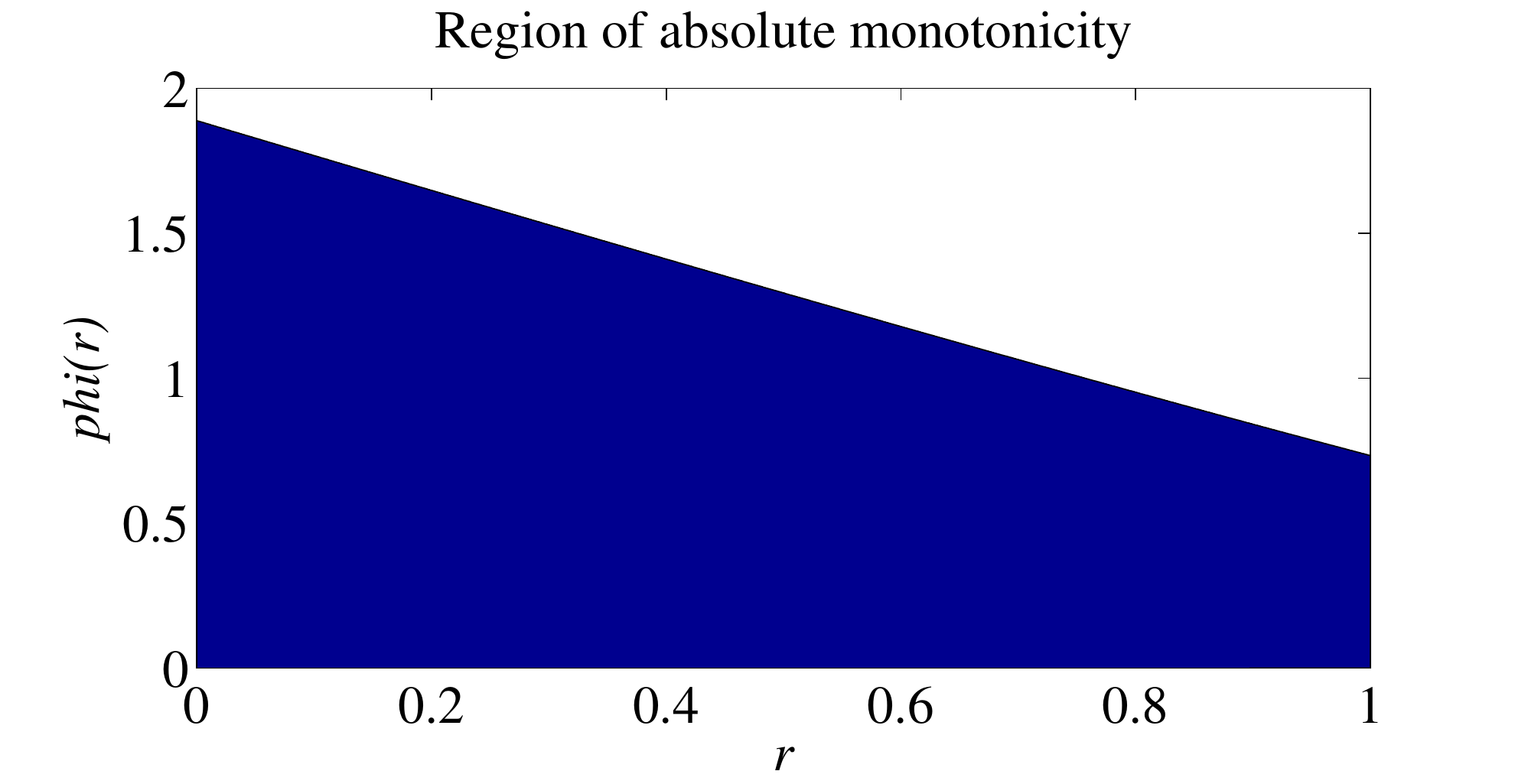}
\caption{Region of absolute monotonicity for method~(\ref{hig2}).\mlabel{kra2}}
\end{center}
\end{figure}

A plot of the stability region is given in Figure~\ref{fig2} (left).
We observe that the stability region is
tangent to the imaginary axis and appears to be unbounded as $\Re(z)\to-\infty$.
Moreover, $\lim_{\Re(z)\to -\infty}R(z) =0$.
\begin{figure}[h]
\begin{center}
\includegraphics[width=4.7cm,angle=270]{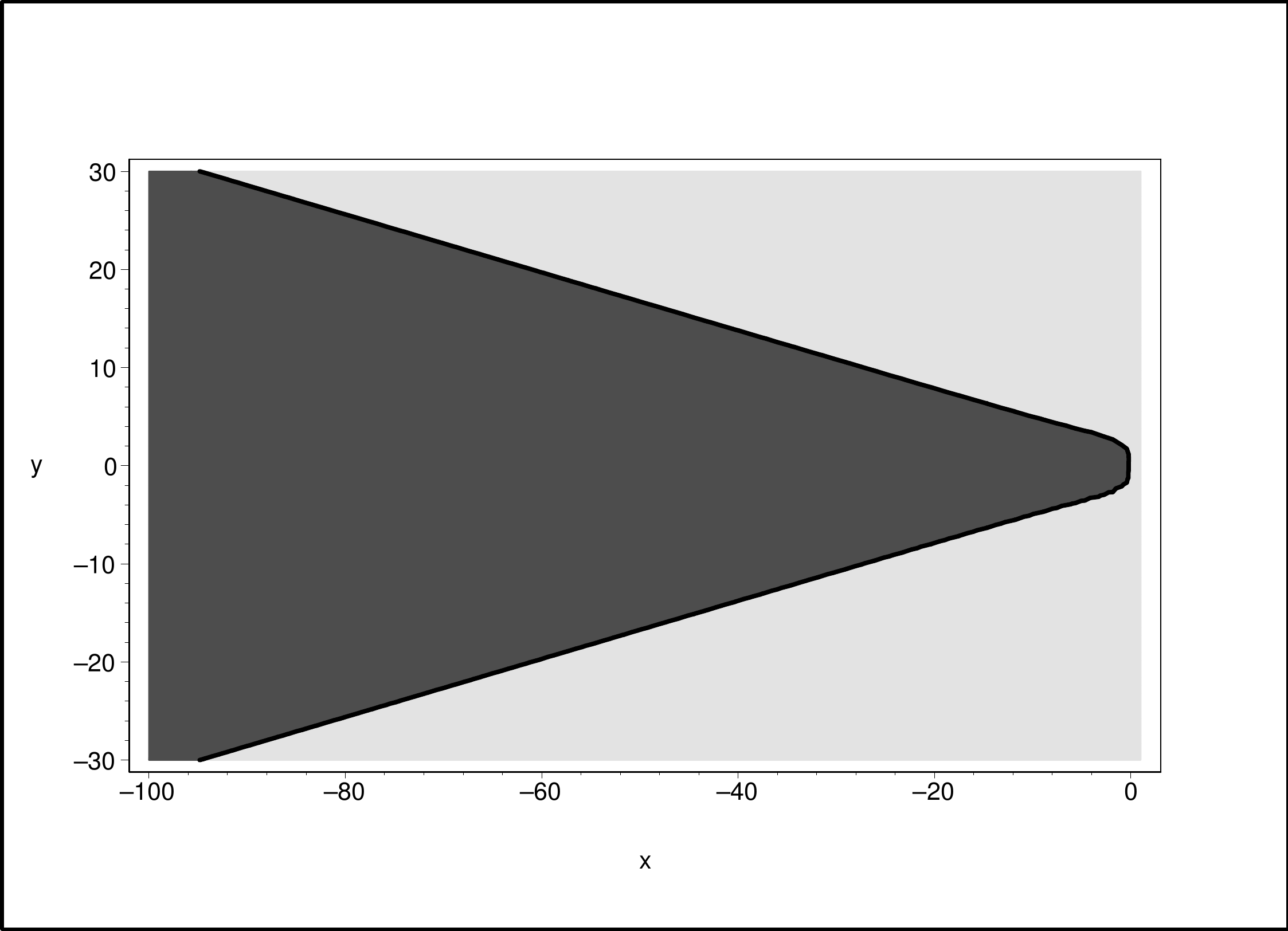}\hspace{0.5cm}
\includegraphics[width=4.7cm,angle=270]{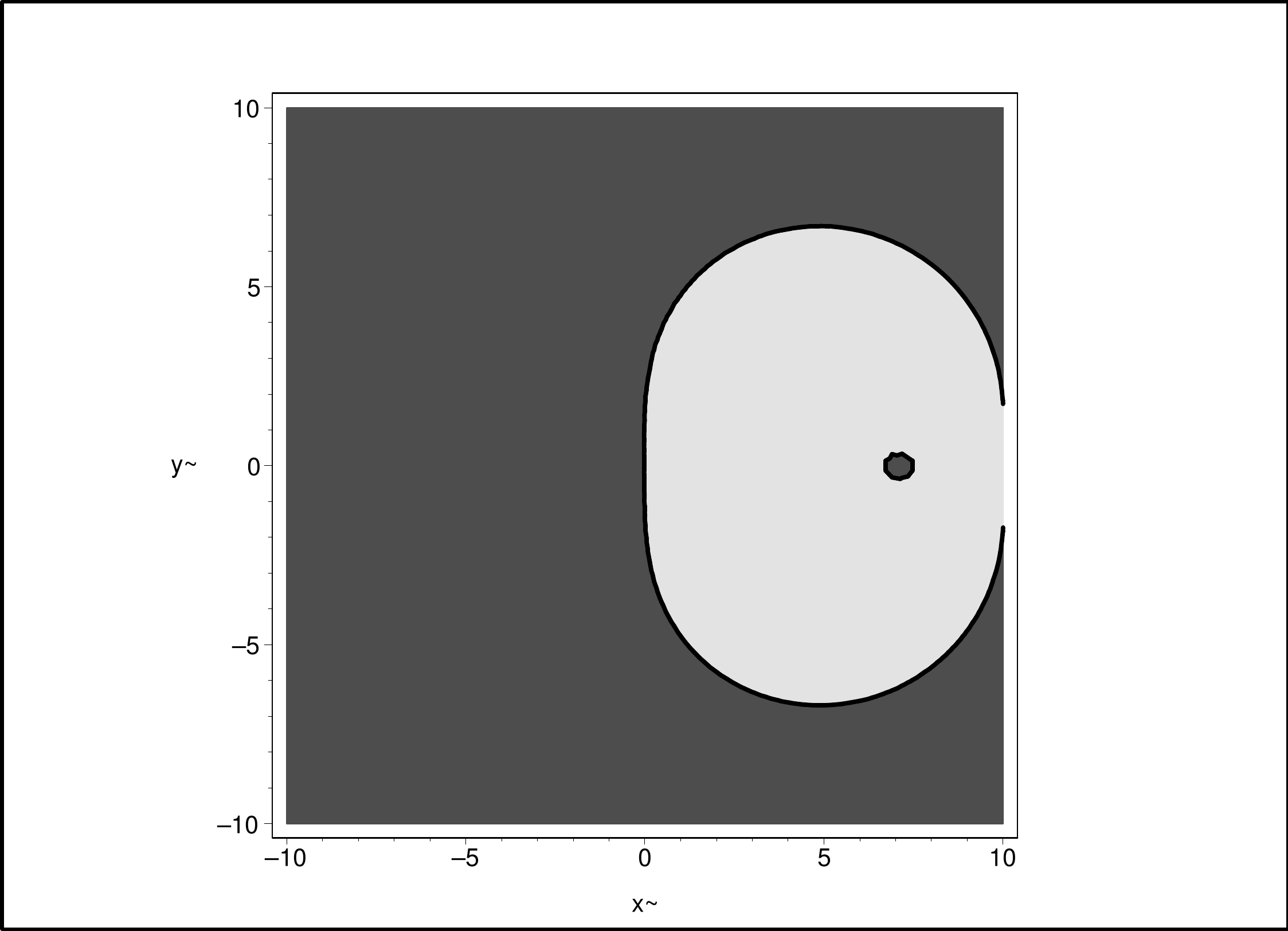}
\caption{Stability region of IMEX method~(\ref{hig2}) (left) and $\tilde A$ (right).\mlabel{fig2}}
\end{center}
\end{figure}
The stability function of the implicit scheme is given by
\begin{equation}\mlabel{stabimplhig2}
R_{\tilde A}(z) = {\frac {-150-40\,z+9\,{z}^{2}}{2 \left( -5+z \right) ^{2} \left( -
3+z \right) }}.
\end{equation}
A plot is shown in Figure~\ref{fig2} (right). Note that the stability region is not connected. The
method is $A$-stable, however.
Moreover, the scheme $\tilde{A}$ satisfies $\lim_{\Re(z)\to-\infty}R_{\tilde A}(z)=0$.

The dissipativity analysis for the implicit scheme defined by $\tilde{A}$ yields the
amplification factors  for the standard
three-point space discretization (\ref{3pt}) and the fourth order stencil (\ref{4th}).
These amplification factors are evaluated at the points $\theta\in\{0,\frac\pi4,\frac\pi2,\pi\}$
in Tables~\ref{evalsfimex23pt} and \ref{evalsfimex2trial}, respectively.
\begin{table}
\begin{center}
\begin{tabular}{||r|c||}
\hline
\multicolumn{1}{||c|}{$\theta$} & \multicolumn{1}{c||}{$g(\theta,\mu)$} \\
\hline\hline
$0$ &  $1$ \\ \hline
$\frac{\pi}{4}$ & $-{\frac {75+20\,\mu\,\sqrt {2}-40\,\mu-27\,{\mu}^{2}+18\,{\mu}^{2}
\sqrt {2}}{ \left( -5+\mu\,\sqrt {2}-2\,\mu \right) ^{2} \left( -3+\mu
\,\sqrt {2}-2\,\mu \right) }}$  \\ \hline
$\frac{\pi}{2}$ & $-{\frac {-75+18\,{\mu}^{2}+40\,\mu}{ \left( 5+2\,\mu \right) ^{2}
 \left( 3+2\,\mu \right) }}$  \\ \hline
$\pi$ & $-{\frac {-75+80\,\mu+72\,{\mu}^{2}}{ \left( 5+4\,\mu \right) ^{2}
 \left( 3+4\,\mu \right) }}$ \\ \hline
\end{tabular}
\caption{Values of $g(\theta,\mu)$ for some $\theta$, implicit scheme in (\ref{hig2}), three point space discretization (\ref{3pt}).
\mlabel{evalsfimex23pt}}
\end{center}
\end{table}

\begin{table}
\begin{center}
\begin{tabular}{||r|c||}
\hline
\multicolumn{1}{||c|}{$\theta$} & \multicolumn{1}{c||}{$g(\theta,\mu)$} \\
\hline\hline
$0$ &  $1$ \\ \hline
$\frac{\pi}{4}$ & $-9\,{\frac {1800-1200\,\mu+640\,\mu\,\sqrt {2}-1059\,{\mu}^{2}+720\,{
\mu}^{2}\sqrt {2}}{ \left( -30-15\,\mu+8\,\mu\,\sqrt {2} \right) ^{2}
 \left( -18-15\,\mu+8\,\mu\,\sqrt {2} \right) }}$  \\ \hline
$\frac{\pi}{2}$ & $-9/2\,{\frac {-450+280\,\mu+147\,{\mu}^{2}}{ \left( 7\,\mu+15 \right)
^{2} \left( 7\,\mu+9 \right) }}$  \\ \hline
$\pi$ & $-9\,{\frac {320\,\mu+384\,{\mu}^{2}-225}{ \left( 15+16\,\mu \right) ^{
2} \left( 9+16\,\mu \right) }}$ \\ \hline
\end{tabular}
\caption{Values of $g(\theta,\mu)$ for some $\theta$, implicit scheme in (\ref{hig2}), fourth order space discretization (\ref{4th}).
\mlabel{evalsfimex2trial}}
\end{center}
\end{table}

The first positive zero of $g(\pi,\mu)$ is $\approx 0.6064$ for the three-point space discretization (\ref{3pt}),
where the function changes its sign, and $|g(\pi,\mu)|$ never exceeds~1.
The first positive zero of $g(\pi,\mu)$ is $\approx 0.4551$ for the fourth order space discretization (\ref{4th}),
where the function changes its sign, and $|g(\pi,\mu)|$ never exceeds~1.

\subsection*{An SSP3(3,3,3) Method}
\cite{higueras09} gives the following SSP3(3,3,3) method with nontrivial region
of absolute monotonicity:
\begin{equation}\mlabel{hig3}
\begin{array}{c|ccc}
0       & 0       & 0       & 0 \\[1mm]
1       & 1       & 0       & 0 \\[1mm]
\frac12 & \frac14 & \frac14 & 0 \\[1mm] \hline\\[-3mm]
A       & \frac16 & \frac16 & \frac23
\end{array}
\hspace*{1cm}
\begin{array}{c|ccc}
0        & 0             & 0          & 0          \\[1mm]
1        & \frac{14}{15} & \frac1{15} & 0          \\[1mm]
\frac12  & \frac7{30}    & \frac15    & \frac1{15} \\[1mm] \hline\\[-3mm]
\tilde A & \frac16       & \frac16    & \frac23
\end{array}
\end{equation}
It holds $\mathcal{R}(A)=1$ and $\mathcal{R}(\tilde
A)=\frac5{47}(13-2\sqrt7)$, and
$$ \mathcal{R}(A,\tilde A) = \{(r,\tilde r): 0 \leq r \leq 1,\ 0 \leq \tilde r \leq \phi(r)\}, $$
where
$$ \phi(r) = \frac{15}{302}\left(28-25r-\sqrt{180-192r+21r^2}\right). $$
A plot of $\mathcal{R}(A,\tilde A)$ is given in Figure~\ref{kra3}.

\begin{figure}[ht]
\begin{center}
\includegraphics[width=6cm,angle=0]{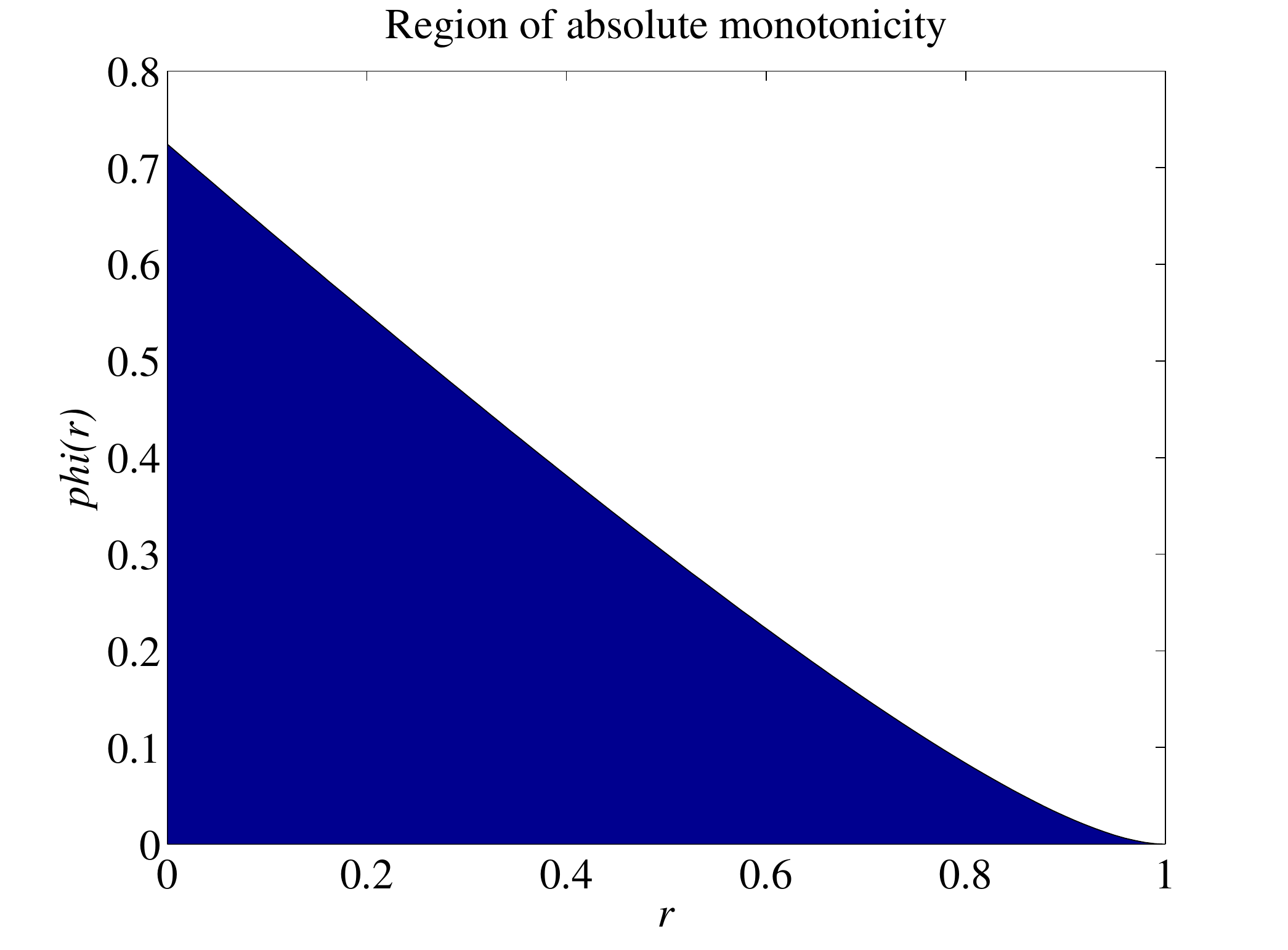}
\caption{Region of absolute monotonicity for method~(\ref{hig3}).\mlabel{kra3}}
\end{center}
\end{figure}

\begin{figure}[t]
\begin{center}
\includegraphics[width=4.7cm,angle=270]{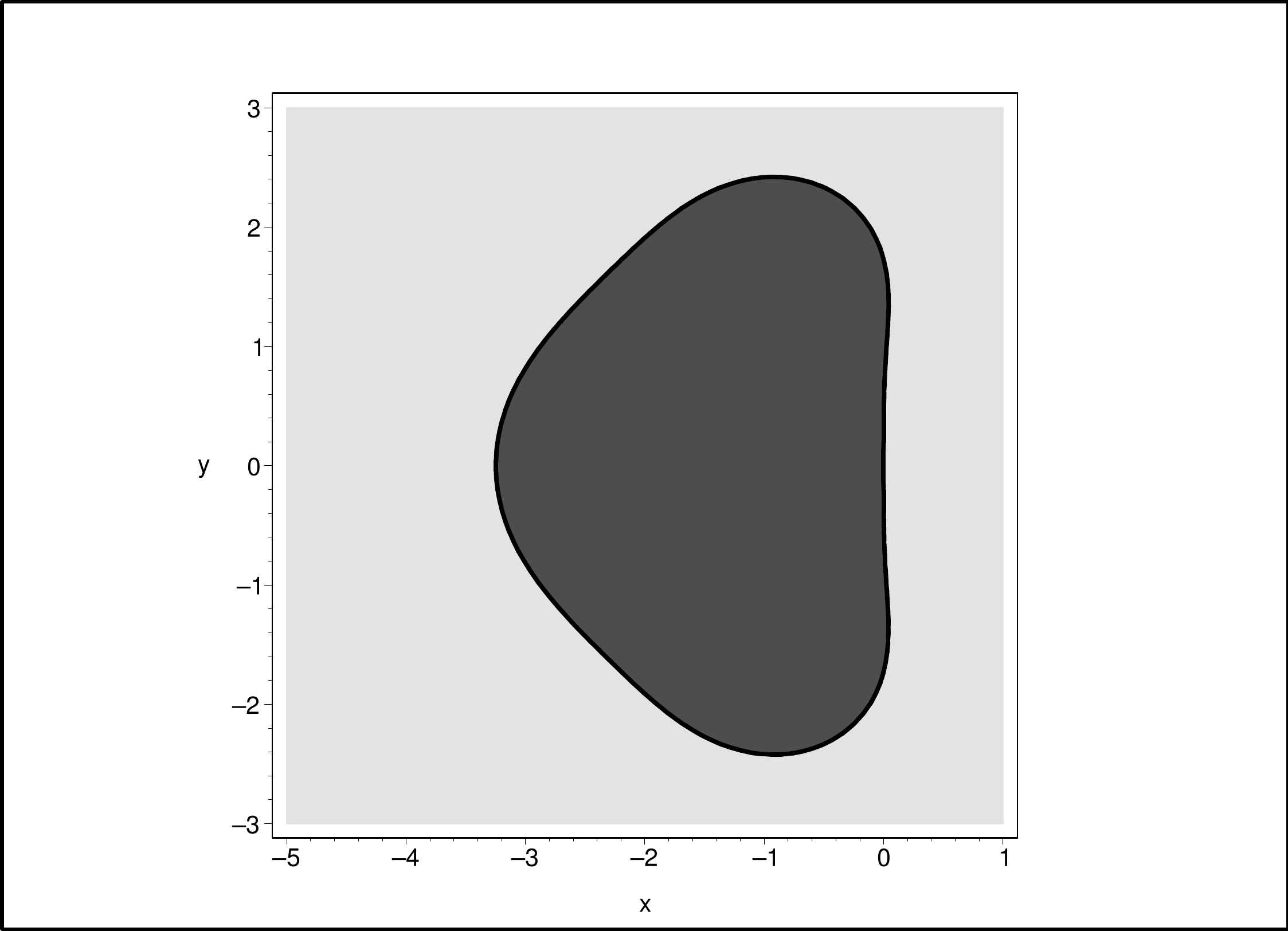} \hspace*{0.5cm}
\includegraphics[width=4.7cm,angle=270]{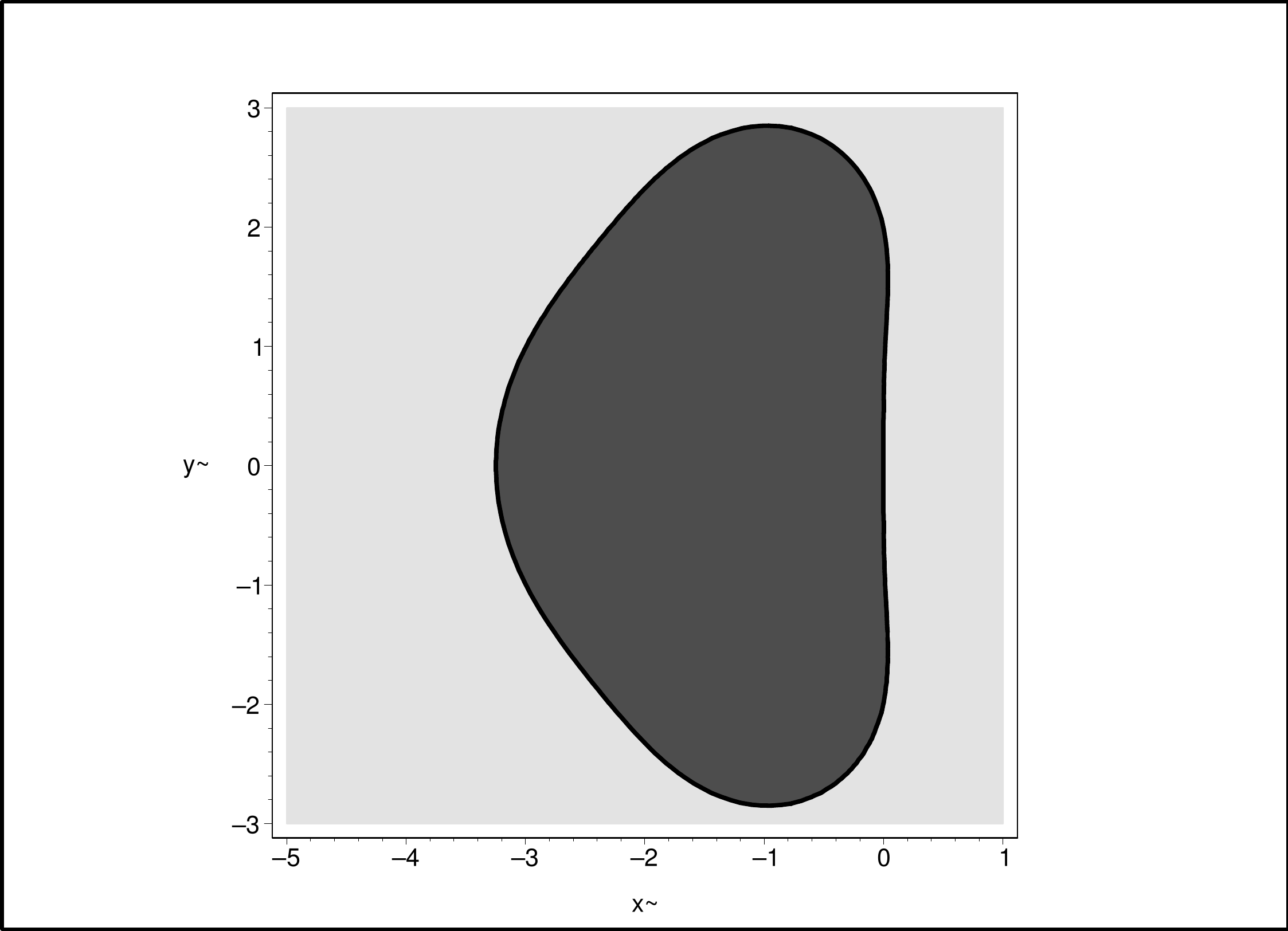}
\caption{Stability regions of IMEX method~(\ref{hig3}) (left) and implicit method $\tilde A$ (right).\mlabel{fig3}}
\end{center}
\end{figure}

The stability region occupies a bounded domain in the
negative half-plane, and slightly overlaps the imaginary axis, see
Figure~\ref{fig3} (left).
Note that $\lim_{\Re(z)\to -\infty}|R(z)| =\infty$.
The stability function of the implicit scheme is given by
\begin{equation}\mlabel{stabimplhig3}
R_{\tilde A}(z) = {\frac {450+390\,z+167\,{z}^{2}+47\,{z}^{3}}{2 \left( -15+z \right) ^{2}}}.
\end{equation}
A plot is shown in Figure~\ref{fig3} (right). The point where the stability region intersects the negative real
half-line is located at $x\approx-3.248$ for the IMEX scheme. The same value is computed
for the stability region of the implicit scheme. However, the stability region of $\tilde A$
extends further along the imaginary axis.
Finally, $\lim_{\Re(z)\to -\infty}|R_{\tilde A}(z)| =\infty$.

The dissipativity analysis for the implicit scheme defined by $\tilde{A}$ yields the
amplification factors for the standard
three-point space discretization (\ref{3pt}) and for the fourth order stencil (\ref{4th}).
These amplification factors are evaluated at the points $\theta\in\{0,\frac\pi4,\frac\pi2,\pi\}$
in Tables~\ref{evalsfimex33pt} and \ref{evalsfimex3trial}, respectively.
\begin{table}
\begin{center}
\begin{tabular}{||r|c||}
\hline
\multicolumn{1}{||c|}{$\theta$} & \multicolumn{1}{c||}{$g(\theta,\mu)$} \\
\hline\hline
$0$ &  $1$ \\ \hline
$\frac{\pi}{4}$ & ${\frac {225+195\,\mu\,\sqrt {2}-390\,\mu+501\,{\mu}^{2}-334\,{\mu}^{2}
\sqrt {2}+329\,{\mu}^{3}\sqrt {2}-470\,{\mu}^{3}}{ \left( -15+\mu\,
\sqrt {2}-2\,\mu \right) ^{2}}}$  \\ \hline
$\frac{\pi}{2}$ & $-{\frac {-225+390\,\mu-334\,{\mu}^{2}+188\,{\mu}^{3}}{ \left( 15+2\,
\mu \right) ^{2}}}$  \\ \hline
$\pi$ & $-{\frac {1504\,{\mu}^{3}+780\,\mu-1336\,{\mu}^{2}-225}{ \left( 15+4\,
\mu \right) ^{2}}}$ \\ \hline
\end{tabular}
\caption{Values of $g(\theta,\mu)$ for some $\theta$, implicit scheme in (\ref{hig3}), three point space discretization.
\mlabel{evalsfimex33pt}}
\end{center}
\end{table}

\begin{table}
\begin{center}
\begin{tabular}{||r|c||}
\hline
\multicolumn{1}{||c|}{$\theta$} & \multicolumn{1}{c||}{$g(\theta,\mu)$} \\
\hline\hline
$0$ &  $1$ \\ \hline
$\frac{\pi}{4}$ & $1/12\,{\frac {97200-210600\,\mu+112320\,\mu\,\sqrt {2}+353706\,{\mu}^{
2}-240480\,{\mu}^{2}\sqrt {2}-429345\,{\mu}^{3}+301928\,{\mu}^{3}
\sqrt {2}}{ \left( -90-15\,\mu+8\,\mu\,\sqrt {2} \right) ^{2}}}$  \\ \hline
$\frac{\pi}{2}$ & $-1/6\,{\frac {-12150+24570\,\mu-24549\,{\mu}^{2}+16121\,{\mu}^{3}}{
 \left( 45+7\,\mu \right) ^{2}}}$  \\ \hline
$\pi$ & $-1/3\,{\frac {-6075+28080\,\mu-64128\,{\mu}^{2}+96256\,{\mu}^{3}}{
 \left( 45+16\,\mu \right) ^{2}}}$ \\ \hline
\end{tabular}
\caption{Values of $g(\theta,\mu)$ for some $\theta$, implicit scheme in (\ref{hig3}), fourth order space discretization.
\mlabel{evalsfimex3trial}}
\end{center}
\end{table}

For the three-point space discretization (\ref{3pt}), the first positive zero of $g(\pi,\mu)$ is $\approx 0.4650$,
where the function changes its sign, and $g(\pi,\mu)=-1$ for $\mu\approx 0.8122$.
The first positive zero of $g(\pi,\mu)$ is $\approx 0.3488$ for the fourth-order space discretization (\ref{4th}),
where the function changes its sign, and $g(\pi,\mu)=-1$ at $\mu=0.6093$.

\section{Numerical Experiments}\mlabel{num}

In this section, we present the results of the simulations performed with the time integrators
discussed in this paper and compare their performance to the classical explicit 2-stages second
order and 3-stages third order methods of Shu \& Osher \cite{shu88a}, and the
3-stages second order method of \cite{kraaijevanger91}. They are the optimum SSP RK methods
for their given number of stages and order and we refer to them by their common technical
specifications `SSPRK(2,2)' and `SSPRK(3,3)', as well as `SSPRK(3,2)', respectively.
The coefficients of the SSPRK(2,2) and SSPRK(3,3) methods were derived for a different purpose 
in \cite{Heun1900} and in \cite{Fehlberg70}, where the third order method
was proposed as an embedding formula for the second order method (see also \cite{Butcher87}).
Shu \& Osher \cite{shu88a} derived a first framework for deducing higher order Runge--Kutta methods
with the total variation diminishing property and first identified the optimal explicit second and
third order methods with two and three stages, respectively. For these SSPRK(2,2) and SSPRK(3,3)
schemes an analysis of their stability, dissipativity and accuracy properties can be found in \cite{kochetal10a}.
For the SSPRK(3,2) method, a similar study is given in the Appendix section further below. Note
that SSPRK(2,2), SSPRK(3,2), and SSPRK(3,3) are the explicit schemes in the IMEX methods
(\ref{hig1}), (\ref{hig2}), and (\ref{hig3}), respectively.

Furthermore, we also consider some non-SSP methods from the literature, both explicit and IMEX,
for the sake of comparison, and also an asymptotically stable SSP IMEX scheme from \cite{parrus05}.

\subsubsection*{Implementation Issues}

To define the test problem, according to \cite{muthsamconvzones} we specify a hydrostatic
configuration which is unstable against convection. The simulation of a single semiconvective
layer requires the mean molecular weight to be linearly (and stably) stratified. As time evolves,
we expect convection to set in and mix the zone completely, although its development is inhibited
by the stable mean molecular weight gradient. A critical quantity in this process is hence the
buoyancy timescale and we return to this topic further below.

The simulations shown here have been performed on the Vienna Scientific Cluster, using 64 CPU cores
in parallel. The spatial resolution is 400$\times$400 grid points. Simulation time is measured in units of
\textit{sound crossing times (scrt)}. One scrt is defined as the time taking an acoustic wave to propagate
from the bottom to the top of the simulated box. In our simulations, $1\, \mathrm{scrt} = 5215.5\, s$ and
the simulation time is 200 scrt.



Restrictions on the time-step $\Delta t$ are imposed by heat diffusion $\tau_T$,
diffusion of the second species $\tau_c$, the viscosity $\tau_{\mathrm{visc}}$
and the velocity of the fluid $\tau_{\mathrm{fluid}}$,
\begin{equation}  \label{delta_t_explicit}
\Delta t = \min\{\tau_c, \tau_T, \tau_{\mathrm{visc}},
\tau_{\mathrm{fluid}}\},
\end{equation}
where
\begin{eqnarray*}
&&\tau_{c} = \frac{C_c}{\kappa_c}\, \min\{(\Delta x)^2, (\Delta y)^2\},\quad 
\tau_{T} = \frac{C_T}{\kappa_T}\, \min\{(\Delta x)^2, (\Delta y)^2\}, \\
&&\tau_{\mathrm{visc}} = \frac{C_\mathrm{visc}\, \rho}{\eta}\, \min\{(\Delta x)^2, (\Delta y)^2\},\quad 
\tau_{\mathrm{fluid}} = \frac{C_\mathrm{fluid}}{\max(\left|\mathbf{u}\right|)}\, \min\{\Delta x, \Delta y\},
\end{eqnarray*}
with (not necessarily equal) Courant numbers (CFL numbers) $C_c$, $C_T$, $C_\mathrm{visc}$, 
$C_\mathrm{fluid}$. This assumes that the time-step limitation due to sound waves has been 
removed by a fractional step approach as mentioned at the end of Section~\ref{sec:antares}.
Otherwise, $\tau_{\mathrm{fluid}}$ additionally
depends on the sound speed. We note that the source term in $F(y(t))$ in (\ref{mnseimex}),
which represents buoyancy forces acting on the flow, can be neglected in the limit where
$\max{\{\Delta x, \Delta y\}} \rightarrow 0$, since its contributions are of lower order (see
\cite{strikwerda04} for a discussion of the treatment of lower order terms in stability analyses).

Due to (\ref{mnseimex}), IMEX methods treat the terms $\nabla \cdot (\rho \kappa_c \nabla c)$
and $\nabla \cdot (K \nabla T)$ implicitly, so the restrictions $\tau_c$ and $\tau_T$ do not
have to hold. Since at least the first part of simulations of semiconvection is usually
dominated by diffusion processes and the Prandtl and Lewis numbers satisfy $\mathrm{Pr} < 1,\
\mathrm{Le}<1$ in a stellar context, the simulations are initially restricted by $\tau_T$. Hence, IMEX methods
can provide the desired computational advantage.

To enhance the stability and the efficiency of our methods we have implemented a heuristic
to adaptively select the time steps.
The most effective criterion for regulating the time steps turned out to be to monitor two-point instabilities
appearing in the conservative variables $(\rho, \rho c, \rho \vec{u}, Et)$. Due to the gravitational
force operating vertically
, such instabilities are prone to appear in the horizontal direction.

To detect the occurrence of such oscillations, for each variable we use the difference between two
grid cells to determine the sign of the corresponding gradient. In case of the density $\rho$  this reads

\begin{eqnarray*}
&& d_1 = \rho_{i,j-1} - \rho_{i,j-2}, \\
&& d_2 = \rho_{i,j} - \rho_{i,j-1}, \\
&& d_3 = \rho_{i,j+1} - \rho_{i,j}, \\
&& d_4 = \rho_{i,j+2} - \rho_{i,j+1},
\end{eqnarray*}

where $1 < i < n_x$ and $1 < j < n_y$, assuming the grid consists of $n_x \times n_y$ points.
If the sign pattern of $(d_1,d_2,d_3,d_4)$ corresponds to $(+,-,+,\pm)$, $(-,+,-,\pm)$, $(\pm, +,-,+,)$ or $(\pm,-,+,-)$,
we have located a two-point instability. Since the time-step is initially chosen to be that one required for
a fully explicit time integration method as in (\ref{delta_t_explicit}), such patterns are smoothed out rapidly, if present
in the initial condition. Consequently, their later occurrence is a good indicator for an instability developing because
of too large a time-step taken during the time integration.

The time-step control permits the occurrence of $n_y \cdot 0.1$ two-point instabilities for fixed $i$.
If this limit is exceeded, the time-step is repeated using a step-size decreased by a factor $\frac{2}{3}$.

To permit the system to readjust after reducing the time-step no modifications of $\Delta t$ are made for the next 15 time-steps
regardless of the number of two-point instabilities. If the number of oscillations still exceeds the given limit after
those 15 time steps, the time-step is again reduced.
If no or very few two-point oscillations are encountered for over 50 successive time-steps, the time-step is
augmented by a factor of $\frac{5}{4}$.

Alternatively to this heuristic control, we could also monitor the rate of change in the solution
to adjust the time-steps. However, the proper rate required to prevent the development of two-point
instabilities for the present application turned out to be too pessimistic to achieve CFL numbers
as high as discussed below. Furthermore, a simple control of the residual did not yield a satisfactory
behaviour. This leaves the heuristic time-step control as the most effective method for the IMEX based
time integration of our simulations.

\subsubsection*{Numerical Results with SSP Schemes}

\begin{table}
\begin{center}
\begin{tabular}{|r||r|r|r|r|r|}
\hline
\multicolumn{1}{|c||}{Method} & \multicolumn{1}{c|}{$\Delta t_{\mathrm{max}}$} & \multicolumn{1}{c|}{$\Delta t_{\mathrm{mean}}$} &
\multicolumn{1}{c|}{$\mathrm{CFL}_{\mathrm{max}}$} & \multicolumn{1}{c|}{$\mathrm{CFL}_{\mathrm{mean}}$} &
\multicolumn{1}{c|}{$\mathrm{CFL}_{\mathrm{start}}$} \\\hline\hline
\multicolumn{6}{|c|}{Singlelayer $\mathrm{Pr}=0.1$,\ $\mathrm{Le}=0.1$,\ $R_{\rho} = 1.1$,\ $\mathrm{Ra}^{*} = 160 000$} \\\hline\hline
SSPRK(2,2) & 3.71 s & 3.71 s & 0.2 &  0.2 & 0.2 \\\hline
SSPRK(3,2) & 9.31 s & 9.31 s & 0.5 &  0.5 & 0.5 \\\hline
IMEX SSP2(2,2,2) & 22.21 s& 11.56 s & 1.20 & 0.62 & 0.2 \\\hline
IMEX SSP2(2,2,2), $\gamma$=0.24 & 36.35 s & 19.44 s & 1.96 & 1.05 & 0.2 \\\hline
IMEX SSP2(2,2,2), $\gamma$=0.24 & 36.35 s & 19.43 s & 1.96 & 1.05 & 0.3 \\\hline
IMEX SSP2(2,2,2), $\gamma$=0.24 & 37.02 s & 19.45 s & 2.00 & 1.05 & 0.4 \\\hline
IMEX SSP2(2,2,2), $\gamma$=0.24 & 35.53 s & 19.47 s & 1.91 & 1.05 & 0.5 \\\hline
IMEX SSP2(3,3,2) & 74.52 s & 57.37 s & 4.02 & 3.09 & 0.4 \\\hline
IMEX SSP2(3,3,2) & 93.15 s & 57.16 s & 5.02 & 3.08 & 0.5 \\\hline
SSPRK(3,3) & 3.71 s & 3.71 s &  0.2 & 0.2 & 0.2  \\\hline
IMEX SSP3(3,3,3) & 15.14 s & 10.14 s & 0.82 & 0.55 & 0.2 \\\hline\hline
\multicolumn{6}{|c|}{Singlelayer $\mathrm{Pr}=0.5$,\ $\mathrm{Le}=0.1$,\ $R_{\rho} = 1.1$,\ $\mathrm{Ra}^{*} = 160 000$} \\\hline\hline
SSPRK(2,2) & 3.72 s &  3.72 s & 0.2 & 0.2 & 0.2  \\\hline
SSPRK(3,2) & 9.31 s & 9.31 s & 0.5 & 0.5 & 0.5 \\\hline
IMEX SSP2(2,2,2) & 23.14 s & 13.84 s & 1.24 & 0.74 & 0.2 \\\hline
IMEX SSP2(2,2,2), $\gamma$=0.24 & 23.13 s& 15.72 s& 1.24 & 0.85 & 0.2 \\\hline
IMEX SSP2(2,2,2), $\gamma$=0.24 & 22.76 s& 15.79 s& 1.22 & 0.85 & 0.3 \\\hline
IMEX SSP2(2,2,2), $\gamma$=0.24 & 20.32 s & 15.81 s & 1.09 & 0.85 & 0.4 \\\hline
IMEX SSP2(2,2,2), $\gamma$=0.24 & 22.81 s& 15.85 s& 1.23 & 0.84 & 0.5 \\\hline
IMEX SSP2(3,3,2) &40.65 s & 33.54 s & 2.19 & 1.80 & 0.4 \\\hline
IMEX SSP2(3,3,2) &40.65 s & 33.87 s & 2.19 & 1.82 & 0.5 \\\hline
SSPRK(3,3) & 3.72 s & 3.72 s & 0.2 & 0.2 & 0.2  \\\hline
IMEX SSP3(3,3,3) & 15.05 s & 9.70 s & 0.81 & 0.52 & 0.2 \\\hline
\end{tabular}
\caption{Comparisons of time steps and CFL-numbers over the first 80 scrt for the case where
$\mathrm{Pr}=0.1$ and over the entire 200 scrt for the case where $\mathrm{Pr}=0.5$ (see also text).} \label{NumericalResults1}
\end{center}
\end{table}

\begin{figure}[!ht]%
  \centering
    \subfigure[\sffamily $\mathrm{CFL}_{\mathrm{start}}=0.2$]{
    \scalebox{0.5}{\includegraphics[width=12cm,angle=270]{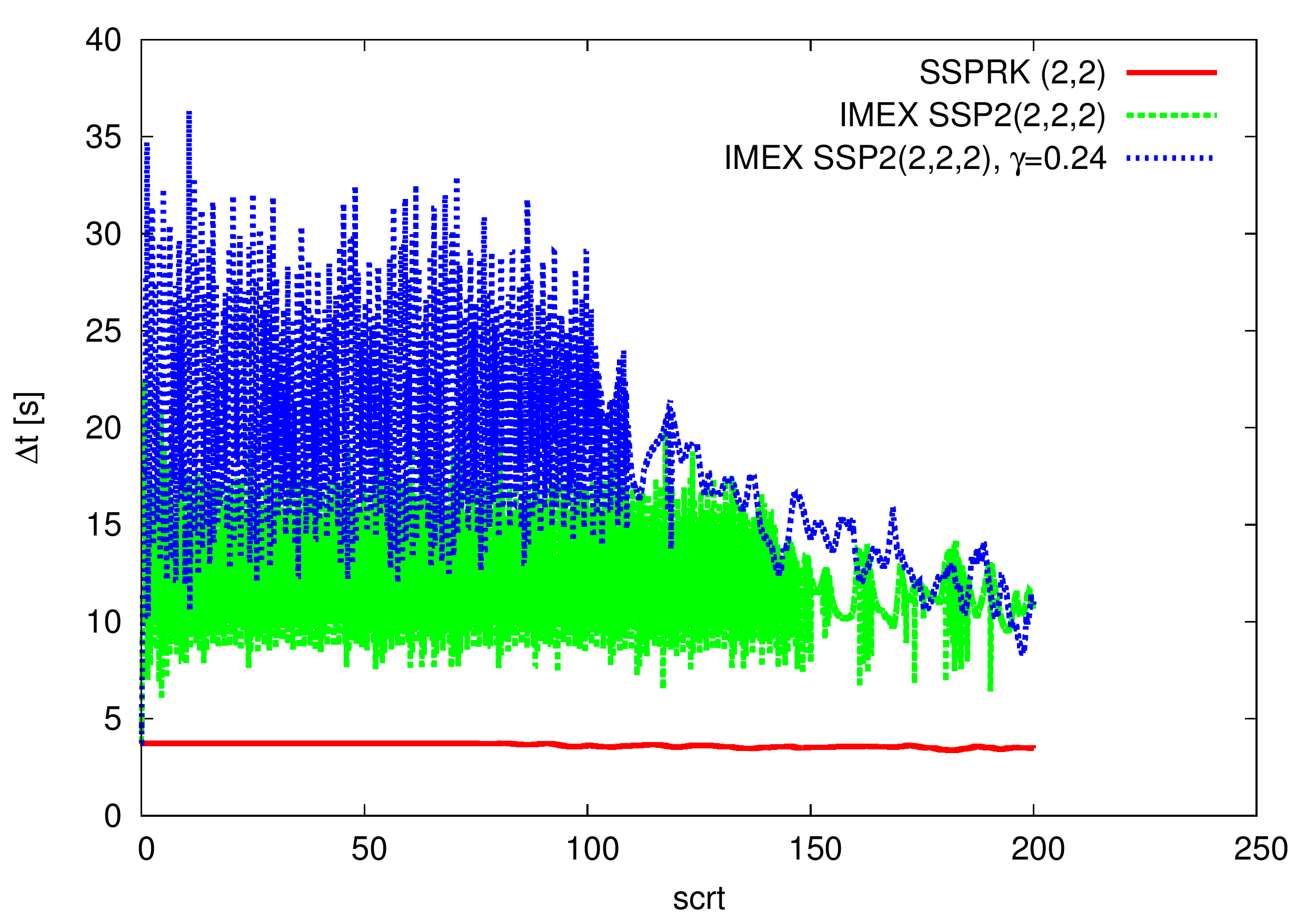}}}
    \label{fig:TauPr0-1ord2co02}
\hspace{0.2cm}
  \subfigure[\sffamily $\mathrm{CFL}_{\mathrm{start}}=0.5$]{%
    \label{fig:Inj16b_0-1}%
    \scalebox{0.5}{\includegraphics[width=12cm,angle=270]{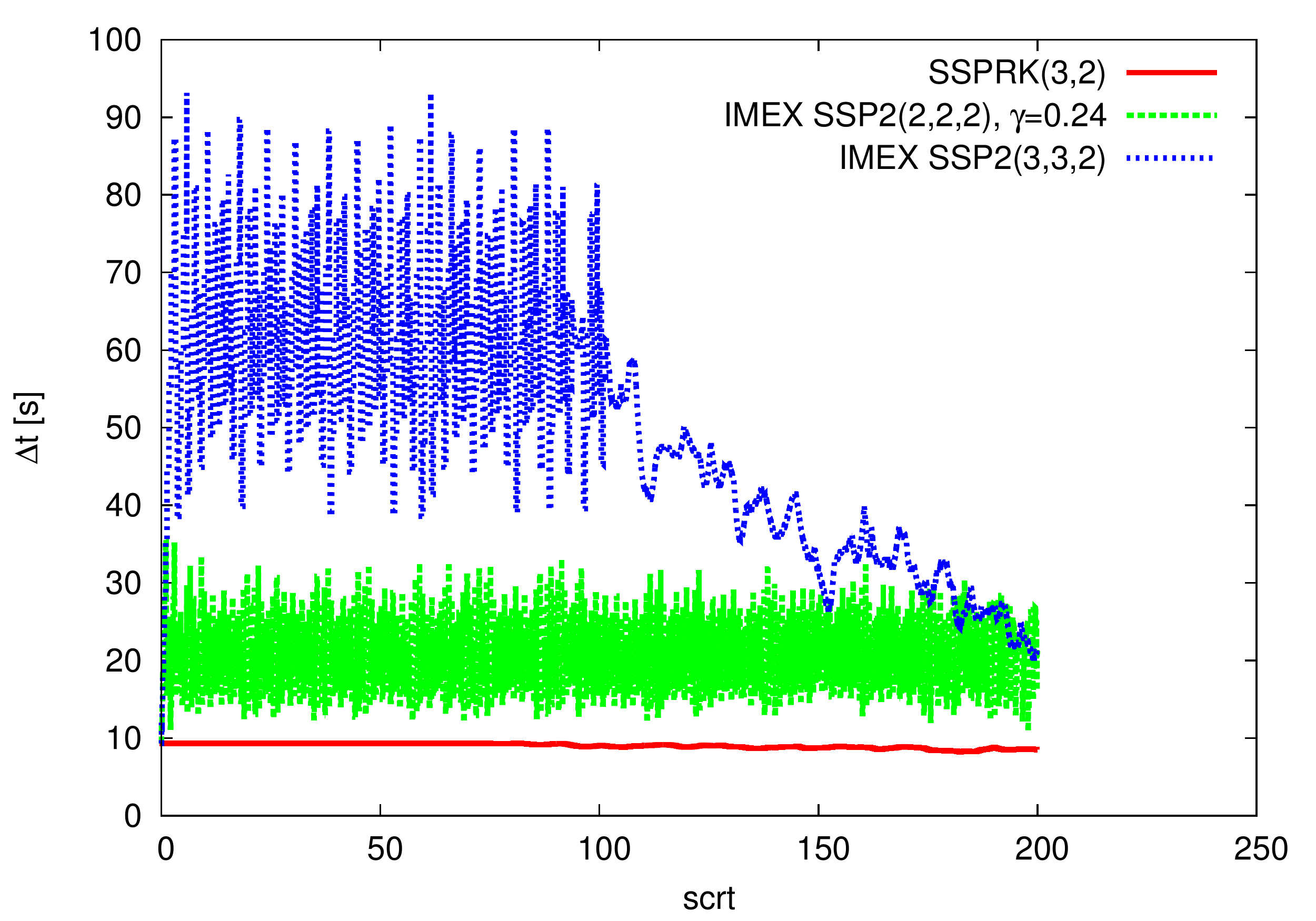}}}
 \label{fig:TauPr0-1ord2co05}%
\hspace{0.2cm}
  \subfigure[\sffamily $\mathrm{CFL}_{\mathrm{start}}=0.2$]{%
    \label{fig:TauPr0-1ord3co02}%
    \scalebox{0.5}{\includegraphics[width=12cm,angle=270]{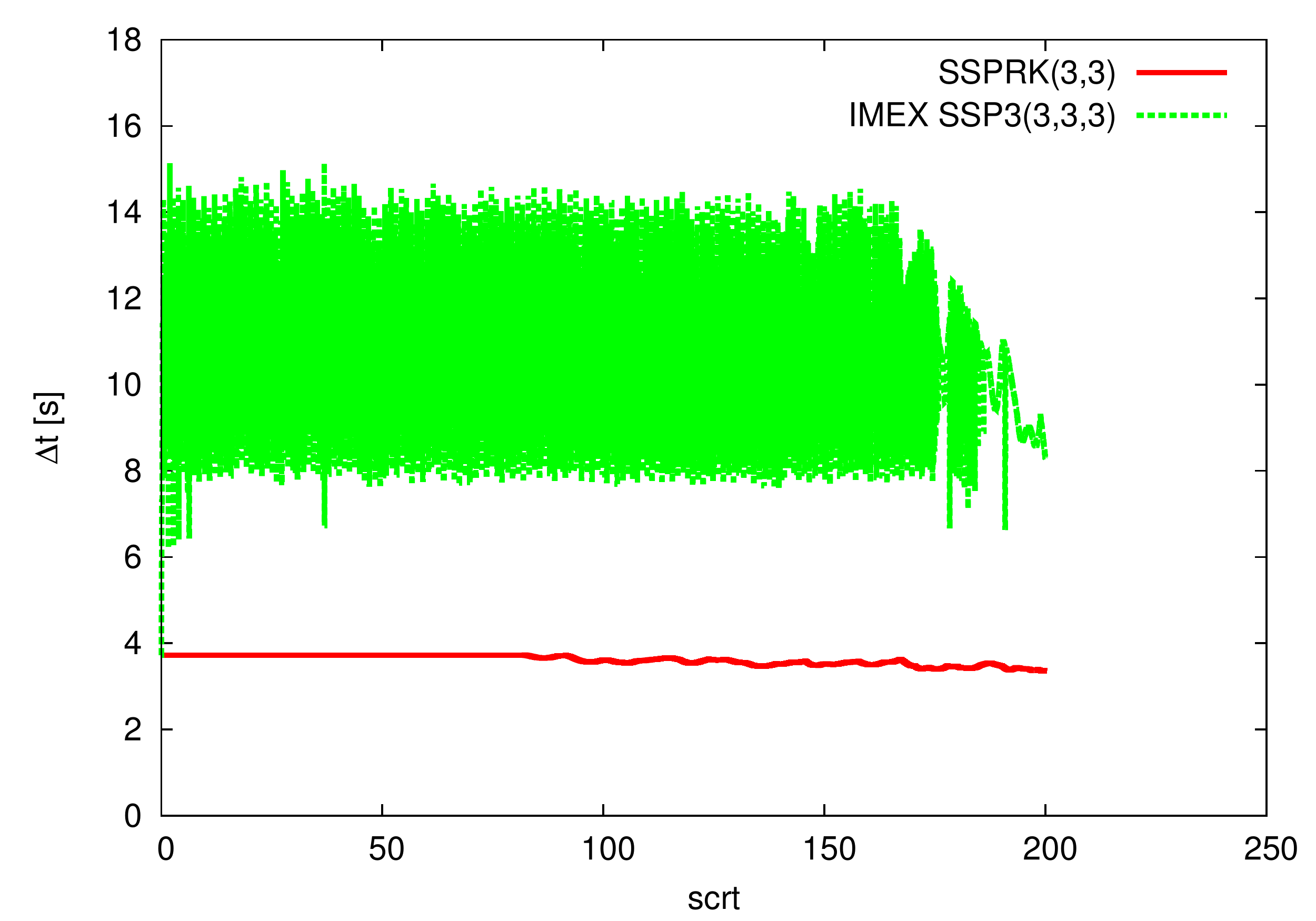}}}
    \\
\caption{Time-step evolution over 200 scrt in Simulation 1 (see text for definitions). Pictures (a) and (b) compare the
time-step $\Delta t$ of the second order schemes whereas (c) shows the evolution of $\Delta t$ using SSPRK(3,3) and
IMEX SSP3(3,3,3).}
\label{fig:simulation1}
\end{figure}

\begin{figure}[!ht]%
  \centering
    \subfigure[\sffamily $\mathrm{CFL}_{\mathrm{start}}=0.2$]{
    \scalebox{0.5}{\includegraphics[width=12cm,angle=270]{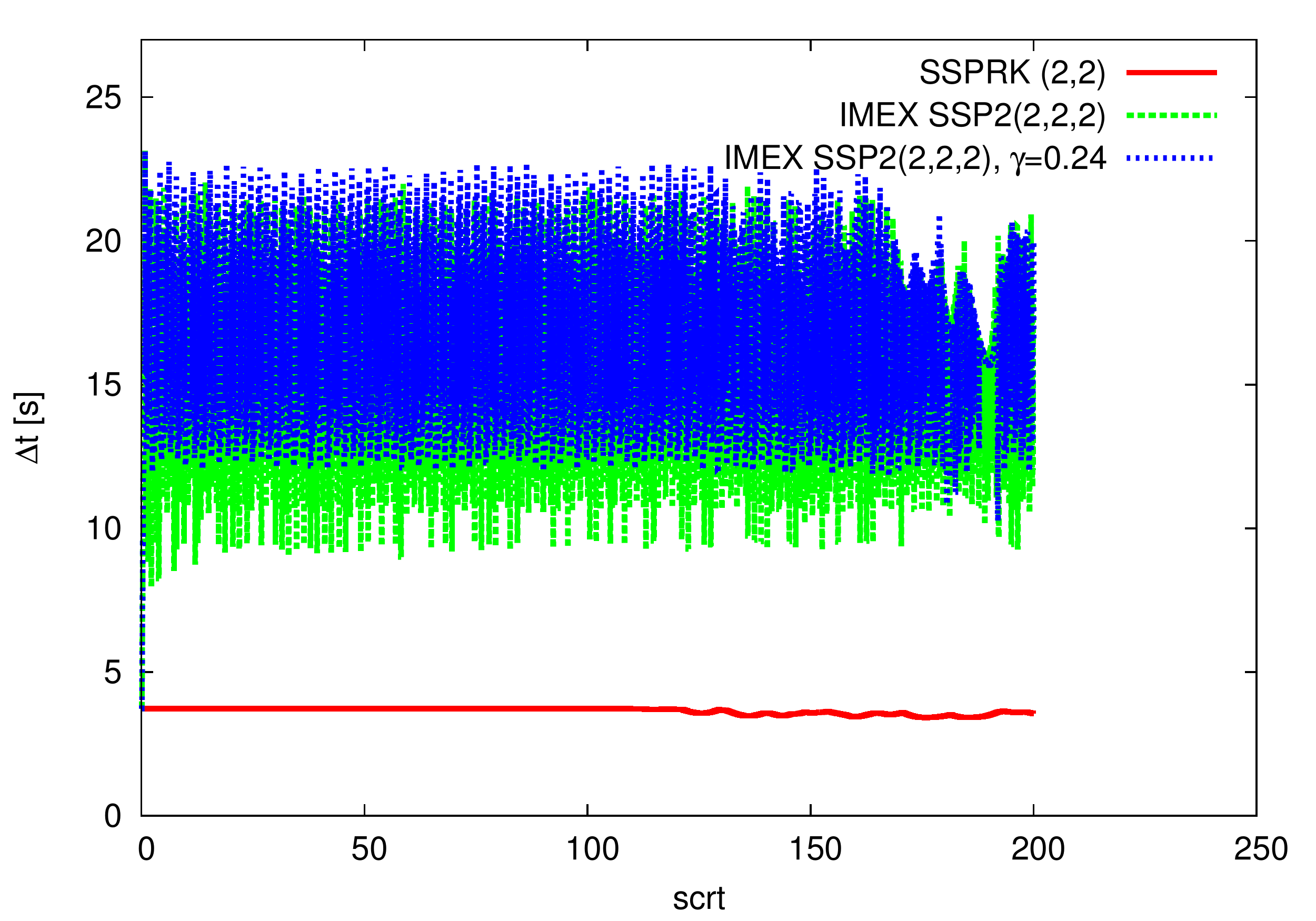}}}
    \label{fig:TauPr0-5ord2co02}
\hspace{0.2cm}
  \subfigure[\sffamily $\mathrm{CFL}_{\mathrm{start}}=0.5$]{%
    \label{fig:Inj16b_0-5}%
    \scalebox{0.5}{\includegraphics[width=12cm,angle=270]{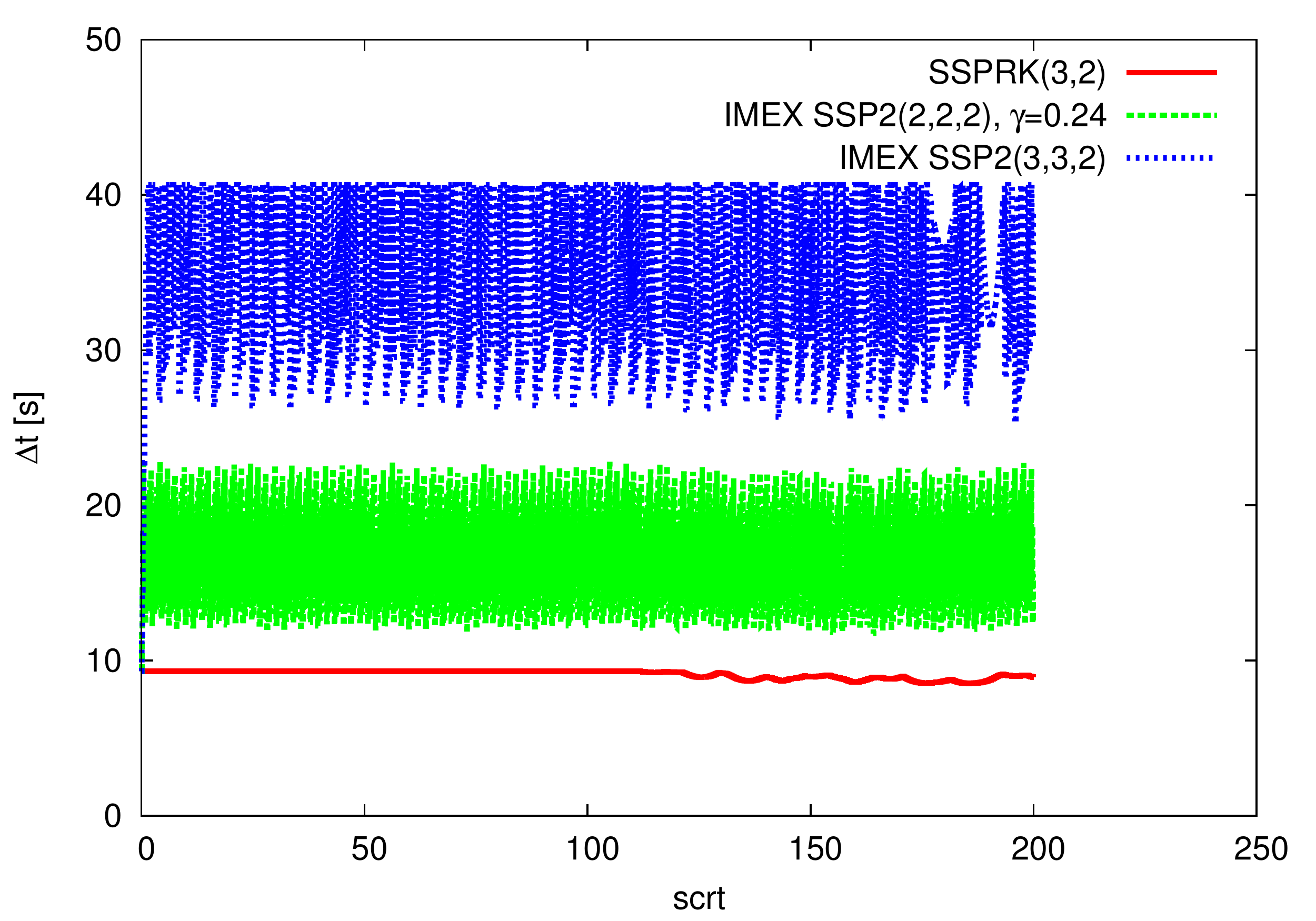}}}
 \label{fig:TauPr0-5ord2co05}%
\hspace{0.2cm}
  \subfigure[\sffamily $\mathrm{CFL}_{\mathrm{start}}=0.2$]{%
    \label{fig:TauPr0-5ord3co02}%
    \scalebox{0.5}{\includegraphics[width=12cm,angle=270]{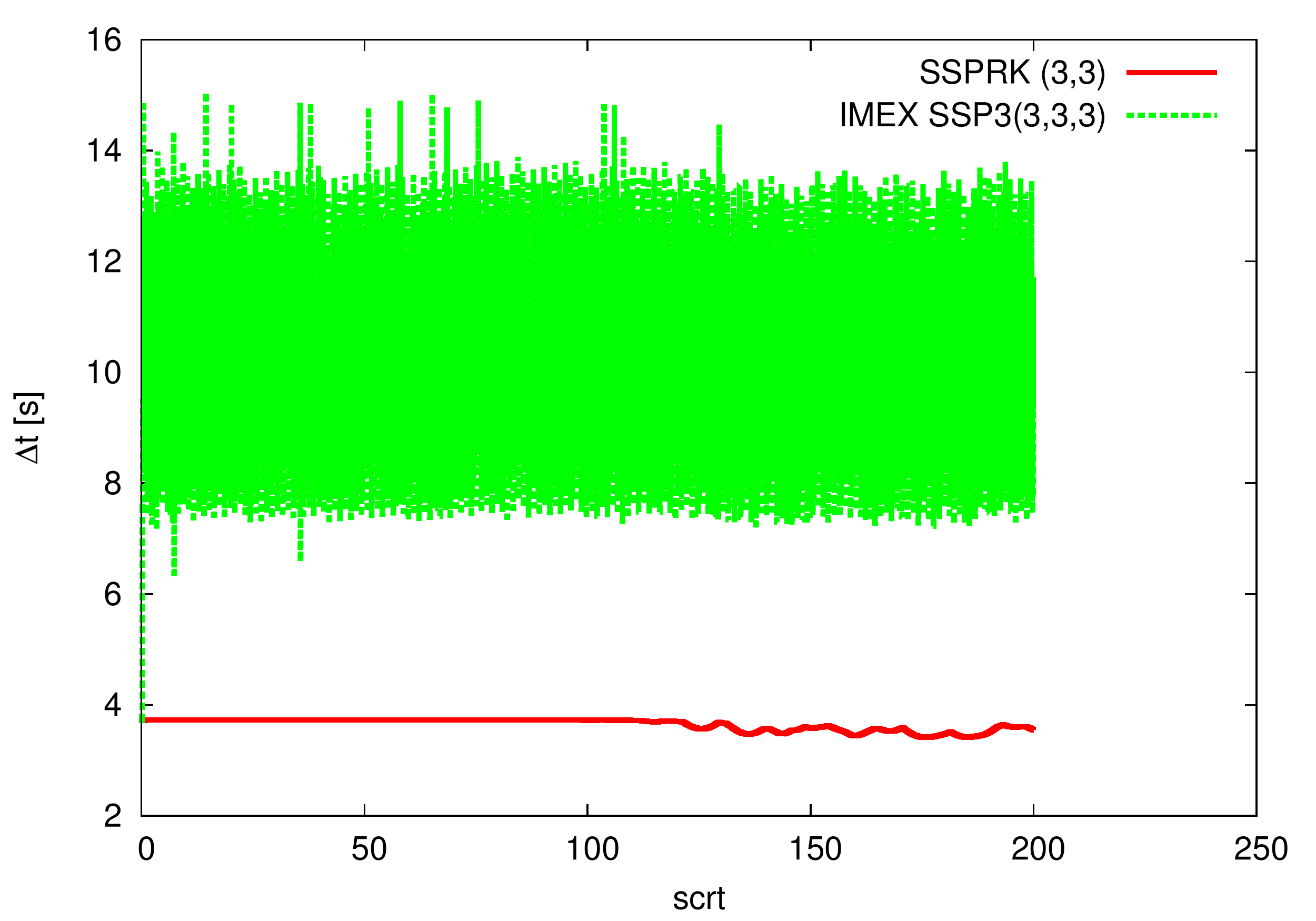}}}
    \\
\caption{Time-step evolution over 200 scrt in Simulation 2 (see text for definitions). Pictures (a) and (b) compare the time-step $\Delta t$
of the second order schemes whereas (c) shows the evolution of $\Delta t$ using SSPRK(3,3) and IMEX SSP3(3,3,3).}
\label{fig:simulation2}
\end{figure}

Tables \ref{NumericalResults1} and \ref{NumericalResults2} sum up the performance achieved with the presented Runge--Kutta schemes.
The evolution of the time steps is illustrated
graphically in Figures~\ref{fig:simulation1} and \ref{fig:simulation2} for the two simulation scenarios we have focused on.
Since the capability of the method is best judged in that part of the simulation where the fluid velocity is too small
to severely limit the time-step $\Delta t$, the CFL-numbers and time-steps listed in Table \ref{NumericalResults1}
have been measured in this regime. For the first test scenario this corresponds only to the first 80 scrt
(note the slope beginning after about 100 scrt in Figure~\ref{fig:simulation1}).

\begin{table}
\begin{center}
\begin{tabular}{|r||r|r|r|}
\hline
\multicolumn{1}{|c||}{Method} & \multicolumn{1}{c|}{$\mathrm{CFL}_{\mathrm{start}}$} & \multicolumn{1}{c|}{computation time} &
\multicolumn{1}{c|}{number of time steps} \\\hline\hline
\multicolumn{4}{|c|}{Singlelayer $\mathrm{Pr}=0.1$,\ $\mathrm{Le}=0.1$,\ $R_{\rho} = 1.1$,\ $\mathrm{Ra}^{*} = 160 000$} \\\hline\hline
SSPRK(2,2) & 0.2 & 7:17:32 & 287 788\\\hline
SSPRK(3,2) & 0.5 & 4:40:24 & 115 691 \\\hline
IMEX SSP2(2,2,2) & 0.2 & 6:04:54 & 92 492 \\\hline
IMEX SSP2(2,2,2), $\gamma$=0.24 & 0.2 & 4:41:15 & 63 586 \\\hline
IMEX SSP2(2,2,2), $\gamma$=0.24 & 0.3 & 4:18:59 & 56 741 \\\hline
IMEX SSP2(2,2,2), $\gamma$=0.24 & 0.4 & 4:10:26 & 54 015 \\\hline
IMEX SSP2(2,2,2), $\gamma$=0.24 & 0.5 & 4:12:38 & 53 941 \\\hline
IMEX SSP2(3,3,2) & 0.4 & 4:31:12 & 27 673 \\\hline
IMEX SSP2(3,3,2) & 0.5 & 4:27:46 & 24 266 \\\hline
SSPRK(3,3) & 0.2 & 10:55:40 & 288 607 \\\hline
IMEX SSP3(3,3,3) & 0.2 & 8:19:39 & 104 789 \\\hline\hline
\multicolumn{4}{|c|}{Singlelayer $\mathrm{Pr}=0.5$,\ $\mathrm{Le}=0.1$,\ $R_{\rho} = 1.1$,\ $\mathrm{Ra}^{*} = 160 000$} \\\hline\hline
SSPRK(2,2) & 0.2 & 7:01:44 & 286 011   \\\hline
SSPRK(3,2) & 0.5 & 4:34:35 & 114 409 \\\hline
IMEX SSP2(2,2,2) & 0.2 & 5:10:01 & 75 384 \\\hline
IMEX SSP2(2,2,2), $\gamma$=0.24 & 0.2  & 4:45:36  & 66 997 \\\hline
IMEX SSP2(2,2,2), $\gamma$=0.24 & 0.3 & 4:42:00 & 68 591 \\\hline
IMEX SSP2(2,2,2), $\gamma$=0.24 & 0.4 & 4:42:24 & 66 847  \\\hline
IMEX SSP2(2,2,2), $\gamma$=0.24 & 0.5 & 4:50:52 & 66 828 \\\hline
IMEX SSP2(3,3,2) & 0.4 & 4:43:27 & 31 225 \\\hline
IMEX SSP2(3,3,2) & 0.5 & 4:43:26 & 31 112 \\\hline
SSPRK(3,3) & 0.2 & 10:45:43 & 286 006  \\\hline
IMEX SSP3(3,3,3) & 0.2 & 7:27:01 &  108 852 \\\hline
\end{tabular}
\caption{Comparisons: computation times and overall number of time steps over 200 scrt.} \label{NumericalResults2}
\end{center}
\end{table}

In Table~\ref{NumericalResults1} we compare the largest possible time steps $\Delta t_{\max}$, the average time steps
$\Delta t_\mathrm{mean}$, the maximal and average CFL numbers resulting from the adaptive step selection, and the initial CFL numbers.
The tests were performed for Prandtl numbers $\mathrm{Pr}=0.1$ and $\mathrm{Pr}=0.5$ distinguishing Simulations 1 and 2, respectively,
a Lewis number $\mathrm{Le}=0.1$, $R_\rho=1.1$, and a modified Rayleigh number $\mathrm{Ra}^{*} = 160000$ related to the
Rayleigh number through $\mathrm{Ra}^{*} = \mathrm{Ra}\cdot \mathrm{Pr}$.

Comparing the performance of the second order schemes it is obvious that the
modification of $\gamma$ in IMEX SSP2(2,2,2) has a stunning effect on the stability of the scheme.
The positive definiteness of dissipation in this method proves most effective in suppressing oscillations,
permitting an average time-step and CFL-number up to two thirds higher than the original IMEX SSP2(2,2,2) method.

Table \ref{NumericalResults1} shows that IMEX SSP2(3,3,2) permits a time-step and CFL number more than twice as
high as IMEX SSP2(2,2,2) even with modified $\gamma$.
However, a comparison of the computation time given in Table~\ref{NumericalResults2} shows that this does not improve by
the same factor as the CFL-number, since as the time-step $\Delta t$ grows, the iterative solver for the generalized elliptic
problems requires more iterations to converge, resulting in an increase in computation time.
A comparison of the computation times shows that, although IMEX SSP2(3,3,2) permits impressively large time steps,
the need to solve three additional generalized Poisson problems related to the third stage takes its toll,
whence the method's performance is inferior to IMEX SSP2(2,2,2) with modified $\gamma$ and in case of Simulation 2,
is not even competitive to SSPRK(3,2), though it still performs better than the  SSPRK(2,2) scheme.
Note that for $\mathrm{Pr}=0.1$, the IMEX methods outperform even the best explicit method, while for the more
moderate $\mathrm{Pr}=0.5$, the best explicit integrator SSPRK(3,2) is slightly
more efficient whereas the classical methods noticeably lag behind.

Interestingly, the initial (preset) CFL number has a negligible influence on the actual CFL number reached
in the diffusive part of the simulation. However, as soon as the fluid velocity seriously restricts $\Delta t$,
a higher initial CFL number leads to a significantly larger average time-step and reduces the required computation time,
since its value is used to define the time-step restriction for the terms integrated with the explicit part of the IMEX scheme.

Once the time-step is limited by $\tau_{\mathrm{fluid}}$, the remaining gain in $\Delta t$ by the IMEX schemes in ANTARES
is essentially due to the optimization of the time-steps by the algorithm explained further above. Since $\Delta t$ is rather
small in that case the convergence of the generalized Poisson solver is fast enough to allow the IMEX schemes to lag only
slightly behind their explicit counterparts.

We point out that optimization of the solver for the generalized Poisson equation (\ref{eq:genellc})
has the potential for a further significant decrease of computation time required, both in the diffusive regime,
but also if the time-step is limited by $\tau_{\mathrm{fluid}}$. This is important since changing
the time integration method during a simulation run is not advisable if both the diffusive and the convective
phase should be interpreted in a consistent manner.

A comparison of SSPRK(3,3) and IMEX SSP3(3,3,3) also shows that the larger time-steps of the semi-implicit method
lead to a significant gain in computational efficiency in case of a third order method.

\begin{figure}[!t]%
  \centering
    \subfigure[\sffamily Simulation 1, first 40 scrt]{
    \scalebox{1.0}{\includegraphics[width=6.6cm,angle=0]{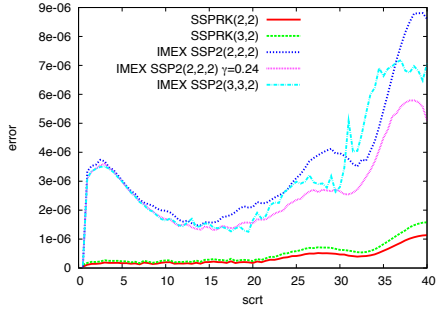}}}
    \label{fig:ErrorSim1_40scrt}
  \subfigure[\sffamily Simulation 1, entire run]{%
    \label{fig:ErrorSim1_200scrt_log}%
\hspace{0.2cm}   
    \scalebox{1.0}{\includegraphics[width=6.6cm,angle=0]{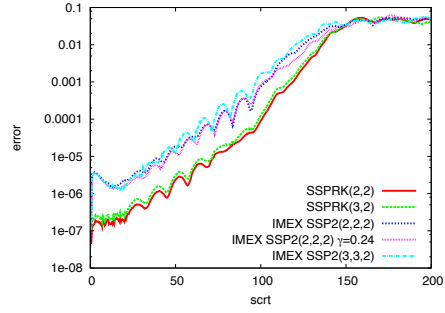}}}
    \subfigure[\sffamily Simulation 2, first 40 scrt]{
    \scalebox{1.0}{\includegraphics[width=6.6cm,angle=0]{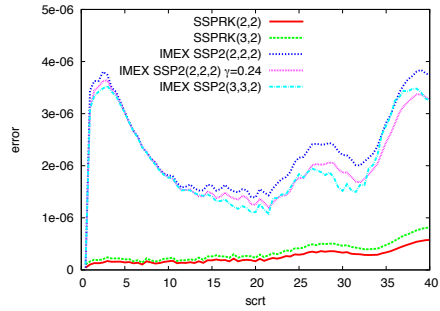}}}
    \label{fig:ErrorSim2_40scrt}
  \subfigure[\sffamily Simulation 2, entire run]{%
    \label{fig:ErrorSim2_200scrt_log}%
    \scalebox{1.0}{\includegraphics[width=6.6cm,angle=0]{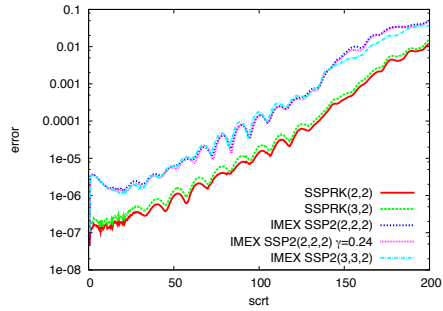}}}    
    \\
\caption{Time development of the root mean square difference per grid point of the mass ratio
             He / (H + He) between a reference solution 
             with the SSPRK(2,2) method and very small time-step (CFL-number 0.025) and various explicit 
             and IMEX SSP methods of second order. In the top row, picture (a) displays the first
             40 scrt on a linear scale for Simulation~1 and picture (b) shows the results for the entire run on 
             a logarithmic scale. In the bottom row, pictures (c) and (d)  show equivalent results
             for the case of Simulation~2.}
\label{fig:l2error}
\end{figure}

We have investigated the accuracy of the time integration with IMEX schemes by a comparison to a reference solution 
obtained with the SSPRK(2,2) method. The reference solution was computed on the same
spatial grid but with a time-step eight times smaller than that one mentioned in Tables~\ref{NumericalResults1} 
and~\ref{NumericalResults2} for this method, i.e.\ for a CFL-number of 0.025. 
Figure~\ref{fig:l2error} shows the root mean square difference of the mass ratio He / (H + He), obtained by
summation over all grid points and normalization relative to their number, between the numerical solutions 
of the second order SSP IMEX methods and the reference solution for each case. 
For the IMEX SSP2(2,2,2) method both the results for the standard choice
of $\gamma = 1-1/\sqrt{2}$ and the best performing value of 0.24 are displayed. We also show the normalized root mean
square differences between the reference solution and the SSP second order explicit methods computed with their
standard CFL number given in Table~\ref{NumericalResults1}.

For both Simulation~1 and~2 one can easily spot the initial increase of the error due to the growing time step for 
the IMEX methods induced by the automatic time step control. A plateau is reached once the time-step stabilizes around 
a typical mean value (cf.\ also Figures~\ref{fig:simulation1} and~\ref{fig:simulation2}). Note that the IMEX SSP2(2,2,2)
method with $\gamma=0.24$ has a smaller error than the original IMEX SSP2(2,2,2) method, as expected from the
error constant shown in Figure~\ref{fig:stabimex1modboundary}. The largest differences occur for the IMEX SSP2(3,3,2)
method which also has the largest mean time-step. The error constants of the different methods (see 
also the summarizing Table~\ref{summary} below) provide a rough measure for comparing simulation runs
with similar time-steps. 

However, a comparison for a given point in time has only limited meaning. One of the reasons is that the solution 
changes its nature as a function of time. Initial vertical oscillations are damped out (at least first few scrt), then the
velocity field slowly starts building up (visible after 15~scrt), followed by the formation of large scale
gravity waves (oscillatory behaviour of the error in the range between 25 and 100~scrt) until the waves start to break and
turbulence sets in. The importance of the contributions of each of the dynamical equations also changes during this
development. The whole process leads to an increasing error until a statistically stationary, turbulent 
state is reached. For Simulation~1 this occurs at around 150~scrt. For Simulation~2 this is just about to occur
shortly after the end of the simulation time of 200~scrt (the delay of the time development in this case is caused
by the larger viscosity of the fluid which follows from the choice of $\mathrm{Pr}$ and $\mathrm{Ra}^{*}$). 
In the turbulent state of the system spatial correlation is lost on very short timescales (a few scrt). Thus, also
the reference solution no longer has a meaning due to the chaotic behaviour of the solution. The
spread of the error and also the saturation value observed during this phase (picture (b) of Figure~\ref{fig:l2error}) 
is set by the Dirichlet vertical boundary conditions on the concentration $c$. The third order methods
IMEX SSP3(3,3,3) and SSPRK(3,3) behave analogously.

\begin{figure}[!t]%
  \centering
    \subfigure[\sffamily SSPRK(2,2)]{
    \scalebox{0.5}{\includegraphics[width=0.95\linewidth,angle=0]{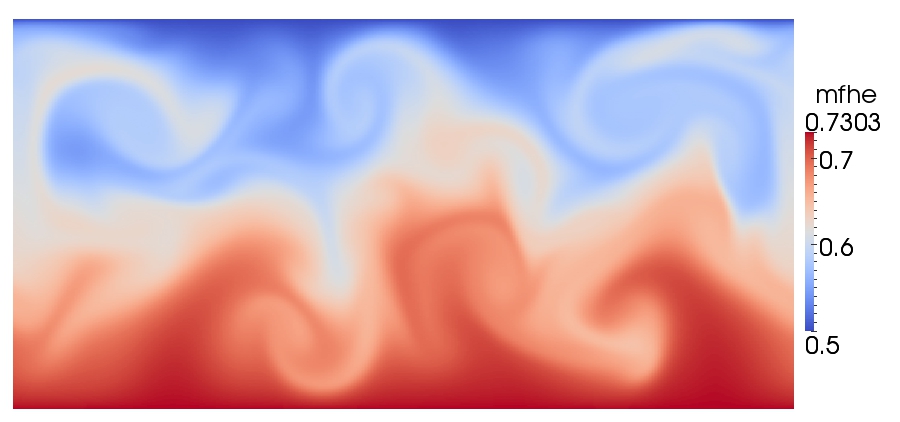}}}
    \label{fig:no1}%
  \subfigure[\sffamily IMEX SSP2(2,2,2)]{%
    \label{fig:no2}%
    \scalebox{0.5}{\includegraphics[width=0.95\linewidth,angle=0]{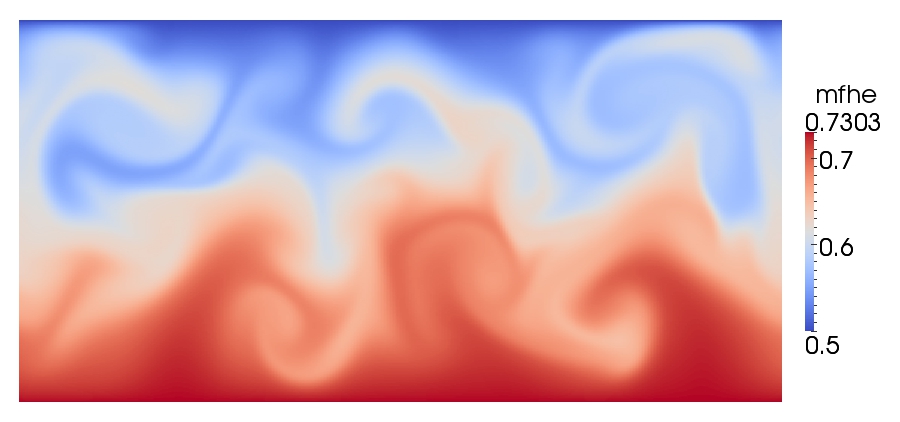}}}
  \subfigure[\sffamily SSPRK(3,2)]{%
    \label{fig:no3}%
    \scalebox{0.5}{\includegraphics[width=0.95\linewidth,angle=0]{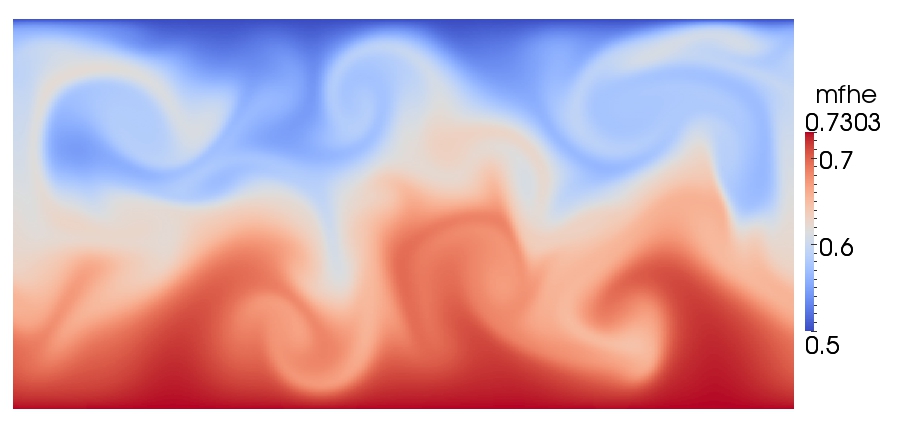}}}
  \subfigure[\sffamily IMEX SSP2(3,3,2)]{%
    \label{fig:no4}%
    \scalebox{0.5}{\includegraphics[width=0.95\linewidth,angle=0]{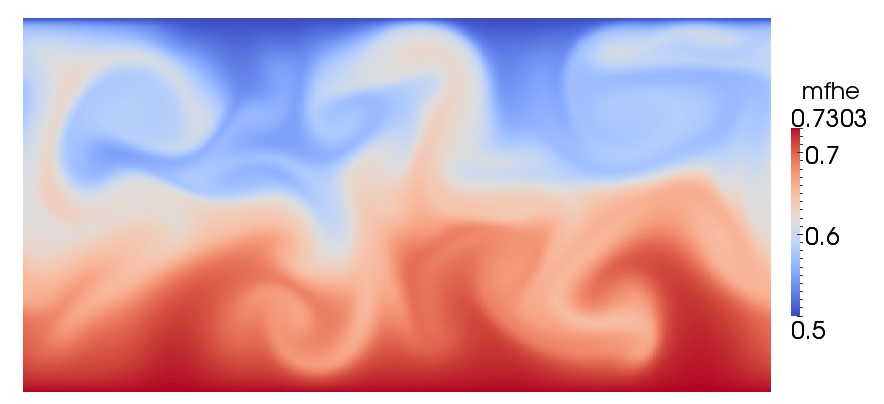}}}
  \subfigure[\sffamily SSPRK(3,3)]{%
    \label{fig:no5}%
    \scalebox{0.5}{\includegraphics[width=0.95\linewidth,angle=0]{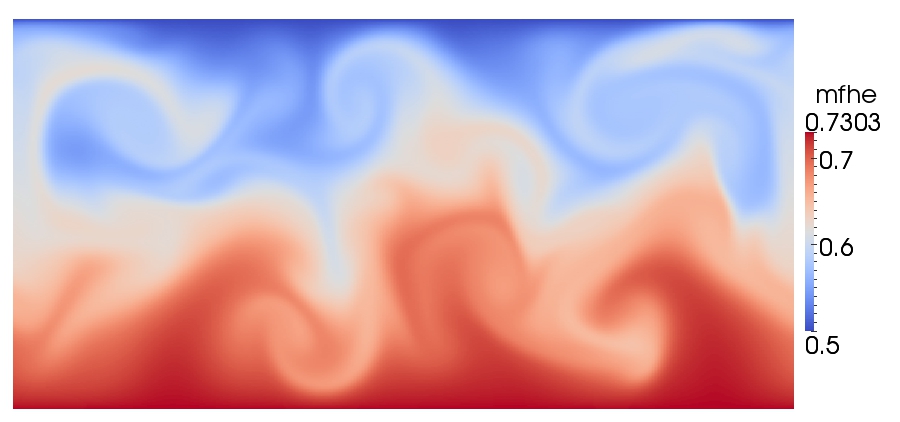}}}
  \subfigure[\sffamily IMEX SSP3(3,3,3)]{%
    \label{fig:no6}%
    \scalebox{0.5}{\includegraphics[width=0.95\linewidth,angle=0]{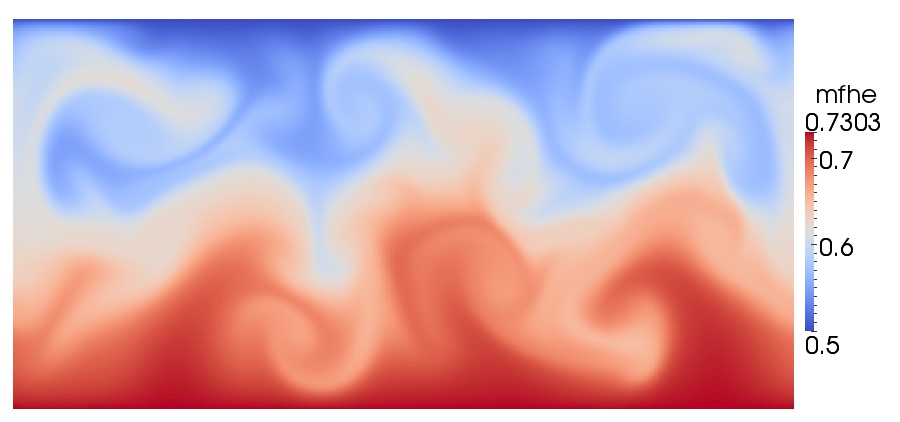}}}
\caption{Simulation~1 at $t = 125$ scrt.\label{figbla1}}
\end{figure}

\begin{figure}[!t]%
  \centering
    \subfigure[\sffamily SSPRK(2,2)]{
    \scalebox{0.5}{\includegraphics[width=0.95\linewidth,angle=0]{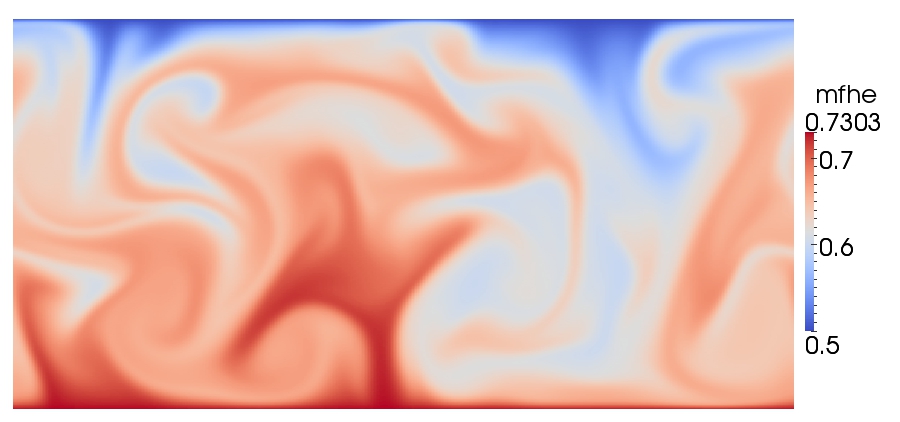}}}
    \label{fig:no7}%
  \subfigure[\sffamily IMEX SSP2(2,2,2)]{%
    \label{fig:no8}%
    \scalebox{0.5}{\includegraphics[width=0.95\linewidth,angle=0]{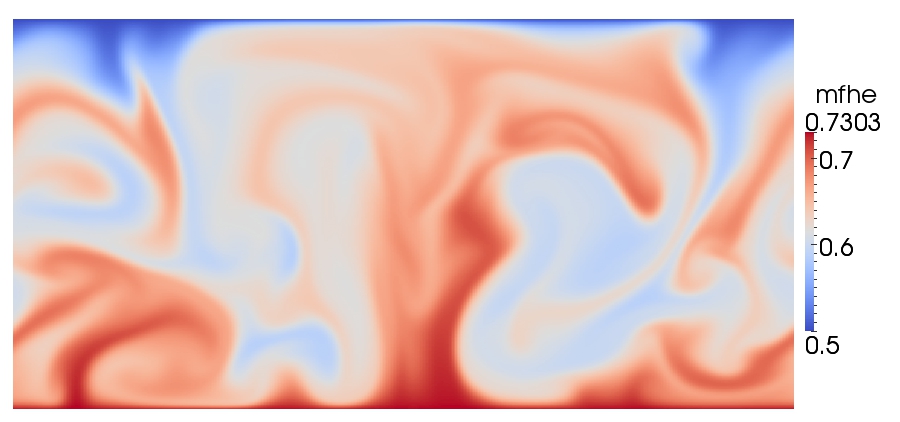}}}
  \subfigure[\sffamily SSPRK(3,2)]{%
    \label{fig:no9}%
    \scalebox{0.5}{\includegraphics[width=0.95\linewidth,angle=0]{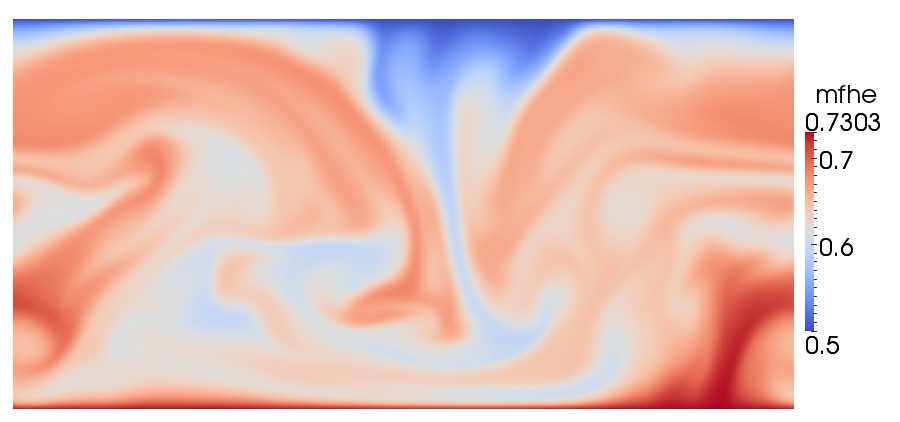}}}
  \subfigure[\sffamily IMEX SSP2(3,3,2)]{%
    \label{fig:no10}%
    \scalebox{0.5}{\includegraphics[width=0.95\linewidth,angle=0]{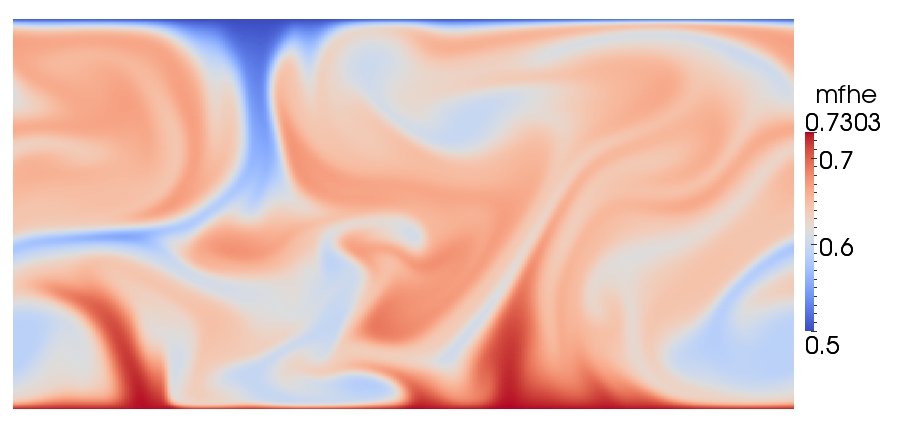}}}
  \subfigure[\sffamily SSPRK(3,3)]{%
    \label{fig:no11}%
    \scalebox{0.5}{\includegraphics[width=0.95\linewidth,angle=0]{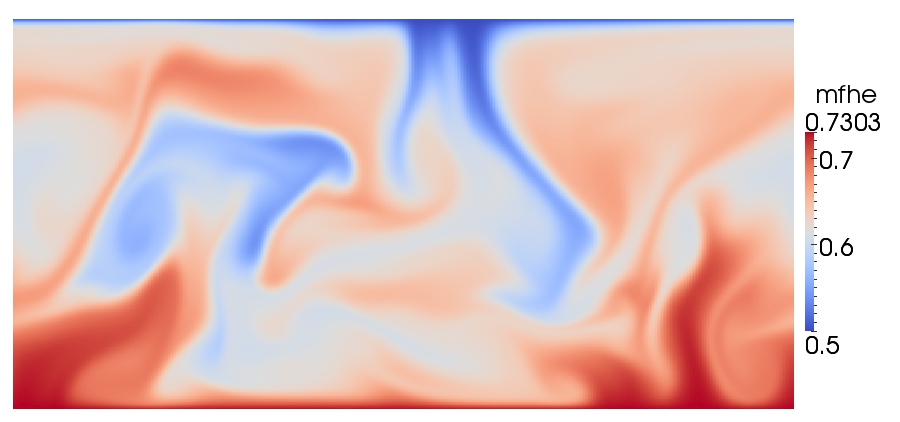}}}
  \subfigure[\sffamily IMEX SSP3(3,3,3)]{%
    \label{fig:no12}%
    \scalebox{0.5}{\includegraphics[width=0.95\linewidth,angle=0]{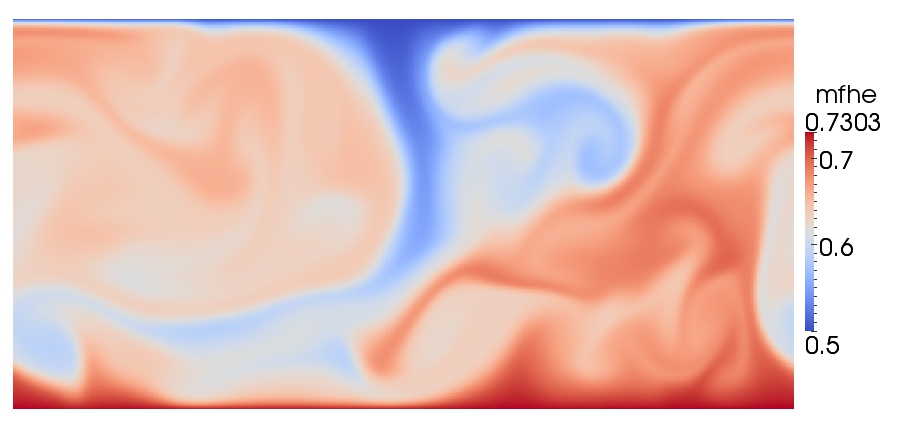}}}
\caption{Simulation results corresponding to those in Figure~\ref{figbla1}
at $t = 200$ scrt.\label{figbla2}}
\end{figure}

To further illustrate that the large time steps of the IMEX methods during the diffusive phase do not degrade 
the accuracy of the time development of the solution, we compare the simulation results obtained by the 
different time integrators in Figures~\ref{figbla1} and \ref{figbla2}. The pictures show the mass ratio of  
He vs.\ He+H at each spatial point at a given instant in time. Figure~\ref{figbla1} demonstrates that after 
the end of the  diffusive phase, just at the onset of turbulence, which occurs at $\sim 125$ scrt for this
problem, the results are quite comparable. The most visible differences can be found for the case of 
the simulation with the largest time-steps (picture (d)). It also shows that the root mean square errors 
displayed in Figure~\ref{fig:l2error} are negligible on a qualitative (and rough quantitative) level as long 
as they are smaller than about 10$^{-3}$. As expected, however, some time after the onset of turbulence
the solutions necessarily have drifted apart. This is demonstrated in Figure~\ref{figbla2}, where the
solutions already look different from each other.

\begin{figure}[!t]%
  \centering
    \subfigure[\sffamily high time resolution, t = 100 scrt]{
    \scalebox{0.5}{\includegraphics[width=0.95\linewidth,angle=0]{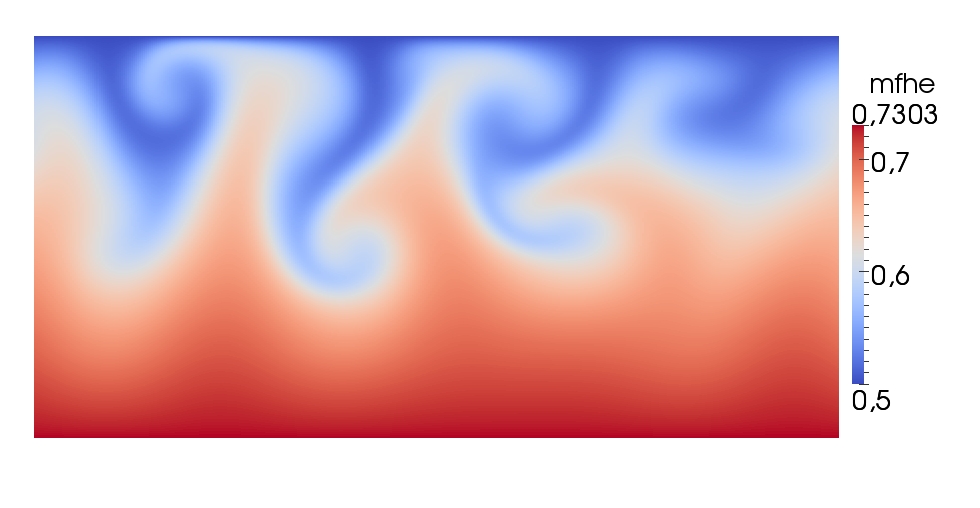}}}
    \label{fig:no13}%
  \subfigure[\sffamily high spatial resolution, t = 100 scrt]{%
    \label{fig:no14}%
    \scalebox{0.5}{\includegraphics[width=0.95\linewidth,angle=0]{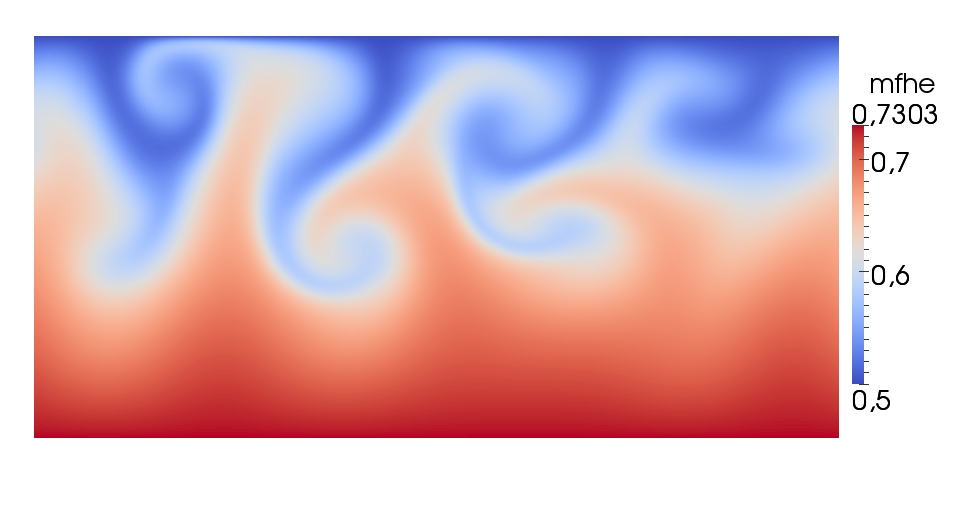}}}
  \subfigure[\sffamily high time resolution, t = 125 scrt]{%
    \label{fig:no15}%
    \scalebox{0.5}{\includegraphics[width=0.95\linewidth,angle=0]{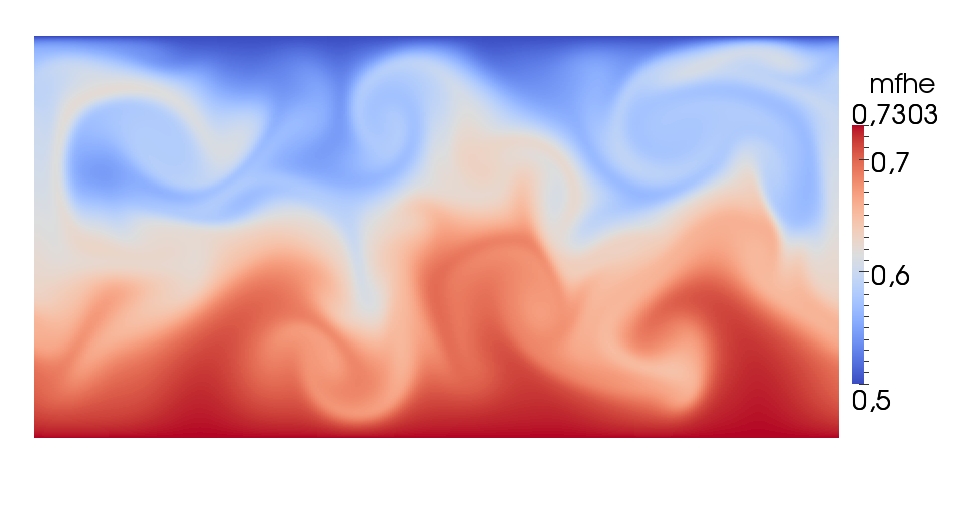}}}
  \subfigure[\sffamily high spatial resolution, t = 125 scrt]{%
    \label{fig:no16}%
    \scalebox{0.5}{\includegraphics[width=0.95\linewidth,angle=0]{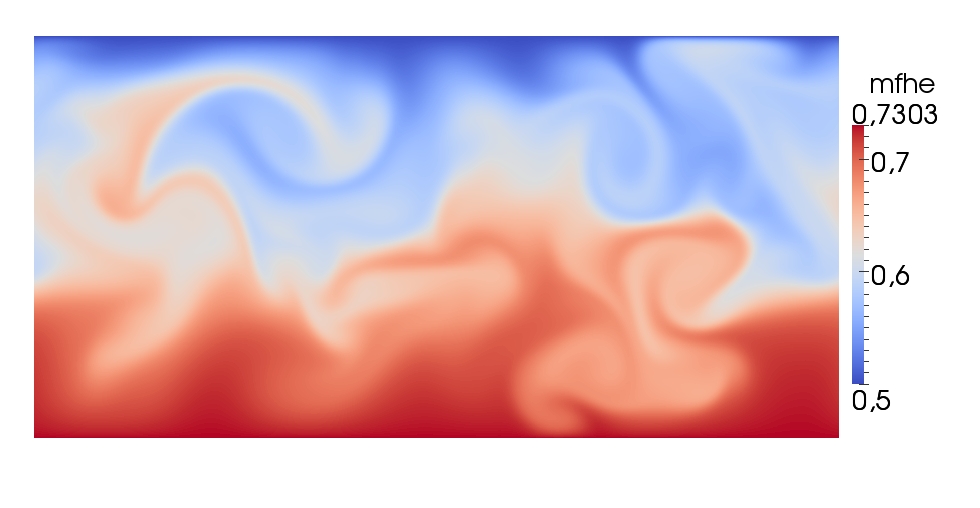}}}
  \subfigure[\sffamily high time resolution, t = 200 scrt]{%
    \label{fig:no17}%
    \scalebox{0.5}{\includegraphics[width=0.95\linewidth,angle=0]{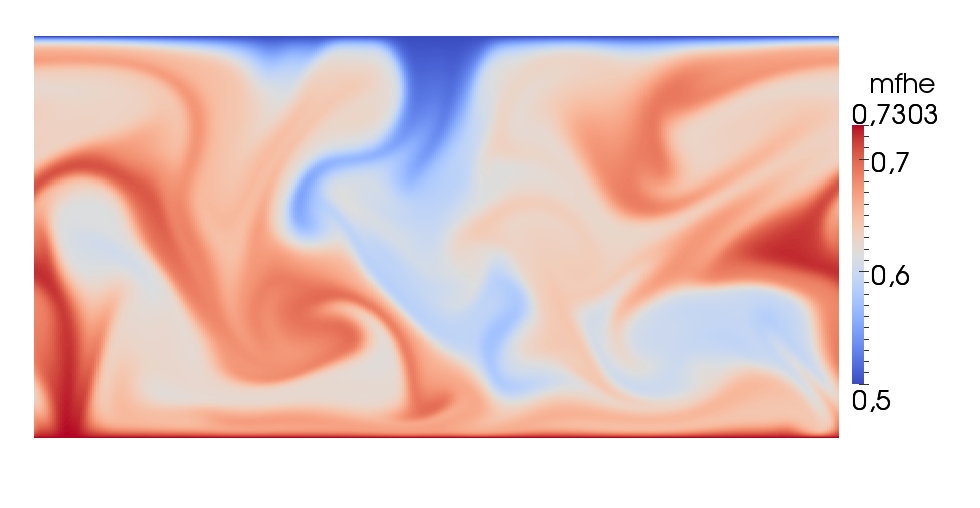}}}
  \subfigure[\sffamily high spatial resolution, t = 200 scrt]{%
    \label{fig:no18}%
    \scalebox{0.5}{\includegraphics[width=0.95\linewidth,angle=0]{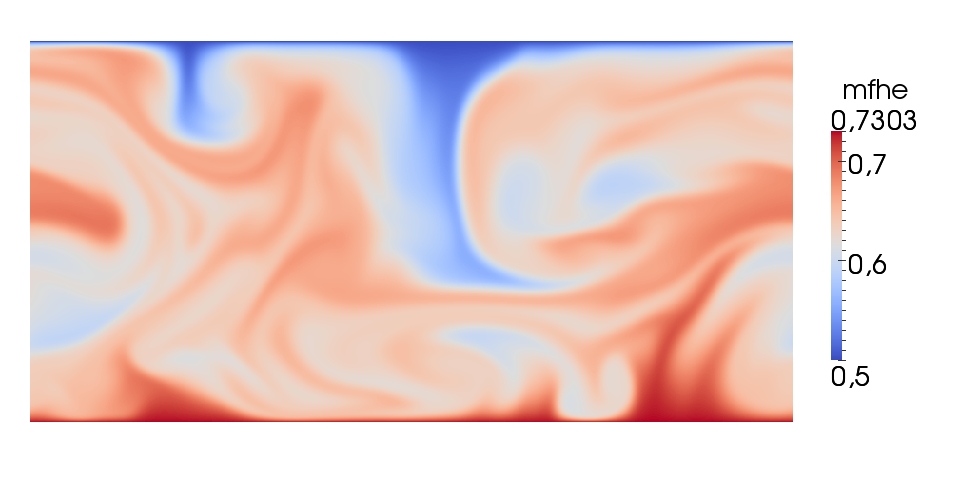}}}
\caption{Simulation~1 with SSPRK(2,2) time integration and high temporal resolution (left column) 
             as well as high spatial (and temporal) resolution (right column). The three different rows show 
             the results at different time $t$ in units of scrt.\label{figbla3}}
\end{figure}

Recalling Figure~\ref{fig:l2error} it is not surprising to find the largest differences
in Figure~\ref{figbla1} for the schemes having the largest time-steps during
the diffusive phase. Still, the large scale structures of the solution begin to diverge only
once the turbulent phase has been reached which in turn for each of the different
time integration methods occurs after about the same integration time $t$. Large
time-steps during the diffusive phase hence yield acceptable accuracy. Indeed,
the spatial resolution is more important than the temporal one. This can be 
demonstrated by using a high resolution grid of 800$\times$800 points. We have
performed such a reference run for the case of Simulation~1 with the SSPRK(2,2)
method for time integration. Note that the doubling of resolution leads to a four times
smaller time-step during the diffusive phase. 
Looking at the first row (pictures (a) and (b) of Figure~\ref{figbla3}) the differences
at 100~scrt, i.e.\ just after the onset of wave-breaking and the beginning of the 
turbulent phase, are still small: in the simulation with higher spatial resolution the
breaking tips are somewhat more pronounced and the contrasts sharper. This
has changed at 125~scrt shown in the second row of Figure~\ref{figbla3}.
The common initial condition may still be inferred, but the simulations have notably
evolved away from each other. Clearly, the spatial resolution is much more
important than the influence of the time steps and the time integration method
chosen, since picture (c) of Figure~\ref{figbla3} is essentially indistinguishable from its counterpart
with standard time resolution, picture (a) of Figure~\ref{figbla1}. Furthermore, at 100~scrt 
the simulation using the IMEX2(3,3,2) method for time integration, which has the 
largest time-step during the diffusive phase, is nearly indistinguishable from the
SSPRK(2,2) run with high time resolution shown in picture (a) of Figure~\ref{figbla3}
(the IMEX results are not shown here for the case of 100~scrt, since the differences to the
latter picture are very difficult to spot). 

We conclude that the spatial resolution is indeed more important than high temporal
resolution and large time steps during the diffusive phase are clearly tolerable for simulations
of astrophysical convective flows, if they do not affect stability. 
Resolutions of 800 grid points per spatial direction in 3D are usually not 
affordable anyway and sometimes even large parametric studies in 2D may still be too
expensive (for instance, for the case of semi-convection and purely explicit time integration
methods). We note that the necessity of a resolution of 400 points, which has been used for 
most of the simulation runs shown here, was calculated following \cite{zaussithesis} and
\cite{zaussispruit11} where in turn the physical arguments of \cite{spruit92} had been used
to estimate the thickness of solute and thermal boundary layers in semi-convection and 
the applicability of this approach to the parameter range we are interested in had been confirmed.
Thus, for the parameters of the more demanding one of our models, Simulation~1, we concluded
the smallest structures of interest, the solutal boundary layers, to span 6 grid points, if
the whole box is discretized by 400 points in each direction. Examples for them are the
top and bottom boundary layers which can easily be seen in Figures~\ref{figbla2} 
and~\ref{figbla3}. Indeed, in the simulation with a high spatial resolution of 800 points in each
direction these layers are hardly any thinner than in the case of 400 points (cf.\ the bottom row
of Figure~\ref{figbla3} which shows the simulations developed well into their final, turbulent state). 

The duration of the diffusive phase is determined by the stability of the stratification,
parametrized through $R_{\rho}$, and related to the buoyancy term, the second term
on the right-hand side of (\ref{mnseimex}) and thus also the second term of $F(y(t))$
in the same equation. Timescales related to this term can be computed for a variety
of physical problems. They include the reciprocal of the growth rate of small density
perturbations when a fluid of higher density $\rho_2$ is layered above fluid of lower density $\rho_1$
(this situation is denoted as \emph{Rayleigh-Taylor instability}, see \cite{Chandra61}). In this case the 
growth rate follows the dispersion relation $\omega^2 = g k (\rho_2-\rho_1) / (\rho_2+\rho_1)$ for a local
gravitational acceleration $g$ and a perturbation with wave number $k$. Its magnitude can be bounded
by the simple relation $\omega^2 = g k$. Similar dispersion relations are found for
gravity waves and growth rates of convective instability (see the classical paper
\cite{cowling41} and also \cite{kipweig94} for a summary). As implied by $\omega^2 = g k$
and the form of (\ref{mnseimex}), one can estimate a buoyancy time-step restriction
by $t_{\mathrm{buoy}} = \min\{(\Delta x)^{1/2}\} / g^{1/2}$ (see also \cite{Kang00}).
For the present simulations we find $t_{\mathrm{buoy}} \sim 360$~s, i.e.\ about 0.069~scrt.
Note that this is about four times larger than the largest time-step reported in Table~\ref{NumericalResults1}.
Though irrelevant in the asymptotic limit and not a constraining quantity here, one might still consider
this term to be integrated implicitly in other cases. However, if the excitation
and breaking of waves is important to describe the onset of the turbulent convective flow, as in the present case, damping
of such waves by an implicit time integration may be undesired. Hence, explicit time integration
of the buoyancy term could be preferred for physical reasons, even if the time-step were actually
constrained by such a splitting for the time integration.

\subsubsection*{Numerical Results with non-SSP Schemes}

From Tables~\ref{NumericalResults1} and~\ref{NumericalResults2} one can readily see that during
the diffusion dominated phase the SSP IMEX schemes achieve time-steps which exceed the
region of absolute monotonicity ensured by (\ref{hig1modkra2}) for IMEX SSP2(2,2,2) and their
counterparts given after (\ref{hig2}) and (\ref{hig3}) for IMEX SSP2(3,3,2) and IMEX SSP3(3,3,3),
respectively. One might hence question whether the property of absolute monotonicity is really
necessary for the time integration of the numerical simulations we have considered above. To show 
that this property is indeed required we have performed several test runs for the case of Simulation~1
with time integration schemes which do not diminish the total--variation norm.

The first candidate we have investigated is the ARK3(2)4L[2]SA scheme proposed in
\cite{kencar03}. This is an L-stable, stiffly accurate third order, additive Runge--Kutta method
with four stages. With its choice of coefficients it belongs to the group of IMEX methods. However,
neither its explicit nor its implicit part are strong--stability--preserving (the Butcher arrays feature
negative coefficients, cf.\ Theorem~4.2 in \cite{kraaijevanger91}), hence also the entire method 
does not fulfill the criteria of Theorem~\ref{am1}
on absolute monotonicity. If the SSP--property were of no importance, this method should be
quite robust. However, it turns out that this is not the case when we use it to integrate Simulation~1
in time. Taking $\mathrm{CFL}_{\mathrm{start}}$ to be 0.2 as for the other IMEX methods presented in
Table~\ref{NumericalResults2}, the time integration with ARK3(2)4L[2]SA  crashes after just 78 
time-steps. Indeed, $\mathrm{CFL}_{\mathrm{start}}$ has to be lowered to 0.1 to successfully launch
the simulation. But over the first 10 scrt $\mathrm{CFL}_{\mathrm{max}}$ is found to never exceed
$\sim 0.2$ and the average $\mathrm{CFL}_{\mathrm{mean}}$ is only $\sim 0.15$. In conclusion the
size of the time-steps achieved with this method have been found to not exceed those achieved with the 
explicit, three-stage, third-order SSPRK(3,3) method of \cite{shu88a}. Compared to IMEX SSP2(3,3,2),
the maximum and mean CFL numbers are 25 and 20 times smaller, respectively. We have thus given
up this simulation run after 10~scrt: evidently, the ARK3(2)4L[2]SA method is not efficient for the
kind of problems we are interested in. The implicit, L-stable and stiffly accurate nature of this scheme is
not sufficient to provide any advantages on its own during the diffusive phase of the simulation.

To further investigate the importance of the strong--stability--preserving property 
we selected the classical, explicit, third order, three--stage Runge--Kutta method
first proposed in \cite{Heun1900} and known as Heun's third order method.
This is a non-SSP scheme since not all of its stages are
used in the final integration which yields $y_\mathrm{new}$, as pointed out in \cite{kraaijevanger91},
where it was used to illustrate the growth of solutions measured in standard norms for
both parabolic and hyperbolic problems in situations where the exact solution is not growing in these
norms.  By comparison the SSPRK(3,3) scheme was found to not exhibit such growth. When
we apply Heun's third order scheme to integrate Simulation~1, we achieve a stable simulation over the
entire extent of 200~scrt with the same average time-step and with the same Courant number as 
for the SSPRK(3,3) scheme. However, during the diffusive phase in Simulation~1 the solution itself 
is slowly growing in time while a velocity field is being built up by the convective instability and at the 
grid scale the dissipation is provided by the parabolic terms during the entire simulation. 

The semi-convection problem discussed here is a rather benign example for numerical simulations in 
astrophysical applications. If we apply the same scheme to a simulation of solar surface convection 
as in \cite{Muthsam2010}, i.e.\ for a case of $219 \times 159$ grid points and a standard, 
moderately low resolution of $18.57 \times 40$~km$^2$, and a standard choice for
the microphysics with non-grey radiative transfer for the calculation of $Q_\mathrm{rad}$,
the differences in stability become apparent. For this physical problem the term
representing viscous dissipation in the momentum equation does not provide sufficient 
dissipation on the grid-scale at any resolution achievable in the foreseeable future and
the dissipation properties of the temporal and spatial discretization of the advection operator
become important (the term implicit large eddy simulation is used in such cases). While the 
SSPRK(3,3) method and also the SSPRK(2,2) method, used together with the spatial 
discretization of \cite{shu97}, have no problems in completing a simulation
of 20~scrt with a CFL number of 0.25, Heun's third order method leads to a crash after
9.2~scrt.\footnote{We would like to thank H.~Grimm--Strele for performing the 2D solar
convection simulations with ANTARES to test the stability properties of Heun's non-SSP explicit 
third order Runge--Kutta scheme.} 
Such failures are usually caused by numerically induced fluctuations in the low density and temperature
region of the simulation box which result in negative or at least unphysically small values, whence
they fall outside the tabulated region of microphysical properties. We consider this finding as sufficient to
exclude non-SSP methods from being recommendable for numerical simulations of stellar convection,
since also here we have chosen a rather benign test case within its class. Numerical simulations of stellar
surface convection in white dwarfs, A-type stars, or Cepheids reach far more extreme conditions with up
to three times higher, super-sonic Mach numbers, density contrasts around shock fronts higher by an order
of magnitude and more, and for the case of A-stars and Cepheids, at lower effective resolution because of 
four to ten times steeper gradients and limited computational resources. 

We finish our study of non-SSP Runge--Kutta schemes with a third case: an IMEX Runge--Kutta method
where both the explicit and the implicit part are strong--stability--preserving, but the combined scheme
does not fulfill the criteria for absolute monotonicity in the sense of Theorem~\ref{am1}. The scheme
was proposed in \cite{parrus05} (Table IV, page 139) and indeed the IMEX SSP2(3,3,2) scheme
(\ref{hig2}) is a modification of that scheme proposed in \cite{higueras06} to obtain a nontrivial region
of absolute monotonicity. The original scheme of \cite{parrus05} differs from that one only by having
the entries $\{1/4,0,0\}$ and $\{0,1/4,0\}$ in the first two rows of the Butcher array of the implicit scheme
instead of $\{1/5,0,0\}$ and $\{1/10,1/5,0\}$, respectively. This scheme performs quite well. Indeed,
during the diffusion dominated phase of Simulation~1 its time-steps are even 6\% to 8\% larger than
those achieved with IMEX SSP2(3,3,2), while during the turbulent phase they fall back to at most the
size achieved by the scheme (\ref{hig2}). Once more, however, we recall that astrophysical simulations
often have to deal with limited resolution at least in part of the simulation domain. To reproduce such
a case we have run Simulation~1 with both (\ref{hig2}) and the original scheme of \cite{parrus05}
for the case of only $100 \times 100$ grid points while leaving everything else unchanged. In this
case the boundary layer due to concentration $c$ is represented vertically only by one to two grid
points (see the discussion on resolution and reference solutions further above). While
during the diffusive phase no major differences become apparent, the behaviour found for the
advection--dominated, turbulent phase was discrepant: whereas nothing suspicious occurred for
the scheme (\ref{hig2}), the time-step of the scheme of  \cite{parrus05} dropped to arbitrarily
small values after $\sim 113$~scrt, which indicates the occurrence of two-point instabilities,
and the simulation had to be terminated. 

We conclude that only Runge--Kutta methods which are strong--stability--preserving  have the necessary
prerequisites for stable time integrations of astrophysical convection simulations. If, in addition,
the time-step is limited by diffusion processes, this limitation can be overcome by IMEX methods
provided their explicit and implicit parts are strong--stability--preserving. To ensure stability also
in cases of low resolution IMEX SSP methods should also have a nontrivial region of absolute
monotonicity as required by Theorem~\ref{am1}.

\section{Conclusions and Outlook}\mlabel{conclusions}

In this paper we have given an extensive discussion of the mathematical properties
and practical usefulness of total--variation--diminishing implicit--explicit
Runge--Kutta methods for the time integration of advection--diffusion
equations arising in the simulation of double--diffusive convection in
astrophysics. In this section, we summarize the results obtained in
Sections~\ref{imex-candidates} and \ref{num} (stability, dissipativity, accuracy and efficiency),
and give a brief outlook on future developments.

\begin{figure}[h!t]%
  \centering
    \subfigure[\sffamily SSP$2(2,2,2)$ (\ref{hig1})]{
    \scalebox{0.25}{\includegraphics[width=12cm,angle=270]{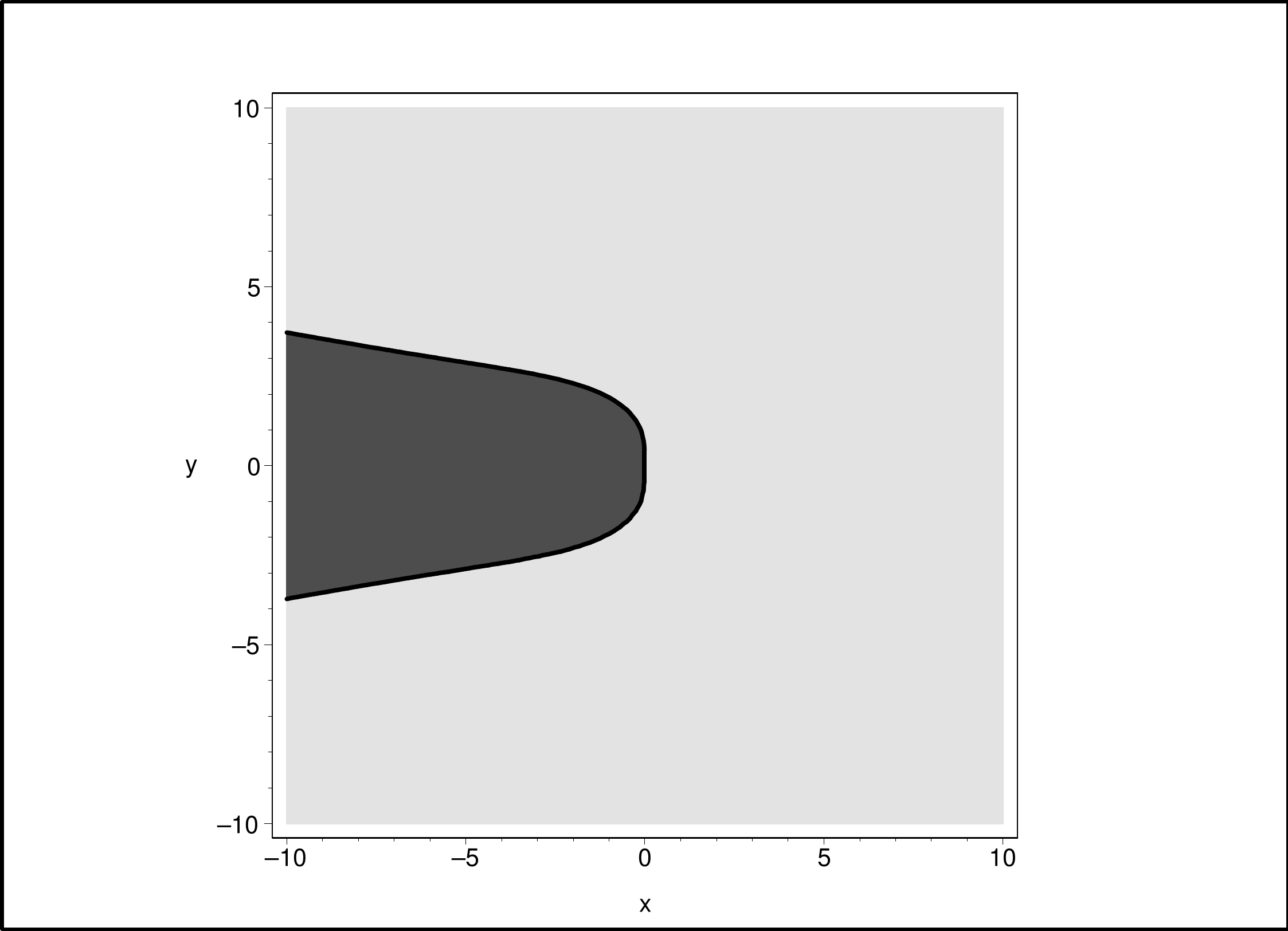}}}
\hspace{0.2cm}
  \subfigure[\sffamily SSP$3(3,3,2)$ (\ref{hig2})]{%
    \scalebox{0.25}{\includegraphics[width=12cm,angle=270]{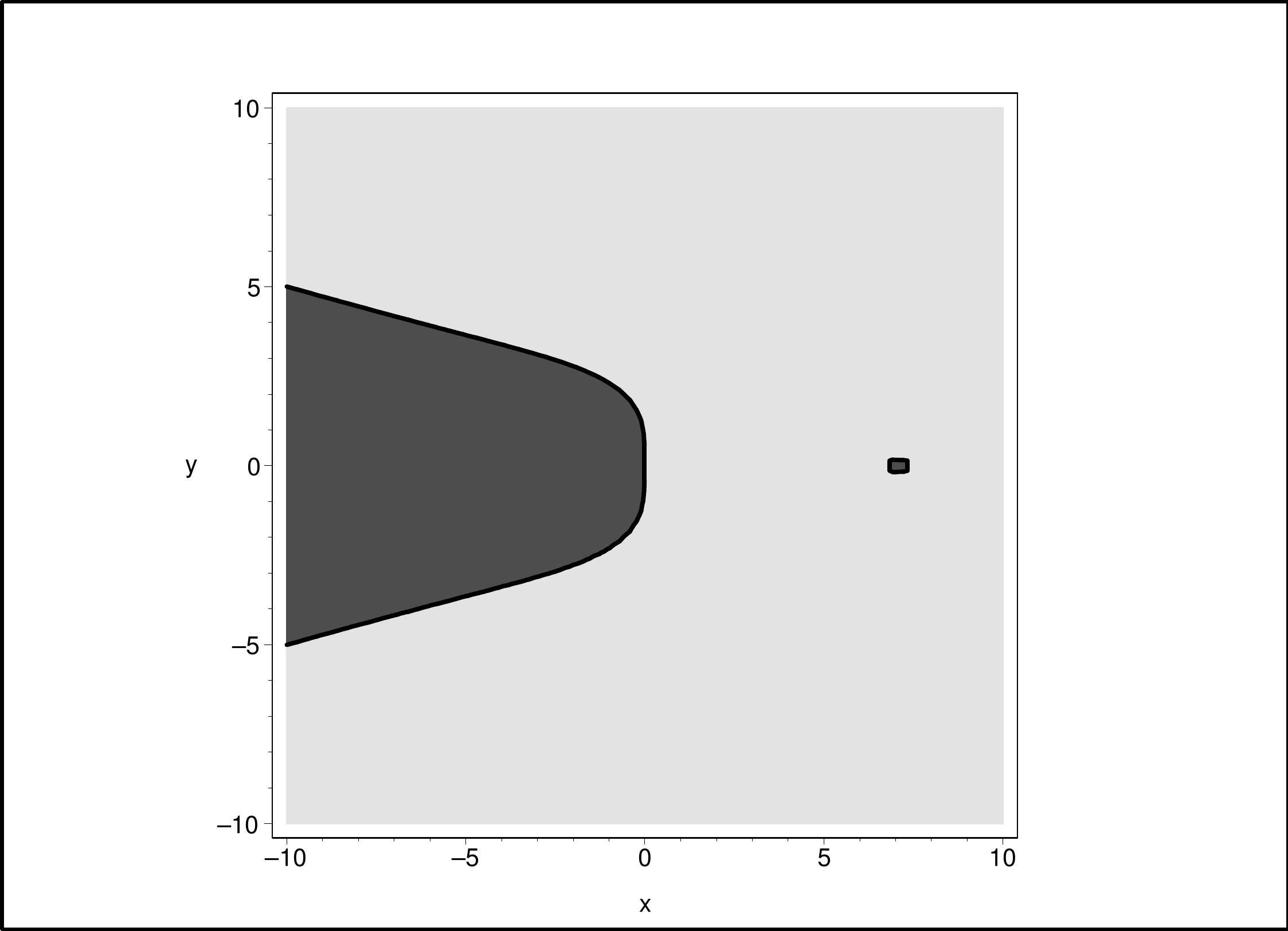}}}
\hspace{0.2cm}
  \subfigure[\sffamily SSP$3(3,3,3)$ (\ref{hig3})]{%
    \scalebox{0.25}{\includegraphics[width=12cm,angle=270]{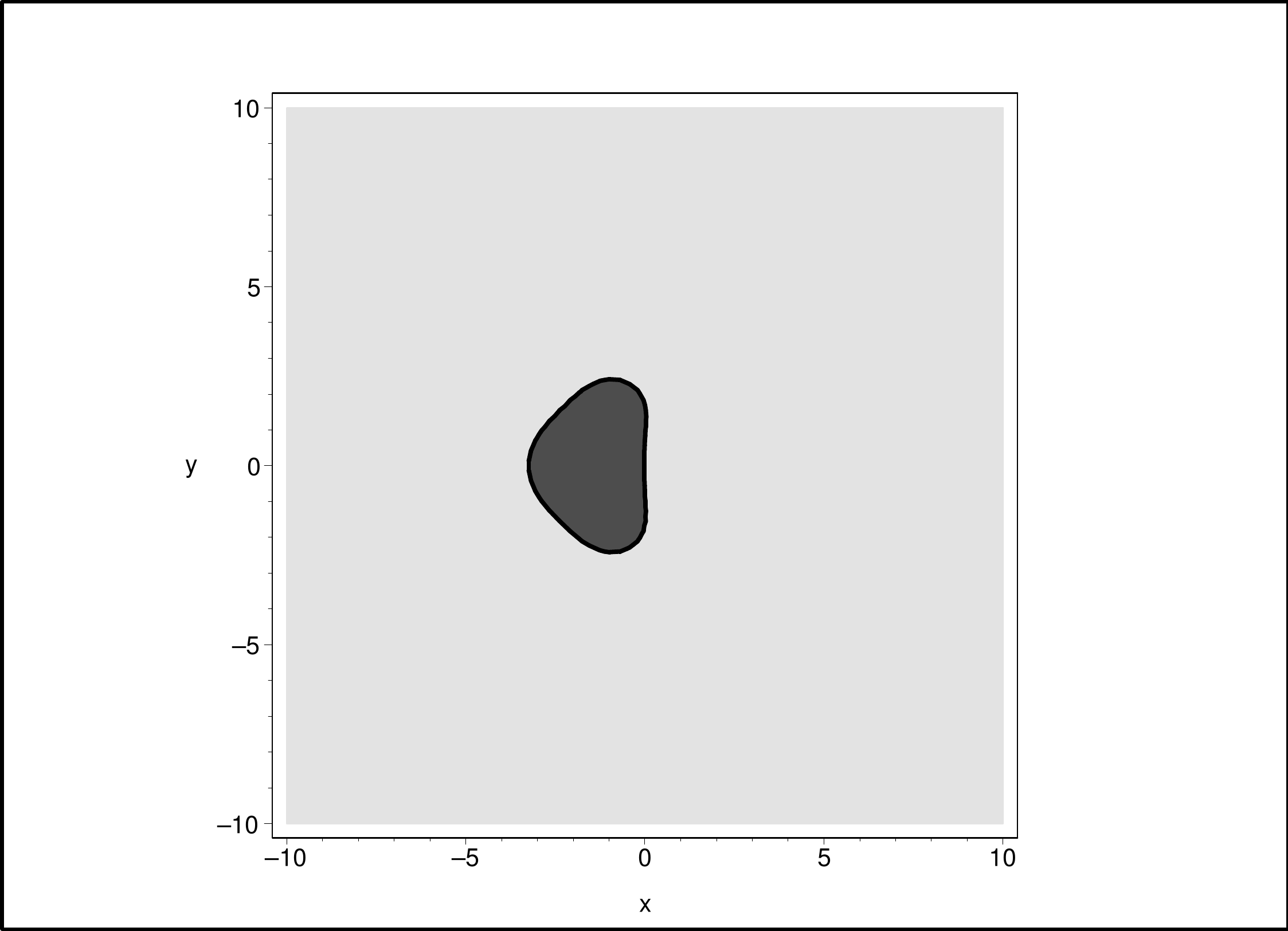}}}
    \\
  \subfigure[\sffamily $\tilde A$ in (\ref{hig1})]{
    \scalebox{0.25}{\includegraphics[width=12cm,angle=270]{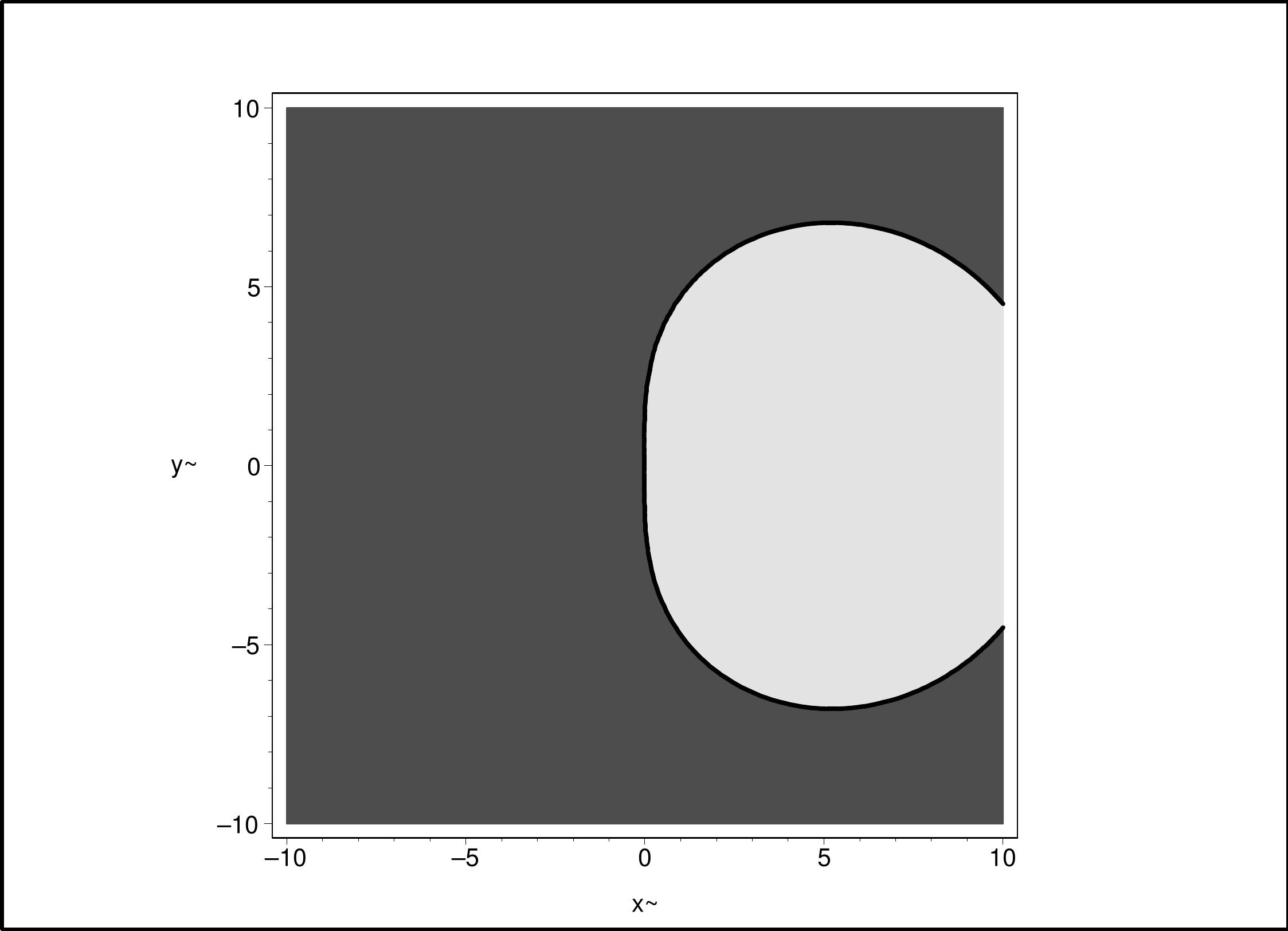}}}
\hspace{0.2cm}
  \subfigure[\sffamily $\tilde A$ in (\ref{hig2})]{%
    \scalebox{0.25}{\includegraphics[width=12cm,angle=270]{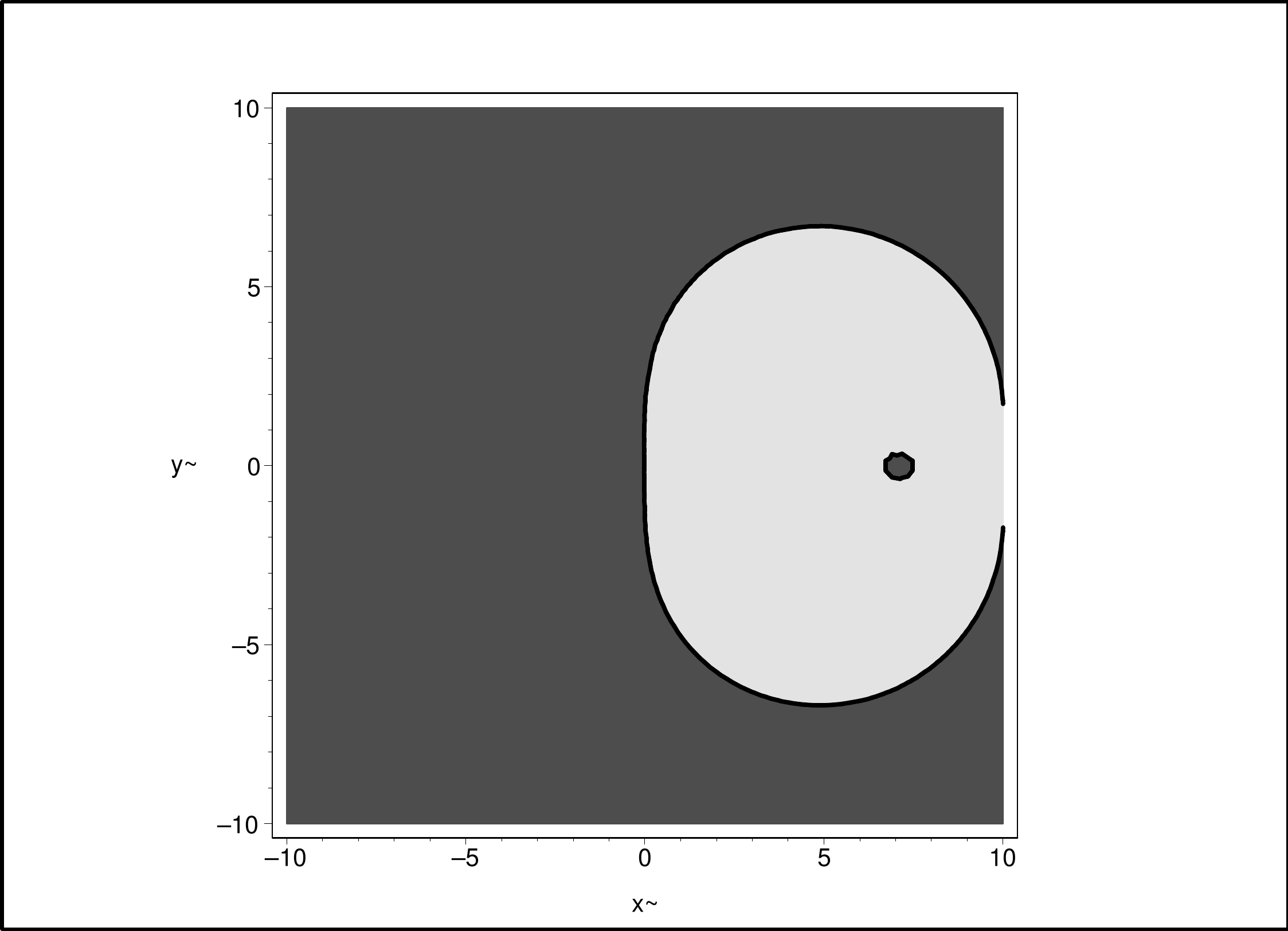}}}
\hspace{0.2cm}
  \subfigure[\sffamily $\tilde A$ in (\ref{hig3})]{%
    \scalebox{0.25}{\includegraphics[width=12cm,angle=270]{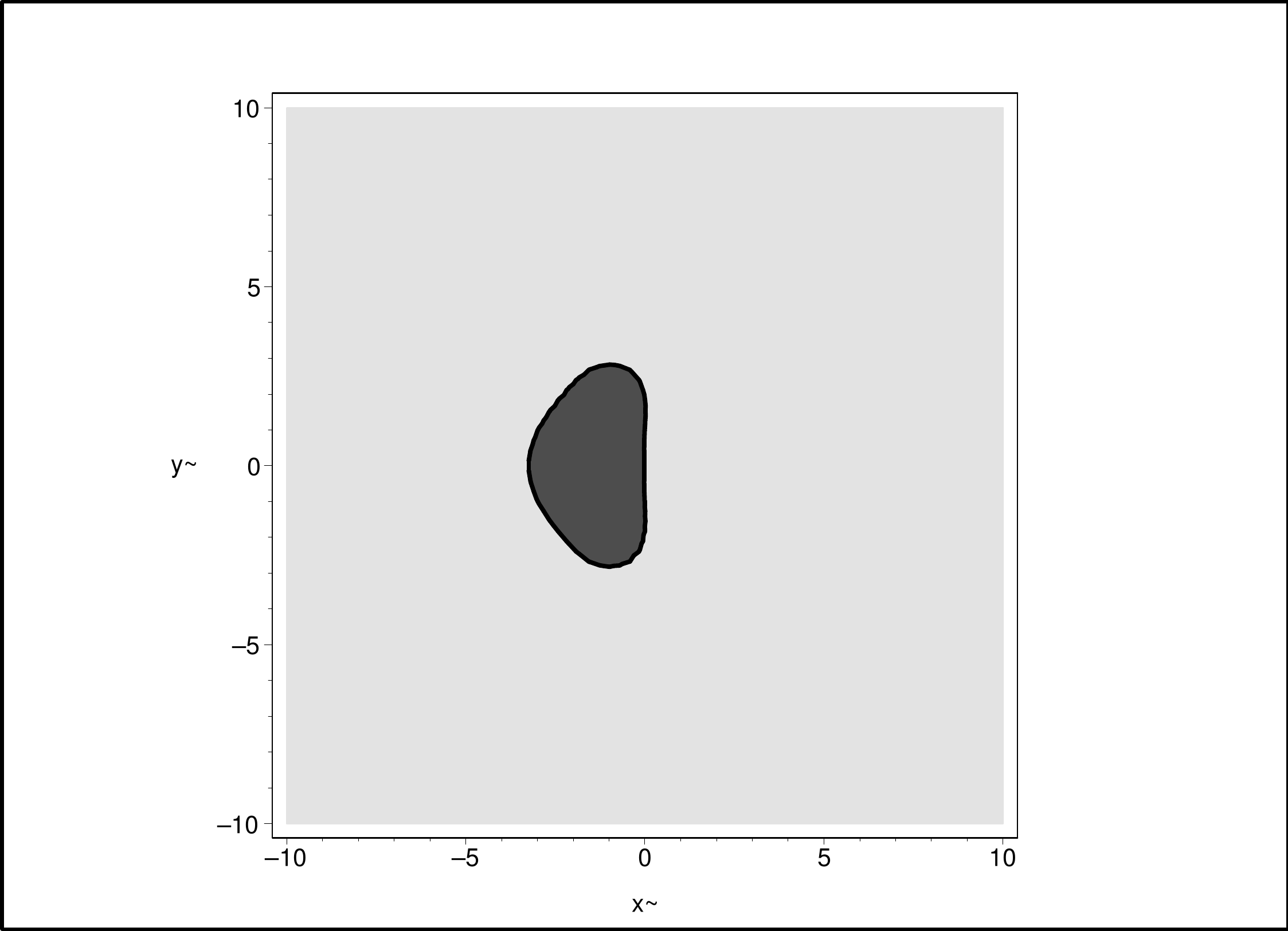}}}
  \subfigure[\sffamily SSPRK(2,2)]{%
    \scalebox{0.25}{\includegraphics[width=12cm,angle=270]{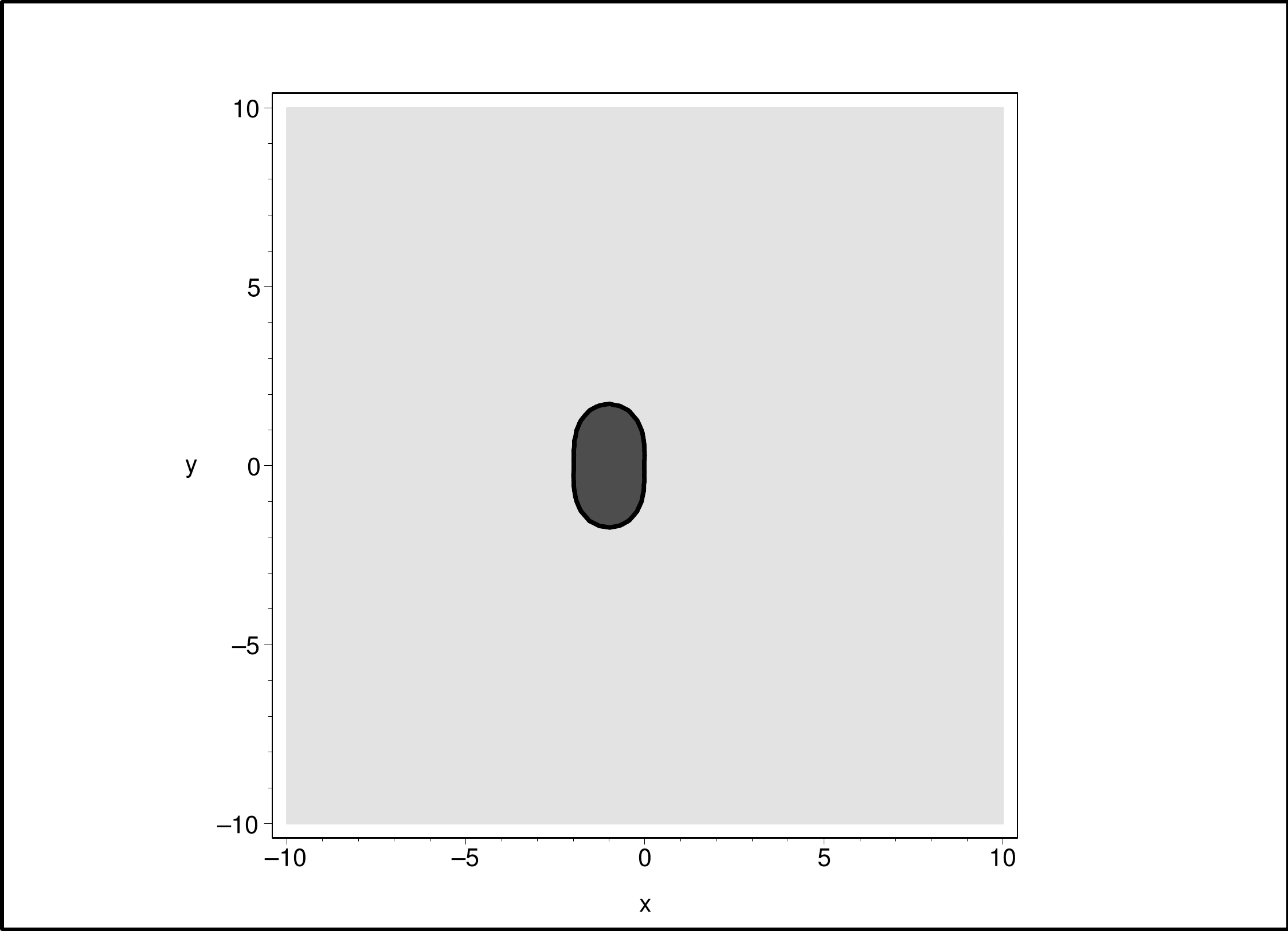}}}
\hspace{0.2cm}
  \subfigure[\sffamily SSPRK(3,2)]{%
    \scalebox{0.25}{\includegraphics[width=12cm,angle=270]{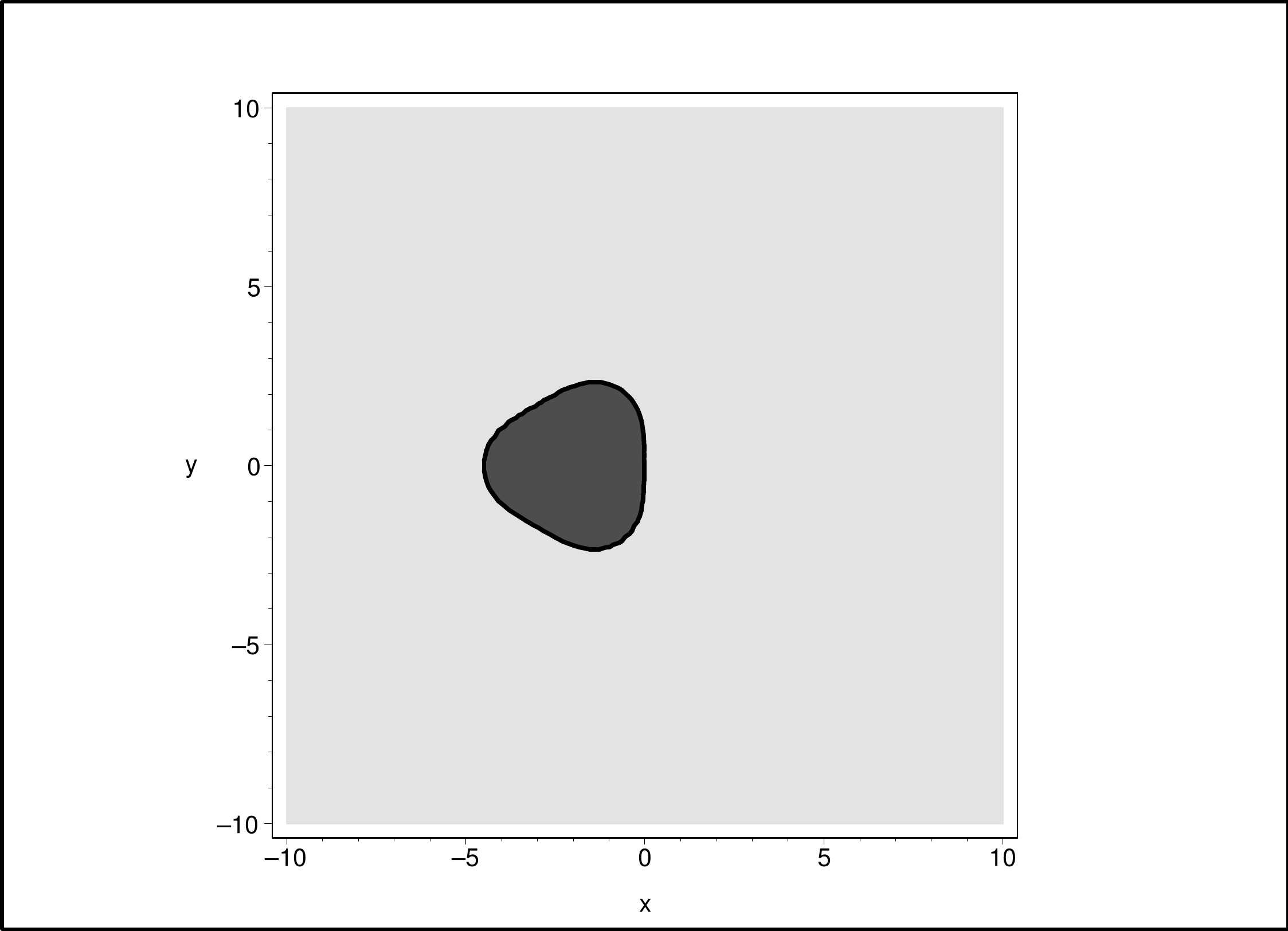}}}
\hspace{0.2cm}
  \subfigure[\sffamily SSPRK(3,3)]{%
    \scalebox{0.25}{\includegraphics[width=12cm,angle=270]{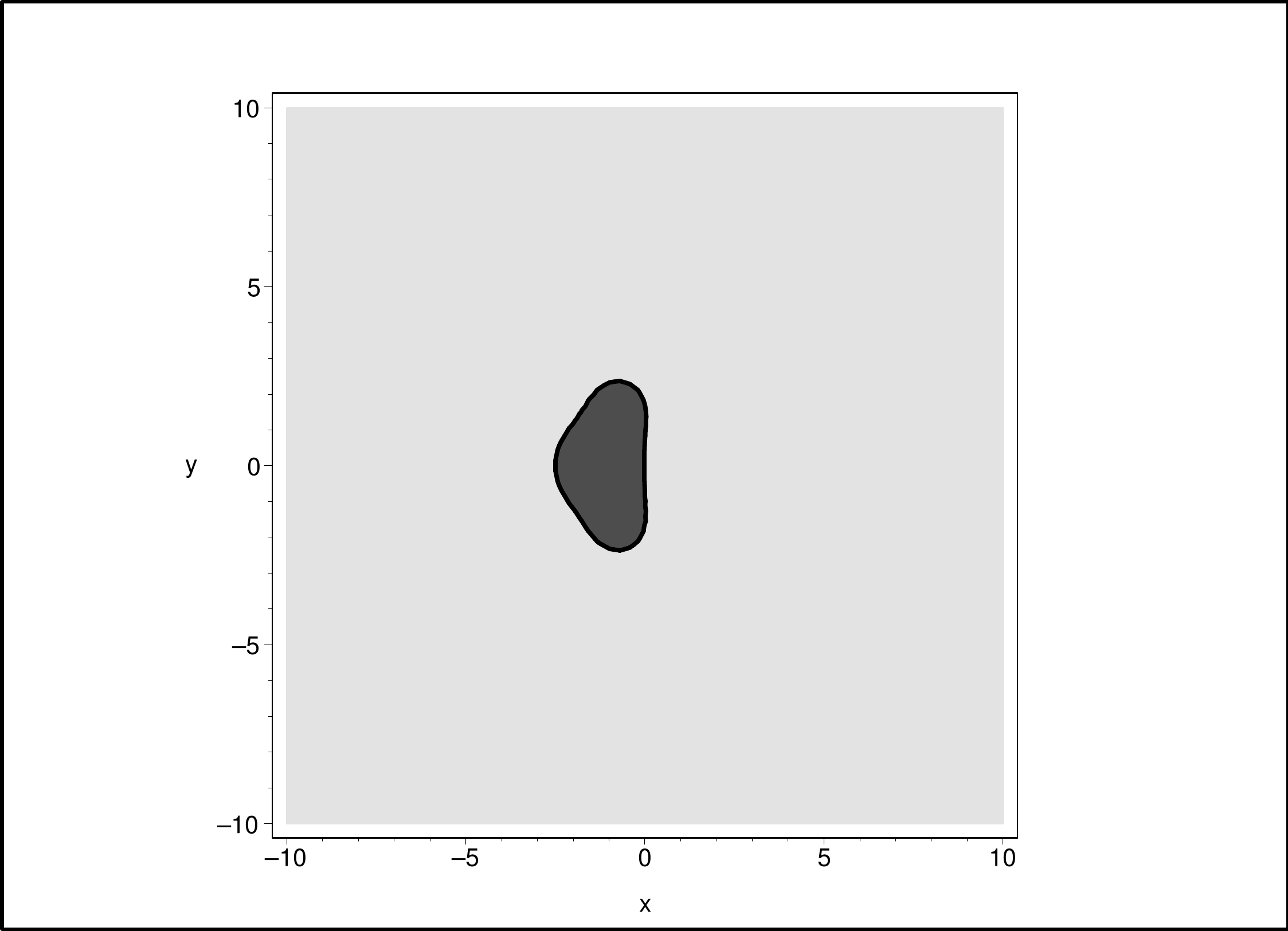}}}
\caption{Stability regions for IMEX methods (first row) and their implicit and explicit sub-parts (second and
third row).
\label{stabcompareimex}}
\end{figure}

The stability regions for the IMEX methods, their implicit sub-parts and
the explicit schemes are given for comparison in Figure~\ref{stabcompareimex}.
The left boundaries of the stability regions $z_\mathrm{left}$, the points where the amplification
factors from the dissipativity analysis become zero and their moduli
exceed 1, and the error constants $C$ for all the methods we have investigated
are summed up in Table~\ref{summary}.

We found that methods introduced in \cite{higueras05,higueras06,higueras09,parrus05}
excel over the classical explicit methods \cite{shu88a}.
It was found that among explicit schemes, only the explicit SSPRK(3,2) scheme first
proposed in \cite{kraaijevanger91} is competitive in situations where
explicit time integration can be expected to yield sufficient efficiency and accuracy.
Examples for such a scenario include simulations of solar granulation
at moderately high spatial resolution, where the time-step limitation
associated with diffusion is negligible (see \cite{Muthsam2010} for results on this
problem obtained with the ANTARES code and \cite{spruitetal90} for a review
on the underlying physics).

From Figure~\ref{stabcompareimex} and Table~\ref{summary}, clearly there is no single scheme
which features the most advantageous properties in all considered aspects.
However, we found in numerical experiments that the most efficient method
seems to be (\ref{hig1}) with the choice $\gamma=0.24$. This value deviates
from the optimal value for strong stability, but leads to a scheme with favourable
dissipativity, stability, and accuracy properties. Depending on the domain of stability required for
a given problem the value of $\gamma$ in (\ref{hig1}) may be optimized
such that it is sufficiently large for stability, but small enough to minimize
the error constant, while showing favourable dissipation properties (a strictly positive
amplification factor with modulus less than 1 for any wave number $k$ other
than zero).

We note that for numerical problems arising from a method of lines
approach to the equations of hydrodynamics, as discussed in this
paper, lower order methods usually have sufficient efficiency to be
competitive, since the spatial discretization limits the overall accuracy.
Hence, the best explicit scheme we have tested for this kind of application
is SSPRK(3,2), as it permits the largest CFL numbers among methods of this class 
at an affordable computational cost and with sufficient accuracy. By comparison, 
the classical methods of second and third order \cite{shu88a} offer the convenience 
of being usable together as an embedding formula. However, this approach
is more than twice as expensive as SSPRK(3,2), as can be seen
from a comparison with SSPRK(3,3) using Table~\ref{NumericalResults2}.

\begin{table}[h]
\begin{center}
\begin{tabular}{|l||r|r|r|r|r|r|r|}
\hline
\multicolumn{1}{|c||}{\raisebox{-1mm}{Method}} &
\multicolumn{1}{c|}{\raisebox{-1mm}{$z_\mathrm{left}$}} &
\multicolumn{1}{c|}{\raisebox{-1mm}{$g_\mathrm{4th}=0$}} &
\multicolumn{1}{c|}{\raisebox{-1mm}{$|g_\mathrm{4th}|=1$}} &
\multicolumn{1}{c|}{\raisebox{-1mm}{$C$}} \\ \hline
IMEX SSP2(2,2,2)                & $-\infty$ & $0.452^\ast$ & ---     & $5.17$ \\\hline
IMEX SSP2(2,2,2), $\gamma=0.24$ & $-50$     &          --- & 9.375   & $2.79$ \\\hline
IMEX SSP2(3,3,2)                & $-\infty$ & $0.455^\ast$ & ---     & $8.05$ \\\hline
IMEX SSP3(3,3,3)                & $-3.248$  & $0.348^\ast$ & $0.609$ & $11.6$ \\\hline
Forward Euler                   & $-2$      & $0.187^\ast$ & $0.375$ & $12.6$ \\\hline
SSPRK(2,2)                      & $-2$      & ---          & $0.375$ & $16.2$ \\\hline
SSPRK(3,2)                      & $-4.519$  & $0.672^\ast$ & $0.847$ & $6.40$\\\hline
SSPRK(3,3)                      & $-2.512$  & $0.299^\ast$ & $0.471$ & $22.8$ \\\hline
\hline
\end{tabular}
\caption{Summary of the analysis of SSP integrators. The asterisk in the third column
              indicates a change of sign at $\mu$  for $g_\mathrm{4th}(\pm\pi,\mu)$. 
              Other details are given in the text.  \mlabel{summary}}
\end{center}
\end{table}

We have also demonstrated that the larger time-steps achieved by SSP IMEX
methods reduce the accuracy of the solution during the diffusive phase of
the semi-convection simulations by an acceptably small amount. For the 
applications shown here, and indeed for a majority of astrophysical fluid dynamical
simulations, the accuracy is limited by spatial resolution (and thus eventually
by existing computational resources) while the time-steps are limited by 
stability. This makes IMEX methods attractive, since quite often the most
severe limitations stem from stiff terms representing diffusion processes
(for restrictions due to sound waves other operator splitting based methods
are existing). However, as we have shown by a comparison with results
from non-SSP methods, it is important that the IMEX methods are
strong--stability--preserving to maximize stable time-steps no matter
whether the constraint is due to the implicitly integrated terms (diffusion) or
the explicitly integrated ones (advection). From that point of view
SSP IMEX methods with a non-trivial region of absolute monotonicity as
defined by Theorem~\ref{am1} are the most robust methods, because
they allow achieving optimally large time-steps also at low resolution. 
We note here that while the region of absolute monotonicity has to be 
observed with respect to the explicitly integrated advection operator of 
the dynamical equations, stable time integration can be performed with step 
sizes falling outside it, if the restriction is due to the implicitly integrated diffusion terms.
Thus, the class of optimal integrators for this kind of problem is probably
larger than that of the SSP IMEX methods. However, none of the other time 
integration methods could significantly outperform them with respect
to the time-steps achievable, and we have always found at least one
case, where competing methods fell substantially short or even failed.

There is potential to further optimize the implementation of IMEX methods.
The additional computational effort due to the implicit subpart is compensated
for by accuracy and stability, but could be reduced in the future by replacing
the solver for the linear equations associated with the arising generalized
Poisson problem by a multigrid solver. In the Boussinesq approximation,
additional solution of a Helmholtz equation is necessary instead. This widens
the choice of fast solvers for the system of linear equations introduced
through implicit time integration. The benefits expected from faster solvers would
allow taking full advantage of the potential of method (\ref{hig2}) that is implied
by the large time steps reported in Table~\ref{NumericalResults1}.
Such an improvement would likewise be useful for the present problem
to minimize the overhead by any of the implicit schemes in the regime where the
time-step is limited by $\tau_{\mathrm{fluid}}$ rather than $\tau_T$.

\begin{ack}
We gratefully acknowledge the help of H.~Grimm--Strele who performed the 2D solar
convection simulations with ANTARES to test the stability properties of Heun's non-SSP explicit 
third order Runge--Kutta scheme. We also thank H.J.~Muthsam for useful discussions
on the ANTARES code as well as J.~Ballot and H.~Grimm--Strele for help in a generic
implementation of the SSPRK(3,2) method which can now be used by all modules of the code.
\end{ack}

\section*{Appendix}

\subsubsection*{The Explicit SSP Scheme of IMEX SSP2(3,3,2)}

We also give the corresponding results for the explicit SSPRK(3,2) scheme $A$ from (\ref{hig2}),
since in Section~\ref{num} we show it to excel in its practical value over the classical explicit
SSP Runge--Kutta schemes schemes \cite{shu88a}. This scheme was first published in \cite{kraaijevanger91} and later
declared the optimal second order scheme with three stages in \cite{spiruu02}
as well as in \cite{ruuspi04}, and independently also in \cite{gottliebsquared03}.

The stability function is
$$R_A(z)= 1 + z + \frac{z^2}{2} + \frac{z^3}{12}.$$
The stability region where $|R(z)|<1$ occupies a bounded region in the
negative half-plane, and is tangent to the imaginary axis, see
Figure~\ref{fig2expl}. Note that $\lim_{\Re(z)\to -\infty}|R(z)| =\infty$.
\begin{figure}[h]
\begin{center}
\includegraphics[width=4.7cm,angle=270]{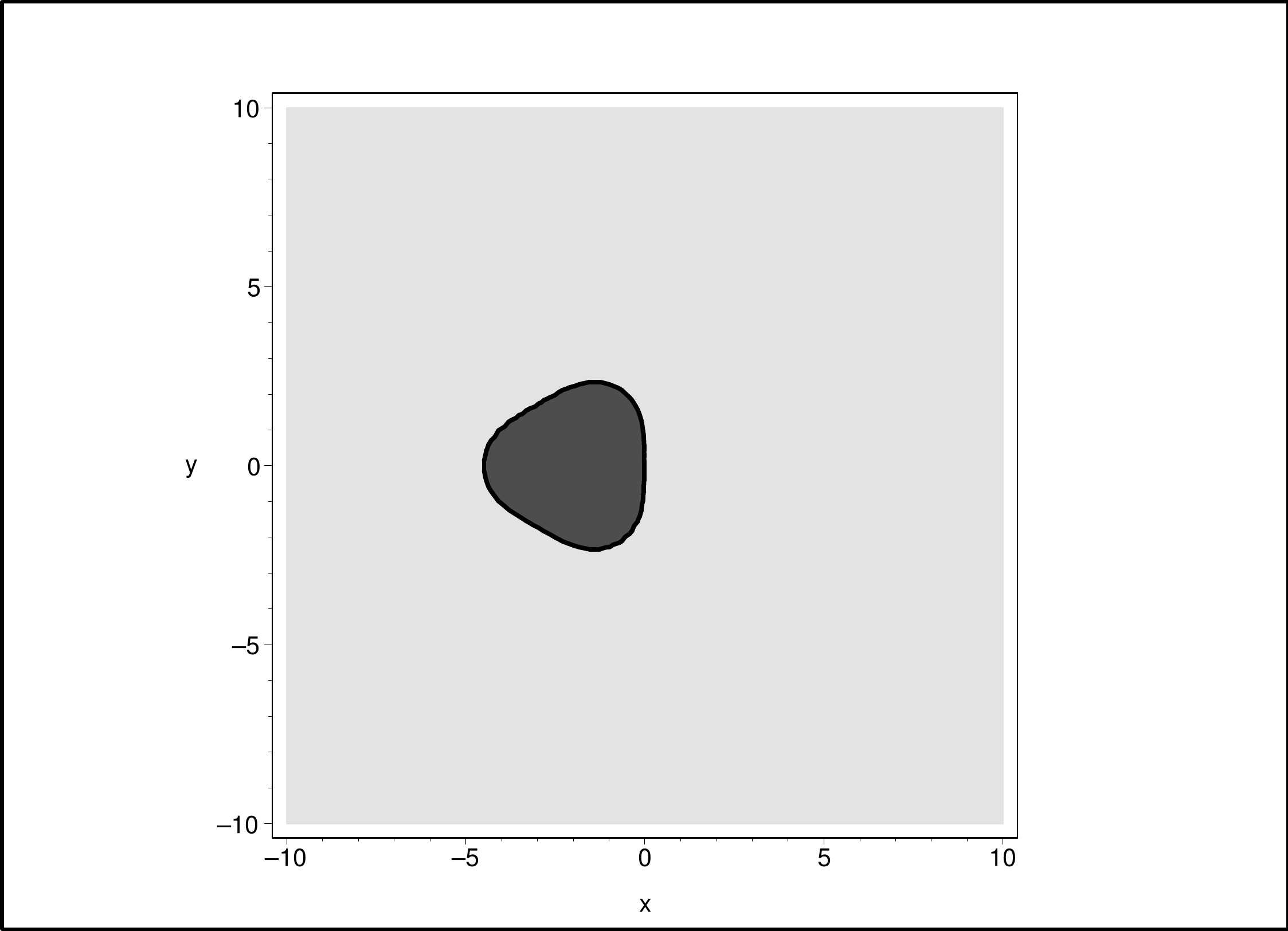}
\caption{Stability region of method SSPRK(3,2) (scheme $A$ from~(\ref{hig2})).\mlabel{fig2expl}}
\end{center}
\end{figure}

The point where the stability region intersects the negative real
half-line is located at $x\approx-4.519$.

The dissipativity analysis for $A$ yields the
amplification factors for the standard
three-point space discretization (\ref{3pt}) and the fourth order stencil (\ref{4th}).
These amplification factors are evaluated at the points $\theta\in\{0,\frac\pi4,\frac\pi2,\pi\}$
in Tables~\ref{evalsfimex2expl3pt} and \ref{evalsfimex2expltrial}, respectively.
\begin{table}
\begin{center}
\begin{tabular}{||r|c||}
\hline
\multicolumn{1}{||c|}{$\theta$} & \multicolumn{1}{c||}{$g(\theta,\mu)$} \\
\hline\hline
$0$ &  $1$ \\ \hline
$\frac{\pi}{4}$ & $1+\mu\,\sqrt {2}-2\,\mu+3\,{\mu}^{2}-2\,{\mu}^{2}\sqrt {2}+7/6\,{\mu}^
{3}\sqrt {2}-5/3\,{\mu}^{3}$  \\ \hline
$\frac{\pi}{2}$ & $1-2\,\mu+2\,{\mu}^{2}-2/3\,{\mu}^{3}$  \\ \hline
$\pi$ & $1-4\,\mu+8\,{\mu}^{2}-16/3\,{\mu}^{3}$ \\ \hline
\end{tabular}
\caption{Values of $g(\theta,\mu)$ for some $\theta$, SSPRK(3,2) scheme (explicit scheme $A$ from (\ref{hig2})),
three point space discretization (\ref{3pt}).
\mlabel{evalsfimex2expl3pt}}
\end{center}
\end{table}

\begin{table}
\begin{center}
\begin{tabular}{||r|c||}
\hline
\multicolumn{1}{||c|}{$\theta$} & \multicolumn{1}{c||}{$g(\theta,\mu)$} \\
\hline\hline
$0$ &  $1$ \\ \hline
$\frac{\pi}{4}$ & $1-5/2\,\mu+4/3\,\mu\,\sqrt {2}+{\frac {353}{72}}\,{\mu}^{2}-10/3\,{\mu
}^{2}\sqrt {2}-{\frac {1015}{288}}\,{\mu}^{3}+{\frac {803}{324}}\,{\mu
}^{3}\sqrt {2}$  \\ \hline
$\frac{\pi}{2}$ & $1-7/3\,\mu+{\frac {49}{18}}\,{\mu}^{2}-{\frac {343}{324}}\,{\mu}^{3}$  \\ \hline
$\pi$ & $1-16/3\,\mu+{\frac {128}{9}}\,{\mu}^{2}-{\frac {1024}{81}}\,{\mu}^{3}$ \\ \hline
\end{tabular}
\caption{Values of $g(\theta,\mu)$ for some $\theta$, SSPRK(3,2) (explicit scheme $A$ from (\ref{hig2})), fourth order space discretization
(\ref{4th}).
\mlabel{evalsfimex2expltrial}}
\end{center}
\end{table}

For the three-point space discretization (\ref{3pt}), the first positive zero of $g(\pi,\mu)$ is $\approx 0.8968$,
where the function changes its sign, and $g(\pi,\mu)=-1$ for $\mu\approx 1.129$.
The first positive zero of $g(\pi,\mu)$ is $\approx 0.6726$ for the fourth-order space discretization (\ref{4th}),
where the function changes its sign, and $g(\pi,\mu)=-1$ at $\mu\approx 0.8474$.

\bibliographystyle{elsart-num-sort}

\bibliography{books,num}

\end{document}